\documentclass[final]{siamltex}

\usepackage{amsmath}
\usepackage{mathrsfs}
\usepackage{amsfonts}
\usepackage{amssymb}
\usepackage{graphicx}
\usepackage{psfrag}
\usepackage{placeins}
\usepackage{epsfig}
\usepackage{amsopn}
\usepackage{amstext}
\usepackage{amssymb}
\usepackage{amscd}
\usepackage{latexsym}
\allowdisplaybreaks
\usepackage{geometry}
\usepackage{multirow}
\usepackage{diagbox}
\usepackage{color}

\newtheorem{rem}[theorem]{Remark}

\newtheorem{thm}{Theorem}[section]

\newcommand{\bbm}{\begin{bmatrix}}
	\newcommand{\ebm}{\end{bmatrix}}

\renewcommand{\theequation}{\arabic{section}.\arabic{equation}}
\usepackage{subfigure}
\usepackage{graphicx}
\usepackage{float}
\usepackage[numbers,sort&compress]{natbib}

\begin{document}

\title{Ensemble Domain Decomposition Algorithm for the Fully-mixed Random Stokes-Darcy Model with the Beavers-Joseph Interface Conditions}
\author{
Feng Shi \thanks{%
College of Science, Harbin Institute of Technology, Shenzhen,
P.R. China. \texttt{shi.feng@hit.edu.cn}. Partially supported by Foundation Research Project of Shenzhen (Grant No. GXWD20201230155427003-20200822102539001).}
\and Yizhong Sun \thanks{%
School of Mathematical Sciences, East China Normal University,
Shanghai, P.R. China. \texttt{bill950204@126.com}.}
\and Haibiao Zheng \thanks{%
	School of Mathematical Sciences, East China Normal University,
	Shanghai Key Laboratory of Pure Mathematics and Mathematical
	Practice, Key Laboratory of Advanced Theory and Application in Statistics and Data Science (East China Normal University), Shanghai, P.R. China. \texttt{hbzheng@math.ecnu.edu.cn}.
	{ Partially supported by  
	Science and Technology
	Commission of Shanghai Municipality (Grant Nos. 22JC1400900, 21JC1402500, 22DZ2229014) and NSF of China (Grant No. 11971174)}.} }

 \maketitle

\begin{abstract}
In this paper, an efficient ensemble domain decomposition algorithm is proposed for fast solving the fully-mixed random Stokes-Darcy model with the physically realistic Beavers–Joseph (BJ) interface conditions. We utilize the Monte Carlo method for the coupled model with random inputs to derive some deterministic Stokes-Darcy numerical models and use the idea of the ensemble to realize the fast computation of multiple problems. One remarkable feature of the algorithm is that multiple linear systems share a common coefficient matrix in each deterministic numerical model, which significantly reduces the computational cost and achieves comparable accuracy with the traditional methods. Moreover, by domain decomposition, we can decouple the Stokes–Darcy system into two smaller sub-physics problems naturally. Both mesh-dependent and mesh-independent convergence rates of the algorithm are rigorously derived by choosing suitable Robin parameters.  Optimized Robin parameters are derived and analyzed to accelerate the convergence of the proposed algorithm. Especially, for small hydraulic conductivity in practice, the almost optimal geometric convergence can be obtained by finite element discretization. Finally, two groups of numerical experiments are conducted to validate the exclusive features of the proposed algorithm.
\end{abstract}

\begin{keywords}
Random Stokes-Darcy Model, Beavers-Joseph Interface Conditions, Ensemble Domain Decomposition, Optimized Schwarz Method, Geometric Convergence.
\end{keywords}
\begin{AMS}
65M55, 65M60
\end{AMS}

\section{Introduction}
Multi-domain, multi-physics coupled problems are significant in many natural and industrial applications, such as groundwater fluid flow in the karst aquifer,  petroleum extraction, industrial filtration, blood flow motion in the arteries, and so on.  A great deal of the mathematically and physically models was constructed for the coupling of surface flows (in free fluid flow) and groundwater flows (in porous media), including the Stokes-Darcy model \cite{DMQ02, LMLayton2003,Badea, Cao10}, Stokes-Darcy-transport/heat model \cite{Zhang20}, dual-porosity-Stokes model \cite{HouJY16}, to name just a few. Typically, the famous Stokes-Darcy model can describe the coupling of one free fluid flow with a single porosity medium flow in two subdomains separated by an interface.

Inspired by the decoupled idea for the Stokes–Darcy model, a natural way is the domain decomposition method (DDM) \cite{Discacciati07, Chen, Cao11, He15, Vassilev1, Sun2021}, because it can decouple the multi-domain, multi-physics problems naturally under suitable interface boundary conditions, and there exits many well-established off-the-shelf and efficient solvers for each decoupled subproblem. Based on the characteristics of easy-to-operation, high precision, and convenient parallel computing, DDM has received extensive attention and applications undoubtedly. In \cite{Cao11}, two Robin-Robin domain decomposition methods for the steady-state Stokes-Darcy system with Beaver-Joseph (BJ) interface conditions are proposed, based on the pioneering work for Beavers-Joseph-Saffman (BJS) interface conditions \cite{Chen}. The authors demonstrated the convergence of DDMs and proved the geometric convergence rate with some suitable choices of Robin parameters, under the assumption that the exchange coefficient $\alpha$ in the BJ interface conditions is sufficiently small.
Sun et al. \cite{Sun2021} first utilized DDM to study the fully-mixed Stokes-Darcy coupled problem with BJS interface conditions, and rigorously proved both mesh-dependent and mesh-independent convergence rates with the suitable choice of Robin parameters based on the introduced modified weak formulation. 

From the early research,  we can conclude that the convergence of DDMs was closely dependent on the selection of Robin parameters. More precisely, the Stokes and Darcy Robin parameters are denoted as $\delta_S$ and $\delta_D$ respectively. Most studies only focused on the case of $\delta_S=\delta_D$, so as to obtain the convergence rate with the dependence of the mesh size. Then some researchers found that the DDMs had geometric convergence for $\delta_S<\delta_D$ \cite{Chen, Cao11, Sun2021}. This choice can save a lot of computing costs so that the DDMs have better potential applications. 
Recently, optimized Schwarz methods have been proposed for the Stokes-Darcy model in \cite{Discacciati, Discacciati18, Gander, Liu22}, which can seek the optimal Robin parameters to speed up convergence. In \cite{Discacciati}, the authors utilized Fourier techniques to explicitly characterize the convergence factor, and then derive the optimal parameters by several practical strategies.
By now the results worked well with moderate hydraulic conductivity $K$. However, for the practical coefficients, such as $10^{-8}\leq K \leq 10^{-2}$, the DDMs seem difficult to converge for $\delta_S\leq \delta_D$. The numerical experiments in \cite{Cao11} could observe the convergence for the practical coefficients while $\delta_S > \delta_D$, but no complete theoretical analysis  supported such observations.  Motivated by the practical applications and analysis difficulty, Liu et al. \cite{Liu21} used the finite element discretization to further improve the convergence results of the DDMs proposed by Cao et al. \cite{Cao11}. They obtained the almost optimal geometric convergence rate in the case of $\delta_S>\delta_D$, particularly for small viscosity and hydraulic conductivity in practice.

It should be noted that in many engineering and geological applications, due to the complexity of the porous media domain and the limitation of measuring instruments, it is not feasible to obtain exact hydraulic conductivity values, whose natural randomness commonly occurs at small scales. So such uncertainties should be taken into account in the numerical simulation. An interesting skill is that the uncertain parameter of interest can be regarded as a random function that is determined by a basic random field with a specified covariance structure (usually determined by experiment). 
In this case, we need to deal with random partial differential equations (PDEs) \cite{ WangZ18, WangZ19}. 

The most popular approach to solving the random PDEs with random inputs is the Monte Carlo method, which transforms the random PDEs into some traditional PDEs. However, in order to get useful statistical information from the solutions, such a method inevitably requires the computation of a large number of realizations with a very slow convergence rate. It is still prohibitively time consuming to use the existing deterministic solvers repetitively.
To overcome such computational challenges, a class of ensemble methods was developed recently \cite{WangZ18, WangZ19, Feng, JiangN14, Max17, Mohebujjaman, Max19, JiangN19, JiangN21}. Corresponding to different physical parameters or body forces, such ensemble algorithms maintained one important feature that several linear systems can share a common coefficient matrix so that both required computational time and storage can be reduced. In \cite{JiangN19}, Jiang et al. first used the ensemble algorithm to solve the non-stationary random Stokes-Darcy equations with BJS interface conditions. They utilized the information of the previous time step to construct an explicit-implicit decoupling ensemble algorithm. It's worth noting that one can use the Gronwall lemma in theoretical analysis for the unsteady explicit-implicit decoupling ensemble algorithm. But for the steady-state system, such a technique fails.  One interesting work in this paper is using the iterative technique to decouple the system for solving the steady-state Stokes-Darcy model.  More importantly, we establish Lemma \ref{abc} to carry out the theoretical analysis instead of such inequality of Gronwall type.

In the present work, we follow the idea of \cite{JiangN19} to propose an efficient ensemble domain decomposition algorithm for  the steady-state fully-mixed random Stokes-Darcy model. Meanwhile, we extend the Robin-type DDMs to solve the Stokes-Darcy model from BJS interface conditions in \cite{Sun2021} to the BJ interface conditions, which have much physical significance but will also bring out several analytical challenges. For the convergence of the proposed parallel ensemble DDM, we obtained the following important results:
\begin{itemize}
    \item $\delta_S<\delta_D$: We prove that our algorithm has a geometric convergence rate in the continuous system. More importantly, through ingenious theoretical analysis, we find that the exchange coefficient $\alpha$ in the BJ interface conditions does not need sufficiently small while the existing studies required $\alpha$ to be small enough, see Remark 4.2.
    \item $\delta_S=\delta_D$: In this case, the proof for convergence is similar to the discussion of $\delta_S<\delta_D$, and the results of this part are consistent with existing studies \cite{Cao10, Cao11}. 
    \item $\delta_S>\delta_D$: We utilize finite element discretization to obtain the almost optimal geometric convergence rate for small hydraulic conductivity in practice. The demonstration is greatly inspired by \cite{Liu21}, but the proof process is much clear. Moreover, the fully-mixed Stokes-Darcy model is more natural than the model in \cite{Liu21}, because engineering applications usually take more concern on the flux or Darcy velocity in the porous media domain. The rigorous analysis illustrates a general guideline of selecting the Robin parameters to achieve the geometric convergence rate.
    \item Optimized Robin parameters: Inspired by the optimized Schwarz methods in \cite{Discacciati, Discacciati18, Gander, Liu22}, we investigate optimized approaches for the fully-mixed random Stokes-Darcy model, and obtain the optimized Robin parameters that can greatly improve the convergence of the proposed ensemble DDM. More importantly, the existing studies have obtained a homogeneous ordinary differential equation through Fourier transformation. Due to the characteristics of the ensemble method, one nontrivial equation with iteration appears in the present work. Then, several analysis skills are developed to overcome the above challenges for deriving the convergence factor.
\end{itemize}

The rest of the paper is organized as follows. The fully-mixed random Stokes-Darcy model is described in Section 2. In Section 3, we propose an efficient ensemble DDM with three Robin-type condition sets. Then, the convergence of ensemble DDM is presented for the continuous system when the Robin parameters satisfy $\delta_S\leq \delta_D$ in Section 4. In particular, the geometric convergence rate of the proposed algorithm is derived in the case of $\delta_S<\delta_D$.  Furthermore, the optimized Robin parameters are analyzed and  derived  to speed up the convergence of the ensemble DDM in Section 5. In Section 6, the finite element approximations are discussed. We prove the almost optimal geometric convergence in the case of $\delta_S>\delta_D$ for small hydraulic conductivity. Finally, two numerical tests are presented to illustrate the exclusive features of the proposed DDM in Section 7. 

\section{Fully-mixed Random Stokes-Darcy Model with BJ Interface Conditions}
Consider the bounded domain ${\Omega }$, which is the union of two bounded subdomains, namely $\Omega _{S},\Omega_{D}\subset R^{d}(d=2~\mathrm{or}~3)$ with an interface $\Gamma$, and $\Omega _{S}\cap \Omega _{D}={\emptyset }$, $\ \overline{\Omega }_{S}\cap\overline{\Omega }_{D}=\Gamma $. These two subdomains $\Omega _{S}$ and $\Omega_{D}$ are usually referred as the free-flow and porous media domains. Denote by $\mathbf{n}_{S}$ and $\mathbf{n}_{D}$ the unit outward normal vectors on $\partial \Omega _{S}$ and $\partial \Omega _{D}$, respectively. It is worth to mention that $\mathbf{n}_{S}=-\mathbf{n}_{D}~\mathrm{on}~\Gamma $. Besides, the unit tangential vectors on the interface $\Gamma $ are represented by $\mathbf{\tau }_{i},i=1,\cdots ,d-1$. And note that $\Gamma_{S}=\partial \Omega _{S}\setminus \Gamma ,\ \Gamma _{D}=\partial \Omega_{D}\setminus \Gamma $, see Fig. 2.1 for a sketch.
\begin{figure}[htbp]
\centering
\includegraphics[width=70mm,height=40mm]{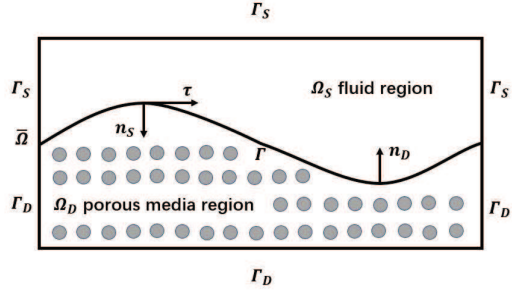}
\caption{  $\overline{\Omega}$ consisting of the fluid region
$\Omega_S$ and the porous media region $\Omega_D$ separated by the
interface $\Gamma$.}
\end{figure}

In the fluid region $\Omega_{S} $, the fluid velocity $\mathbf{u}_{S}$ and kinematic pressure $p_{S}$
are assumed to satisfy the Stokes equations:
\begin{eqnarray}
 -\nabla \cdot \mathbb{T}(\mathbf{u}_{S},p_S)&=&\mathbf{f}_{S}\ \ \ \ \mathrm{in}~\Omega _{S},
\label{Stokes1} \\
\nabla \cdot \mathbf{u}_{S} &=&0\ \ \ \hspace{0.3mm} \ \ \mathrm{in}~\Omega _{S},\label{Stokes2}
\end{eqnarray}
where $\mathbb{T}(\mathbf{u}_{S},p_S)=-p_s\mathbb{I}+2\nu\mathbb{D}(\mathbf{u}_{S})$ indicates the stress tensor, herein $\mathbb{D}( \mathbf{u}_{S})=\frac{1}{2}(\nabla \mathbf{u}_{S}+(\nabla \mathbf{u}_{S})^{T})$ denotes the deformation tensor, and $\nu $ represents the kinematic viscosity of the fluid flow. Besides $\mathbf{f}_S$ is the given external body force.

The porous media flow in $\Omega_{D} $
is governed by the following mixed Darcy equations for the
fluid velocity $\mathbf{u}_{D}$ and the
piezometric head $\phi_{D}$:
\begin{eqnarray}
\mathbf{u}_{D}&=&-\mathbb{K}(\mathbf{x}) \nabla \phi_{D}\ \ \ \ \  \
\mathrm{in}~\Omega_{D} ,  \label{Darcy1} \\
 \nabla \cdot \mathbf{u}_{D} &=& {f}_{D}\ \ \ \ \ \ \ \ \ \hspace{9.8mm} %
\mathrm{in}~\Omega _{D}, \label{Darcy2}
\end{eqnarray}
where $\mathbb{K}(x) $ is the hydraulic conductivity tensor and is physically impossible to determine its parameter values. In addition, ${f}_{D} $ denotes a sink/source term and $\int_{\Omega_D} f_D=0$. The piezometric head $\phi_{D}$ is defined by $\phi_D=z+\frac{p_D}{\rho g}$, where $z$ is the height,  $p_D$ indicates the dynamic pressure, $g$ denotes the gravitational acceleration, and $\rho$ represents density.

We assume that the fluid velocity $\mathbf{u}_{S} $ and the porous media velocity $\mathbf{u}_{D} $ satisfy homogeneous Dirichlet boundary conditions on the impermeable boundary except for the interface $\Gamma$, i.e., no-slip conditions $%
\mathbf{u}_{S}=0~~\mathrm{on}~\Gamma _{S} $, and $
\mathbf{u}_{D}\cdot\mathbf{n}_D=0~~\mathrm{on}~\Gamma _{D} $.
On the interface $\Gamma $, some coupling interface conditions are essential, including conservation of mass, the balance of forces, and tangential conditions on the fluid region's velocity. Also, we impose the original and more physically realistic Beavers-Joseph (BJ) as the tangential conditions on the interface $\Gamma $, see \cite{Beavers, Jones}. In this paper, the interface coupling conditions are assumed as follows:
\begin{eqnarray}
\mathbf{u}_{S}\cdot \mathbf{n}_{S}+\mathbf{u}_{D}\cdot
\mathbf{n}_{D}
&=&0\ \ \ \ \ \ \ \ \ \ \ \ \ \ \hspace{7.5mm} \mathrm{on}~\Gamma ,  \label{interface1} \\
-\mathbf{n}_{S}\cdot (\mathbb{T}(\mathbf{u}_{S},p_S)\cdot \mathbf{n}_{S})
&=&  g(\phi_{D}-z) \ \ \ \ \ \ \  \ \  \mathrm{on}~\Gamma, \label{interface2}\\
-\tau_{i}\cdot(\mathbb{T}(\mathbf{u}_{S},p_S)\cdot \mathbf{n}_{S})
&=&\frac{\alpha}{\sqrt{\tau_{i}\cdot\mathbb{K}\tau_{i}}} \tau _{i} \cdot (\mathbf{u}_{S}-\mathbf{u}_{D})
\hspace{3mm}  1\leq i \leq d-1\hspace{3mm} \mathrm{on}~\Gamma,\label{interface3BJ}
\end{eqnarray}%
where the exchange coefficient $\alpha$ represents an experimentally determined positive parameter depending on the Darcy properties. We shall also assume that all material and fluid parameters above are uniformly positive and bounded. In particular, we assume
$0\leq k^{\min }\leq \lambda (\mathbb{K}^{-1}(\mathbf{x}))\leq k^{\max}<\infty$.

Since usually the hydraulic conductivity tensor $\mathbb{K}(\mathbf{x}) $ is physically impossible to determine, we will further study the random Stokes-Darcy model with a random hydraulic conductivity tensor $\mathbb{K}(\mathbf{x},\omega) $. Let  $(\Pi,\mathcal{F}, \mathcal{P})$ be a complete probability space. Here $\Pi$ is the set of outcomes, $\mathcal{F} \in 2^{ \Pi} $ is the $\sigma$-algebra of events, and $\mathcal{P}: \mathcal{F} \rightarrow [0,1]$ is a probability measure.
The random Stokes-Darcy system reads: Find the functions $\mathbf{u}_{S}: \ \Omega_S \times \Pi \rightarrow \mathbb{R}^d \ (d=2,3)$,  ${p}_{S}: \ \Omega_S \times  \Pi \rightarrow \mathbb{R} $, $\mathbf{u}_{D}: \ \Omega_D \times  \Pi \rightarrow \mathbb{R}^d \ (d=2,3)$, and  ${\phi}_{D}: \ \Omega_D \times \Pi \rightarrow \mathbb{R} $, such that it holds $\mathcal{P} - a.e.$ in $ \Pi $:
\begin{eqnarray}
	-\nabla \cdot \mathbb{T}(\mathbf{u}_{S}(\mathbf{x},\omega),p_{S}(\mathbf{x},\omega))&=& \mathbf{f}_{S}(\mathbf{x},\omega) \ \hspace{16mm} \ \ \mathrm{in}~\Omega _{S}\times  \Pi,
	\label{StStokes1} \\
	\nabla \cdot \mathbf{u}_{S}(\mathbf{x},\omega) &=&0\ \ \hspace{24mm} \ \ \ \mathrm{in}~\Omega _{S}\times  \Pi,\label{StStokes2}\\
	\mathbb{K}^{-1}(\mathbf{x},\omega) \mathbf{u}_{D}(\mathbf{x},\omega)&=& \nabla \phi_{D}(\mathbf{x},\omega)\ \ \ \ \ \ \hspace{8.4mm} 
	\mathrm{in}~\Omega_{D}\times \Pi,  \label{StDarcy1} \\
	\nabla \cdot \mathbf{u}_{D}(\mathbf{x},\omega) &=&  {f}_{D}(\mathbf{x},\omega) \ \ \ \ \ \ \ \ \ \hspace{8.3mm} %
	\mathrm{in}~\Omega _{D}\times \Pi, \label{StDarcy2}
\end{eqnarray}
where  $\mathbf{f}_{S}:\  \Omega_S \times \Pi \rightarrow \mathbb{R}^d, \ {f}_{D}(\mathbf{x},\omega):\ \Omega_D \times \Pi \rightarrow \mathbb{R} $, and $\mathbb{K}(\mathbf{x},\omega) $ is continuous and bounded.

One of the most classical approaches to solving random PDEs is the Monte Carlo method. The main idea behind this method is that the identically distributed approximations of the solutions are computed through repeated sampling of the input parameters and then by solving the corresponding deterministic PDEs using standard numerical methods. Finally, the approximate solutions are further analyzed to obtain useful statistical information. The detailed computation procedures can be summarized as follows:
\vspace{0.1cm}

\begin{enumerate}
	\item Generate a number of independently, identically distributed (i.d.d) samples for the random hydraulic conductivity $\mathbb{K}(\mathbf{x},\omega_j)$ and the random forces $\mathbf{f}_{S}(\mathbf{x},\omega_j), \ f_D(\mathbf{x},\omega_j)$, where $ j=1,\cdots,J$;
	\item Apply standard numerical method to solve for approximate solutions $\mathbf{u}_{S}(\mathbf{x},\omega_j)$, $p_{S}(\mathbf{x},\omega_j)$, $\mathbf{u}_{D}(\mathbf{x},\omega_j)$, $\phi_{D}(\mathbf{x},\omega_j), \ j=1,\cdots, J $;
	\item Output required statistical information, such as the expectation of the free fluid velocities $\mathbf{u}_{S}(\mathbf{x},\omega): \ E[\mathbf{u}_{S}(\mathbf{x},\omega)] \approx \frac{1}{J} \sum_{j=1}^{J}\mathbf{u}_{S}(\mathbf{x},\omega_j)$.
\end{enumerate}
\vspace{0.1cm}

\begin{rem}
	When using a standard numerical method to solve the corresponding deterministic PDEs, we usually need to solve linear systems of the form
	\begin{eqnarray*}
		\begin{aligned}
			&A_j(\mathbf{x})\left[\begin{array}{l}
				\mathbf{u}_{S}(\mathbf{x},\omega_j) \\
				p_{Sj}(\mathbf{x},\omega_j)
			\end{array}
			\right]=\left[RHS_{j}(\mathbf{x})\right], \hspace{3mm}
			&B_j(\mathbf{x})\left[\begin{array}{l}
				\mathbf{u}_{Dj}(\mathbf{x},\omega_j) \\
				\phi_{Dj}(\mathbf{x},\omega_j)
			\end{array}
			\right]=\left[RHS^*_{j}(\mathbf{x})\right] \hspace{3mm}
			j=1,\cdots,J.
		\end{aligned}
\end{eqnarray*}
  As shown above, to better estimate the uncertainty and sensitivity in the solution, we shall select more samples. However, both deriving and solving such algebraic equations with varying stiffness matrices will essentially increase the computational cost. 
\end{rem}

When regarding random hydraulic conductivity $\mathbb{K}(\mathbf{x},\omega_j)$ as $\mathbb{K}_j(\mathbf{x})$, similarly denoting forces $\mathbf{f}_{S}(\mathbf{x},\omega_j)$ and $f_{D}(\mathbf{x}, \omega_j)$
as $\mathbf{f}_{Sj}(\mathbf{x})$ and $f_{Dj}(\mathbf{x})$, we can get $J$ Stokes-Darcy systems. 
Then we have an ensemble of $J$ Stokes-Darcy systems corresponding to $J$ different parameter sets $(\mathbf{f}_{Sj}(\mathbf{x}),f_{Dj}(\mathbf{x}),\mathbb{K}_j(\mathbf{x}))$, $j=1,...,J$ to be computed as follows:
\begin{eqnarray}
	-\nabla \cdot \mathbb{T}(\mathbf{u}_{Sj},p_{Sj})&=&\mathbf{f}_{Sj}\ \hspace{16mm} \ \ \mathrm{in}~\Omega _{S},
	\label{EStokes1} \\
	\nabla \cdot \mathbf{u}_{Sj} &=&0\ \ \hspace{16.5mm} \ \ \ \mathrm{in}~\Omega _{S},\label{EStokes2}\\
	\mathbb{K}_j^{-1}(\mathbf{x}) \mathbf{u}_{Dj}&=& -\nabla \phi_{Dj}\ \ \ \ \ \ \hspace{5.5mm} 
	\mathrm{in}~\Omega_{D} ,  \label{EDarcy1} \\
	\nabla \cdot \mathbf{u}_{Dj} &=& {f}_{Dj}\ \ \ \ \ \ \ \ \ \hspace{8mm}
	\mathrm{in}~\Omega _{D}. \label{EDarcy2}
\end{eqnarray}

\section{ Ensemble Domain Decomposition Method}
In this paper, we mainly focus on proposing a novel numerical method for the second procedure of the Monte Carlo method stated above. In order to solve the ensemble of $J$ Stokes-Darcy systems by domain decomposition method, one natural thought is to consider Robin-type conditions for $J$ Stokes equations and $J$ Darcy systems, since Robin-type conditions can embody both the Neumann- and Dirichlet-type conditions in (\ref{interface1})-(\ref{interface3BJ}) \cite{Lions, Chen}.

To this end,  we will study two Robin-type condition sets for $J$ Stokes equations.
For a given constant $\delta_S>0$, two corresponding functions $g_{Sj},\ g_{Sj,\tau}$ are defined on $\Gamma$:
\begin{eqnarray}
&&g_{Sj}=-\mathbf{n}_{S}\cdot (\mathbb{T}(\mathbf{u}_{Sj},p_{Sj})\cdot \mathbf{n}_{S})
-\delta_S{\mathbf{u}}_{Sj}\cdot \mathbf{n}_{S},\label{Robintype1} \\
&&g_{Sj,\tau}=-\sum_{i=1}^{d-1}(\tau_{i}\cdot \mathbb{T}(\mathbf{u}_S,p_S)\cdot \mathbf{n}_{S})
-\sum_{i=1}^{d-1}\frac{\alpha}{\sqrt{\tau_i\cdot \mathbb{K}_j\tau_i}}\mathbf{u}_{Sj} \cdot \tau _{i}.\label{Robintype2}
\end{eqnarray}
Similarly, we can propose a Robin-type condition set for the porous media flow subproblem of $J$ Darcy systems.
For a given constant $\delta_D>0$, a function $g_{Dj}$ on $\Gamma$ is constructed as:
\begin{eqnarray}
g_{Dj}&=&g\phi_{Dj}-\delta_D\mathbf{u}_{Dj} \cdot\mathbf{n}_D.\label{Robintype3}
\end{eqnarray}

The following Lemma will describe the equivalence of the original interface conditions
(\ref{interface1})-(\ref{interface3BJ}) and the above Robin-type
conditions (\ref{Robintype1})-(\ref{Robintype3}).

\begin{lemma}
The interface  conditions
(\ref{interface1})-(\ref{interface3BJ}) are equivalent to the Robin-type
conditions (\ref{Robintype1})-(\ref{Robintype3}) if and
only if $g_{Sj}, g_{Sj,\tau} $ and $\ g_{Dj}$ satisfy the following  compatibility conditions on $\Gamma$:
\begin{eqnarray}
&&g_{Dj} =g_{Sj}+(\delta_S+\delta_D)
 {\mathbf{u}}_{Sj}\cdot
\mathbf{n}_{S}+gz, \label{comp1}\\
&&g_{Sj} =g_{Dj}+(\delta_S+\delta_D)
{\mathbf{u}}_{Dj}\cdot \mathbf{n}_{D}-gz,\label{comp2}\\
&&g_{Sj,\tau}=-\sum_{i=1}^{d-1}\frac{\alpha}{\sqrt{\tau_i\cdot\mathbb{K}_j\tau_i}}\mathbf{u}_{Dj} \cdot \tau_{i}.\label{comp3}
\end{eqnarray}
\end{lemma}

\begin{proof}
The proof is much trivial, and the reader is referred to \cite{Chen, Cao11, Sun2021} for the main derivation.
\end{proof}

Before proposing our ensemble domain decomposition method, we also need to introduce some notations, Sobolev spaces, and norms. For the fluid domain $\Omega_S $ and the porous media domain $\Omega_D $, the inner products are denoted by $(\cdot, \cdot )_{S}$ and $(\cdot, \cdot )_{D}$ respectively, and the corresponding $L^2$-norms are denoted by $||\cdot ||_S$ and $||\cdot ||_D$. Moreover, $\langle \cdot, \cdot \rangle$ is defined as the $L^2$ inner product on the interface $\Gamma$, and the related $L^{2}(\Gamma )$ norm is denoted by $||\cdot ||_{\Gamma }$.

Some useful Sobolev spaces are introduced as in \cite{Sun2021}
\begin{eqnarray*}
H(\mathrm{div};\Omega _{D})&:=&\{\mathbf{v}_{D}\in L^{2}(\Omega
_{D})^{d}:\nabla \cdot \mathbf{v}_{D}\in L^{2}(\Omega _{D})\}.\\
\mathbf{X}_{S}&:=&\Big\{\mathbf{v}_{S}\in H^{1}(\Omega _{S})^{d}:\mathbf{v}_{S}=0~%
\mathrm{on}~\Gamma _{S}; ~~ \int_{\Gamma}\mathbf{v}_{S} \cdot \mathbf{n}_S=0 ~\mathrm{on}~\Gamma\Big\},\hspace{17mm}
Q_{S}:=L_0^{2}(\Omega _{S}), \\
\mathbf{X}_{D}&:=&\Big\{\mathbf{v}_{D}\in H(\mathrm{div};\Omega
_{D}): \mathbf{v}_D \cdot \mathbf{n}_D=0~%
\mathrm{on}~\Gamma _{D};~~ \int_{\Gamma}\mathbf{v}_{D} \cdot \mathbf{n}_D=0 ~\mathrm{on}~\Gamma \Big\}, \hspace{3mm} Q_{D}:=L_0^{2}(\Omega _{D}),
\end{eqnarray*}
equipped with the norms as
\begin{eqnarray*}
||\mathbf{v}_{S}||_{1}=\sqrt{||\mathbf{v}_{S}||_S^{2}+||\nabla \mathbf{v}_{S}||_{S}^{2}} \ \ \ \forall\  \mathbf{v}_{S}\in \mathbf{X}_{S}, \hspace{4mm}
||\mathbf{v}_{D}|| _{\mathrm{div}}=\sqrt{||\mathbf{v}_{D}||_D ^{2}+
||\nabla \cdot \mathbf{v}_{D}||_D^{2}} \ \ \ \forall\  \mathbf{v}_{D}\in
\mathbf{X}_{D}.
\end{eqnarray*}


With the above notations, the weak formulation of the decoupled Stokes and Darcy model with Robin-type boundary conditions can be written as follows: For the function sets
$g_{Sj}, g_{Sj,\tau}, \ g_{Dj}\in L^2(\Gamma)$, find $({\mathbf{u}}_{Sj},{p_{Sj}}; {\mathbf{u}}_{Dj},{{\phi_{Dj}}})
\in  (\mathbf{X}_{S}, Q_{S}; \mathbf{X}_{D},Q_{D})$ satisfying the compatibility conditions (\ref{comp1})-(\ref{comp3}) on the interface $\Gamma$,
such that for all $(\mathbf{v}_S,q;\mathbf{v}_D,\psi)\in (\mathbf{X}_{S},Q_S$; $\mathbf{X}_D,Q_D)$
\begin{eqnarray}
	a_{S}(\mathbf{u}_{Sj},\mathbf{v}_{S})-b_S(p_{Sj},\mathbf{v}_S)
	&+&\delta_S\langle\mathbf{u}_{Sj}\cdot\mathbf{n}_S,\mathbf{v}_{S}\cdot\mathbf{n}_S\rangle+\sum_{i=1}^{d-1} \xi_{i,j}
	\langle\mathbf{u}_{Sj}\cdot \tau_i,\mathbf{v}_{S}\cdot \tau_i\rangle\nonumber\\
	&=&(\mathbf{f}_{Sj},\mathbf{v}_S)_{S}-\langle g_{Sj},\mathbf{v}_S\cdot \mathbf{n}_S\rangle-\sum_{i=1}^{d-1}\langle g_{Sj,\tau},\mathbf{v}_S\cdot \tau_{i}\rangle,\label{StokesR1} \\
	b_S(q,\mathbf{u}_{Sj})&=&0,\label{StokesR2}\\
	a_{Dj}(\mathbf{u}_{Dj},\mathbf{v}_{D})-b_D(\phi_{Dj},\mathbf{v}_D)
	&+&\delta_D\langle\mathbf{u}_{Dj}\cdot\mathbf{n}_D,\mathbf{v}_D\cdot\mathbf{n}_D\rangle\nonumber\\
	&=&k_j^{\mathrm{min}} g(f_{Dj},\mathrm{div} \hspace{0.5mm} \mathbf{v}_D)_{D}
	-\langle g_{Dj},\mathbf{v}_D\cdot\mathbf{n}_D\rangle,\label{DarcyR1} \\
	b_D(\psi,\mathbf{u}_{Dj})&=&g(f_{Dj},\psi)_{D},\label{DarcyR2}
\end{eqnarray}
where $\xi_{i,j} = \frac{\alpha}{\sqrt{\tau_i\cdot\mathbb{K}_j\tau_i}}$, $k_j^{\min}$ and $k_j^{\max}$ are the minimum and maximum eigenvalues of $\mathbb{K}_j^{-1}(\mathbf{x}) $.
Hereafter, the bilinear forms are used:
\begin{eqnarray*}
	a_S(\mathbf{u}_{S},\mathbf{v}_{S})&=&2\nu(\mathbb{D}(\mathbf{u}_S),\mathbb{D}(\mathbf{v}_S))_{S}, \hspace{47.5mm}
	b_S(q,\mathbf{v}_S)=(q,\triangledown\cdot \mathbf{v}_S)_{S},\\
	a_{Dj}(\mathbf{u}_{D},\mathbf{v}_{D})&=&g(\mathbb{K}_j^{-1}\mathbf{u}_D,\mathbf{v}_D)_{D}
	+k_{j}^{\min} g(\mathrm{div}\hspace{0.5mm}\mathbf{u}_D,\mathrm{div}\hspace{0.5mm}\mathbf{v}_D)_{D}, \hspace{12mm}
	b_D(\psi,\mathbf{v}_D)=g(\psi,\triangledown\cdot \mathbf{v}_D)_{D}.
\end{eqnarray*}

It is worth mentioning that the well-posedness for the above weak formulation of $J$ Stokes-Darcy decoupled models can be verified easily. Let

\begin{eqnarray*}
	&&\overline{\xi}_{i}=\frac{1}{J}\sum_{j=1}^J \xi_{i,j}, \hspace{6mm}
	\overline{k}^{\min}=\frac{1}{J}\sum_{j=1}^J k_j^{\min}, \hspace{6mm}
	\overline{\mathbb{K}}=\frac{1}{J}\sum_{j=1}^J \mathbb{K}_j^{-1},\\
	&&\overline{a_{D}}(\mathbf{u}_{D},\mathbf{v}_{D})=g(\overline{\mathbb{K}}\hspace{0.5mm}\mathbf{u}_D,\mathbf{v}_D)_{D}
	+\overline{k}^{\min} g(\mathrm{div}\hspace{0.5mm}\mathbf{u}_D,\mathrm{div}\hspace{0.5mm}\mathbf{v}_D)_{D}.
\end{eqnarray*}
Then, we can propose the following parallel ensemble domain
decomposition method.

\textbf{\emph{{Ensemble DDM Algorithm}}}

1. Initial values of $g_{Sj}^0$, $g_{Sj,\tau}^0$ and $g_{Dj}^0$ are guessed. $\mathbf{u}_{Sj}^0, \mathbf{u}_{Dj}^0$ should also be given to accomplish the idea of ensemble. Both groups of the initial values are possibly taken as zeros.

2. For $n=1,2,\cdots,$ independently solve the Stokes and Darcy equations with Robin-type boundary conditions.  Namely, find $(\mathbf{v}_S,q;\mathbf{v}_D,\psi)\in (\mathbf{X}_{S},Q_S;\mathbf{X}_D,Q_D)$, solve  $(\mathbf{u}_{Sj}^{n},p_{Sj}^{n};\mathbf{u}_{Dj}^{n},\phi_{Dj}^{n}) \in  (\mathbf{X}_{S}, Q_{S};
\mathbf{X}_{D},Q_{D})$ by solving
\begin{eqnarray}\label{decoupled-1}
	&&\hspace{-8mm}a_{S}(\mathbf{u}_{Sj}^n,\mathbf{v}_{S})
	-b_S(p_{Sj}^n,\mathbf{v}_S)
	+\delta_S\langle\mathbf{u}_{Sj}^n\cdot\mathbf{n}_S,\mathbf{v}_{S}\cdot\mathbf{n}_S\rangle
	+\sum_{i=1}^{d-1} \overline{\xi}_{i}
		\langle\mathbf{u}_{Sj}^n\cdot \tau_i,\mathbf{v}_{S}\cdot \tau_i\rangle \nonumber\\
 	&&\hspace{-2mm}=(\mathbf{f}_{Sj},\mathbf{v}_S)_{S}-\langle g_{Sj}^{n-1},\mathbf{v}_S\cdot \mathbf{n}_S\rangle-\sum_{i=1}^{d-1}\langle g_{Sj,\tau}^{n-1},\mathbf{v}_S\cdot \tau_{i}\rangle +\sum_{i=1}^{d-1} (\overline{\xi}_{i}-\xi_{i,j})
		\langle\mathbf{u}_{Sj}^{n-1}\cdot \tau_i,\mathbf{v}_{S}\cdot \tau_i\rangle,\\
	&&\hspace{-8mm}b_S(q,\mathbf{u}_S^n)=0,\label{decoupled-2}\\
	&&\hspace{-8mm}\overline{a_{D}}(\mathbf{u}_{Dj}^n,\mathbf{v}_{D})-b_D(\phi_{Dj}^n,\mathbf{v}_D)
	+\delta_D\langle\mathbf{u}_{Dj}^n\cdot\mathbf{n}_D,\mathbf{v}_D\cdot\mathbf{n}_D\rangle\nonumber\\
	&&\hspace{-2mm}= k_j^{\mathrm{min}} g(f_{Dj},\mathrm{div} \hspace{0.5mm} \mathbf{v}_D)_{D}
	-\langle g_{Dj}^{n-1},\mathbf{v}_D\cdot\mathbf{n}_D\rangle +\Big[a_{Dj}(\mathbf{u}_{Dj}^{n-1},\mathbf{v}_{D})-\overline{a_{D}}(\mathbf{u}_{Dj}^{n-1},\mathbf{v}_{D})\Big],\label{decoupled-3} \\
	&&\hspace{-8mm}b_D(\psi,\mathbf{u}_{Dj}^n)=g(f_{Dj},\psi)_{D}.\label{decoupled-4}
\end{eqnarray}

3. Update $g_{Sj}^{n},\ g_{Sj,\tau}^{n}$ and $g_{Dj}^{n}$ in the following manner:
\begin{eqnarray}
	&&g_{Dj}^{n} =g_{Sj}^{n-1}+(\delta_S+\delta_D)
	{\mathbf{u}}_{Sj}^n\cdot
	\mathbf{n}_{S}+gz,\label{decoupled-comp1}\\
	&&g_{Sj}^{n} =g_{Dj}^{n-1}+(\delta_S+\delta_D)
	{\mathbf{u}}_{Dj}^n\cdot \mathbf{n}_{D}-gz,\label{decoupled-comp2}\\
	&&g_{Sj,\tau}^{n}=-\sum_{i=1}^{d-1}\xi_{i,j}\mathbf{u}_{Dj}^n \cdot \tau_{i}.\label{decoupled-comp3}
\end{eqnarray}

Noting that $J$ Stokes equations and $J$ Darcy equations are decoupled and independent systems. So for any given initial guesses and fixed $n$, the existence and uniqueness of the solution   $(\mathbf{u}_{Sj}^{n},p_{Sj}^{n})$ and $(\mathbf{u}_{Dj}^{n},\phi_{Dj}^{n})$ to each system follow immediately.

\begin{rem}
	In each iteration step, we only  compute the solutions of two linear systems (for $(\mathbf{u}_{Sj}, p_{Sj})$ and $(\mathbf{u}_{Dj}, \phi_{Dj})$ respectively), which share the same coefficient matrix of the form
	\begin{eqnarray*}
		\begin{aligned}
			&A\left[\begin{array}{l l l}
				\mathbf{u}_{S1} \\
				p_{S1}
			\end{array}
			|\cdots| 
			\begin{array}{l}
				\mathbf{u}_{SJ} \\
				p_{SJ}
			\end{array}\right]=\left[RHS_{1}|\cdots| RHS_{J}\right], \hspace{4mm}
			&B\left[\begin{array}{l l l}
				\mathbf{u}_{D1} \\
				\phi_{D1}
			\end{array}|\cdots| \begin{array}{l}
				\mathbf{u}_{DJ} \\
				\phi_{DJ}
			\end{array}\right]=\left[RHS^*_{1}|\cdots| RHS^*_{J}\right].
		\end{aligned}
	\end{eqnarray*}
	Hence the coefficient matrices $A$ and $B$ only need to use once efficient iterative solvers or direct solvers such as LU factorization for fast computation. Moreover, the domain decomposition method will also bring out better parallel efficiency. 
\end{rem}

\section{ The Convergence of Ensemble Domain Decomposition Method}
In this section, we will demonstrate the convergence of the parallel ensemble domain decomposition method by applying the elegant energy method. Moreover, we will try to overcome the assumption that the exchange coefficient $\alpha$ in the BJ interface conditions is sufficiently small, and prove that the constructed algorithm has a mesh-independent convergence rate in the case of  $\delta_S<\delta_D $. This section will focus on the cases of $\delta_S<\delta_D$ and $\delta_S=\delta_D $. The convergence of the case $\delta_S>\delta_D$ is much more difficult to prove in the continuous form and will be demonstrated in the next section. 

The main result concerning the convergence of our algorithm is listed in the following theorem.
\begin{theorem}
Assume that $(\mathbf{u}_{Sj}^n,p_{Sj}^n;\mathbf{u}_{Dj}^n,\phi_{Dj}^n)$ and $(\mathbf{u}_{Sj},p_{Sj};\mathbf{u}_{Dj},\phi_{Dj})$ are the solutions of Ensemble DDM Algorithm and the DDM weak formulation (\ref{StokesR1})-(\ref{DarcyR2}), respectively. Then if $\delta_S\leq\delta_D $, $(\mathbf{u}_{Sj}^n,p_{Sj}^n;\mathbf{u}_{Dj}^n,\phi_{Dj}^n)$  will converge to $(\mathbf{u}_{Sj},p_{Sj};\mathbf{u}_{Dj},\phi_{Dj})$.
\end{theorem}

\begin{proof}
Define the following error functions:
\begin{eqnarray*}
&&\mathbf{e}^n_{Sj}=\mathbf{u}_{Sj}-\mathbf{u}^n_{Sj},\hspace{3mm}\mathbf{e}^n_{Dj}=\mathbf{u}_{Dj}-\mathbf{u}^n_{Dj},\hspace{8mm}
  {\varepsilon}_{Sj}^{n}=p_{Sj}-{p}_{Sj}^{n},\hspace{4mm}
{\varepsilon}_{Dj}^{n}=\phi_{Dj}-{\phi}_{Dj}^{n},\\
&&  {\eta}_{Sj}^{n}=g_{Sj}-{g}^{n}_{Sj},\hspace{4mm}
{\eta}_{Sj,\tau}^{n}=g_{Sj,\tau}-{g}^{n}_{Sj,\tau}, \hspace{4mm}
{\eta}_{Dj}^{n}=g_{Dj}-{g}^{n}_{Dj}.
\end{eqnarray*}
Suppose that $\tilde{k}_j^{\mathrm{max}}$ is the maximum absolute value among all eigenvalues of the matrix $\mathbb{K}_j^{-1}-\overline{\mathbb{K}}$, and define the following quantities:
\begin{eqnarray*}
	&& E_{{\xi}_{i,j}}^{\max}=\max|{\xi_{i,j}}-\overline{\xi}_i|,  \hspace{8mm} E_{{\xi}_{j}}^{\max}=\sum_{i=1}^{d-1} E_{{\xi}_{i,j}}^{\max}, \hspace{8mm}
	\overline{\xi}=\sum_{i=1}^{d-1} \overline{\xi}_i, \hspace{8mm}
	\xi_{j}=\sum_{i=1}^{d-1}{\xi_{i,j}}, \\
	&& E_{k_j}^{\mathrm{max}}= \max(\tilde{k}_j^{\mathrm{max}}, | k_j^{\mathrm{min}}-\overline{k}^{\mathrm{min}}| ), \hspace{10mm}
		\overline{k}^{\max}=\frac{1}{J}\sum_{j=1}^J k_j^{\max}.
\end{eqnarray*}
Then for all $(\mathbf{v}_S,q;\mathbf{v}_D,\psi)\in (\mathbf{X}_{S},Q_S;\mathbf{X}_D,Q_D)$, subtract
(\ref{decoupled-1})-(\ref{decoupled-4}) from
(\ref{StokesR1})-(\ref{DarcyR2}):
\begin{eqnarray}\label{err-1}
&&a_{S}(\mathbf{e}_{Sj}^n,\mathbf{v}_{S})
-b_S(\varepsilon_{Sj}^n,\mathbf{v}_S)
+\delta_S\langle\mathbf{e}_{Sj}^n\cdot\mathbf{n}_S,\mathbf{v}_{S}\cdot\mathbf{n}_S\rangle
+\sum_{i=1}^{d-1} \overline{\xi}_{i}
\langle\mathbf{e}_{Sj}^n\cdot \tau_i,\mathbf{v}_{S}\cdot \tau_i\rangle \nonumber\\
&&\hspace{15.5mm}=-\langle\eta_{Sj}^{n-1},\mathbf{v}_S\cdot \mathbf{n}_S\rangle-\sum_{i=1}^{d-1}\langle\eta_{Sj,\tau}^{n-1},\mathbf{v}_S\cdot \tau_{i}\rangle +\sum_{i=1}^{d-1} (\overline{\xi}_{i}-\xi_{i,j})
\langle\mathbf{e}_{Sj}^{n-1}\cdot \tau_i,\mathbf{v}_{S}\cdot \tau_i\rangle,\\
&&b_S(q,\mathbf{e}_{Sj}^n)=0,\label{err-2}\\
&&\overline{a_{D}}(\mathbf{e}_{Dj}^n,\mathbf{v}_{D})-b_D(\varepsilon_{Dj}^n,\mathbf{v}_D)
+\delta_D\langle\mathbf{e}_{Dj}^n\cdot\mathbf{n}_D,\mathbf{v}_D\cdot\mathbf{n}_D\rangle\nonumber\\
&&\hspace{17mm}=-\langle\eta_{Dj}^{n-1},\mathbf{v}_D\cdot\mathbf{n}_D\rangle
+\Big[a_{Dj}(\mathbf{e}_{Dj}^{n-1},\mathbf{v}_{D})-\overline{a_{D}}(\mathbf{e}_{Dj}^{n-1},\mathbf{v}_{D})\Big],\hspace{1mm}\label{err-3} \\
&&b_D(\psi,\mathbf{e}_{Dj}^n)=0.\label{err-4}
\end{eqnarray}
Along the interface $\Gamma$, the error functions can be updated as follows
\begin{eqnarray}
&&{\eta}_{Dj}^{n} ={\eta}_{Sj}^{n-1}+(\delta_S+\delta_D)
{\mathbf{e}}_{Sj}^n\cdot
\mathbf{n}_{S},\label{errep}\\
&&{\eta}_{Sj}^{n} ={\eta}_{Dj}^{n-1}+(\delta_S+\delta_D)
{\mathbf{e}}_{Dj}^n\cdot \mathbf{n}_{D},\label{erref}\\
&&{\eta}_{Sj,\tau}^{n}=-\sum_{i=1}^{d-1} \xi_{i,j}
\mathbf{e}_{Dj}^n \cdot \tau_{i}.\label{erref2}
\end{eqnarray}
Equation (\ref{errep}) can lead to
\begin{eqnarray}\label{etaD2}
||\eta_{Dj}^{n}||^2_{\Gamma} =||\eta_{Sj}^{n-1}||^2_{\Gamma}+2(\delta_S+\delta_D)\langle\eta_{Sj}^n,\mathbf{e}_{Sj}^{n}\cdot
\mathbf{n}_{S}\rangle+(\delta_S+\delta_D)^2||\mathbf{e}_{Sj}^{n}\cdot\mathbf{n}_{S}||^2_{\Gamma}.
\end{eqnarray}
Choosing $(\mathbf{v}_{S},q)=(\mathbf{e}_{Sj}^{n},\varepsilon_{Sj}^{n})$ in (\ref{err-1})-(\ref{err-2}) and together with (\ref{erref2}), we can get
\begin{eqnarray}\label{aS+bS}
&&\hspace{-5mm}a_{S}(\mathbf{e}_{Sj}^n,\mathbf{e}_{Sj}^n)+\delta_S||\mathbf{e}_{Sj}^{n}\cdot\mathbf{n}_{S}||^2_{\Gamma}
+\sum_{i=1}^{d-1}\overline{\xi}_i||\mathbf{e}_{Sj}^n\cdot \tau_i||_{\Gamma}^2\nonumber\\
&&=-\langle{\eta}_{Sj}^{n-1},\mathbf{e}_{Sj}^n\cdot \mathbf{n}_S\rangle
+\sum_{i=1}^{d-1}\xi_{i,j}\langle\mathbf{e}_{Dj}^{n-1}\cdot \tau_{i},\mathbf{e}_{Sj}^n\cdot \tau_{i}\rangle +\sum_{i=1}^{d-1}(\overline{\xi}_i-\xi_{i,j})\langle\mathbf{e}_{Sj}^{n-1}\cdot \tau_{i},\mathbf{e}_{Sj}^n\cdot \tau_{i}\rangle.
\end{eqnarray}
Combining (\ref{etaD2}) and (\ref{aS+bS}), an important relation can be derived as
\begin{eqnarray}\label{errorD}
||\eta_{Dj}^{n}||^2_{\Gamma}
&=&||\eta_{Sj}^{n-1}||^2_{\Gamma}
+(\delta_D^2-\delta_S^2)||\mathbf{e}_{Sj}^{n}\cdot\mathbf{n}_{S}||^2_{\Gamma}
-2(\delta_S+\delta_D)a_S(\mathbf{e}_{Sj}^n,\mathbf{e}_{Sj}^n)
\nonumber\\
&&-2(\delta_S+\delta_D)\sum_{i=1}^{d-1}\overline{\xi}_i
||\mathbf{e}_{Sj}^n\cdot \tau_{i}||_{\Gamma}^2
+2(\delta_S+\delta_D)\sum_{i=1}^{d-1}{\xi}_{i,j}
\langle\mathbf{e}_{Dj}^{n-1}\cdot \tau_{i},\mathbf{e}_{Sj}^n\cdot \tau_{i}\rangle\nonumber\\
&&-2(\delta_S+\delta_D)\sum_{i=1}^{d-1}({\xi}_{i,j}-\overline{\xi}_i)
\langle\mathbf{e}_{Sj}^{n-1}\cdot \tau_{i},\mathbf{e}_{Sj}^n\cdot \tau_{i}\rangle.
\end{eqnarray}
Similarly,  using (\ref{erref}), taking $(\mathbf{v}_{D},\psi)=(\mathbf{e}_{Dj}^{n},\varepsilon_{Dj}^{n})$ in (\ref{err-3})-(\ref{err-4}) to get another important relation
\begin{eqnarray}\label{errorS}
||\eta_{Sj}^{n}||^2_{\Gamma}&=& ||\eta_{Dj}^{n-1}||^2_{\Gamma}
+(\delta_S^2-\delta_D^2)||\mathbf{e}_{Dj}^{n}\cdot\mathbf{n}_{D}||^2_{\Gamma}
-2(\delta_S+\delta_D)\overline{a_D}(\mathbf{e}_{Dj}^n,\mathbf{e}_{Dj}^n)\nonumber\\
&&-2(\delta_S+\delta_D)\Big[a_{Dj}(\mathbf{e}_{Dj}^{n-1},\mathbf{e}_{Dj}^n)-\overline{a_D}(\mathbf{e}_{Dj}^{n-1},\mathbf{e}_{Dj}^n)\Big].
\end{eqnarray}

Most of the literature usually assumes that the positive parameter $\alpha$ is small enough in the cases of $\delta_S<\delta_D$ and $\delta_S=\delta_D$. However, we will show how to overcome this assumption in the case of $\delta_S<\delta_D$ through appropriate parameter selection.

For the case $\delta_S<\delta_D$, we will give the convergence proof and more importantly remove the assumption in most of the developed methods, that the positive parameter $\alpha$ is small enough. The main idea is to utilize one interface term and the following important inequality (\ref{inf2div}).

Noting that there exists a continuous and linear mapping $H(\mathrm{div};\Omega):\rightarrow H^{-1/2}(\partial\Omega_D)$ (\cite{GR86}, Corollary 2.8 pp.29) and clearly $L^2(\partial \Omega_D)\subset H^{-1/2}(\partial \Omega_D)$. Then an important inequality follows:
\begin{eqnarray}\label{inf2div}
||\mathbf{e}_{Dj}^n||_{\mathrm{div}}\leq
||\mathbf{e}_{Dj}^n\cdot\mathbf{n}_D||_{H^{-1/2}(\Gamma)}\leq
||\mathbf{e}_{Dj}^n\cdot\mathbf{n}_D||_{\Gamma}.
\end{eqnarray}
Adding equations (\ref{errorD})-(\ref{errorS}) together, then summing the resulting equation over $n$ from $n=1$ to $N$, we can derive
\begin{eqnarray}\label{etaS+etaD}
	||\eta_{Sj}^{N}||^2_{\Gamma}+||\eta_{Dj}^{N}||^2_{\Gamma}
	&=&||\eta_{Sj}^{0}||^2_{\Gamma}+||\eta_{Dj}^{0}||^2_{\Gamma}
	+(\delta_D^2-\delta_S^2)\sum_{n=1}^{N}||\mathbf{e}_{S}^{n}\cdot\mathbf{n}_{S}||^2_{\Gamma}+(\delta_S^2-\delta_D^2)\sum_{n=1}^{N}||\mathbf{e}_{D}^{n}\cdot\mathbf{n}_{D}||^2_{\Gamma}\nonumber\\
	&&-2(\delta_S+\delta_D) \sum_{n=1}^{N} \Big[ a_S(\mathbf{e}_{Sj}^n,\mathbf{e}_{Sj}^n)
	+\overline{a_D}(\mathbf{e}_{Dj}^n,\mathbf{e}_{Dj}^n)
	+\sum_{i=1}^{d-1} \overline{\xi}_i||\mathbf{e}_{Sj}^n\cdot \tau_i||_{\Gamma}^2\nonumber\\
	&&-\sum_{i=1}^{d-1} \xi_{i,j} \langle\mathbf{e}_{Dj}^{n-1}\cdot \tau_i,\mathbf{e}_{Sj}^n\cdot \tau_i\rangle
	+\sum_{i=1}^{d-1} (\xi_{i,j}-\overline{\xi}_i) \langle\mathbf{e}_{Sj}^{n-1}\cdot \tau_i,\mathbf{e}_{Sj}^n\cdot \tau_i\rangle\nonumber\\
	&&+g((\mathbb{K}_j^{-1}-\overline{\mathbb{K}})\mathbf{e}_{Dj}^{n-1},\mathbf{e}_{Dj}^{n})_{D}
	+g(k_j^{\mathrm{min}}-\overline{k}^{\mathrm{min}})(\mathrm{div}\hspace{0.5mm}\mathbf{e}_{Dj}^{n-1},\mathrm{div}\hspace{0.5mm}\mathbf{e}_{Dj}^{n})_{D}\Big].
\end{eqnarray}
Then applying Korn's inequality, for positive constant $C_1$, we can derive
\begin{eqnarray}\label{bilinear-ieq}
	a_S(\mathbf{e}_{Sj}^n,\mathbf{e}_{Sj}^n)\geq C_1\nu ||\mathbf{e}_{Sj}^n||_1^2,\hspace{8mm}
	\overline{a_D}(\mathbf{e}_{Dj}^n,\mathbf{e}_{Dj}^n)\geq g \overline{k}^{\min} ||\mathbf{e}_{Dj}^n||_{\mathrm{div}}^2.
\end{eqnarray}
Recall some trace inequalities \cite{trace}, which are useful in our analysis. There exist constants $C_{\mathrm{tr}}$, $C'_{\mathrm{tr}}$, $C''_{\mathrm{tr}}$ which only depend on the domain $\Omega_D$ or $\Omega_S$, such that for any $\mathbf{v}_D\in \mathbf{X}_D$ or $\mathbf{v}_S\in \mathbf{X}_S$:
\begin{eqnarray}\label{trace}
	||\mathbf{v}_D||_{H^{-\frac{1}{2}}(\Gamma)}\leq C_{\mathrm{tr}} ||\mathbf{v}_D||_{\mathrm{div}}, \hspace{3mm} ||\mathbf{v}_S||_{H^{\frac{1}{2}}(\Gamma)}\leq C'_{\mathrm{tr}} ||\mathbf{v}_S||_{1}, \hspace{3mm} 
	||\mathbf{v}_S||_{\Gamma}\leq C''_{\mathrm{tr}} ||\mathbf{v}_S||_{S}^{\frac{1}{2}}||\mathbf{v}_S||_{1}^{\frac{1}{2}}.\hspace{3mm}
\end{eqnarray}
By the Cauchy-Schwarz inequality,  trace inequalities (\ref{trace}), and Young's inequality, we arrive at
\begin{eqnarray}\label{term1}
	&&-\sum_{i=1}^{d-1} \xi_{i,j} \langle \mathbf{e}_{Dj}^{n-1}\cdot \tau_i,\mathbf{e}_{Sj}^n\cdot \tau_i\rangle \geq -\sum_{i=1}^{d-1} \xi_{i,j} ||\mathbf{e}_{Dj}^{n-1}\cdot \tau_i||_{H^{-\frac{1}{2}}(\Gamma)}||\mathbf{e}_{Sj}^n\cdot \tau_i||_{H^{\frac{1}{2}}(\Gamma)}\nonumber\\
	&&\hspace{8mm}\geq -C_{\mathrm{tr}}C'_{\mathrm{tr}}\xi_j ||\mathbf{e}_{Dj}^{n-1}||_{\mathrm{div}}
	||\mathbf{e}_{Sj}^{n}||_{1}
	\geq-\frac{C_1\nu}{2}||\mathbf{e}_{Sj}^{n}||_{1}^2-\frac{(C_{\mathrm{tr}}C'_{\mathrm{tr}}\xi_j)^2}{2C_1\nu}||\mathbf{e}_{Dj}^{n-1}||_{\mathrm{div}}^2,\\ \label{term2}
	&&\sum_{i=1}^{d-1} (\xi_{i,j}-\overline{\xi}_i) \langle\mathbf{e}_{Sj}^{n-1}\cdot \tau_i,\mathbf{e}_{Sj}^n\cdot \tau_i\rangle \geq -\sum_{i=1}^{d-1}\frac{E_{{\xi}_{i,j}}^{\max}}{2} ||\mathbf{e}_{Sj}^{n-1}\cdot \tau_i||_{\Gamma}^2-\sum_{i=1}^{d-1}\frac{E_{{\xi}_{i,j}}^{\max}}{2} ||\mathbf{e}_{Sj}^{n}\cdot \tau_i||_{\Gamma}^2,\\ \label{term3}
	&&g((\mathbb{K}_j^{-1}-\overline{\mathbb{K}})\mathbf{e}_{Dj}^{n-1},\mathbf{e}_{Dj}^{n})_{D}
	+g(k_j^{\mathrm{min}}-\overline{k}^{\mathrm{min}})(\mathrm{div}\hspace{0.5mm}\mathbf{e}_{Dj}^{n-1},\mathrm{div}\hspace{0.5mm}\mathbf{e}_{Dj}^{n})_{D}
	\nonumber\\
	&&\hspace{8mm} \geq -\frac{gE_{k_j}^{\max}}{2}||\mathbf{e}_{Dj}^{n-1}||_{\mathrm{div}}^2
	-\frac{gE_{k_j}^{\max}}{2}||\mathbf{e}_{Dj}^{n}||_{\mathrm{div}}^2.
\end{eqnarray}
Since $\delta_S<\delta_D$, clearly we have $\delta_S^2-\delta_D^2<0$. Substituting the inequalities (\ref{bilinear-ieq})-(\ref{term3}) and Poincar$\acute{\mathrm{e}}$ inequality with a constant $C_2>0$ into the equation (\ref{etaS+etaD}) and using (\ref{inf2div}), we can deduce
\begin{eqnarray}\label{convergence}
0&\leq& ||\eta_{Sj}^{N}||^2_{\Gamma}+||\eta_{Dj}^{N}||^2_{\Gamma}
= ||\eta_{Sj}^{0}||^2_{\Gamma}+||\eta_{Dj}^{0}||^2_{\Gamma}
+\frac{\delta_S^2-\delta_D^2}{2}\sum_{n=1}^{N}||\mathbf{e}_{Dj}^{n}\cdot\mathbf{n}_{D}||^2_{\Gamma} -2(\delta_S+\delta_D) \sum_{n=1}^{N}\nonumber\\
&&\Big[ a_S(\mathbf{e}_{Sj}^n,\mathbf{e}_{Sj}^n)
+\overline{a_D}(\mathbf{e}_{Dj}^n,\mathbf{e}_{Dj}^n)
+\sum_{i=1}^{d-1} \overline{\xi}_i||\mathbf{e}_{Sj}^n\cdot \tau_i||_{\Gamma}^2
-\frac{\delta_D-\delta_S}{2}||\mathbf{e}_{Sj}^{n}\cdot\mathbf{n}_{S}||^2_{\Gamma}\nonumber\\
&&\hspace{2mm}-\sum_{i=1}^{d-1} \xi_{i,j} \langle\mathbf{e}_{Dj}^{n-1}\cdot \tau_i,\mathbf{e}_{Sj}^n\cdot \tau_i\rangle
+\sum_{i=1}^{d-1} (\xi_{i,j}-\overline{\xi}_i) \langle\mathbf{e}_{Sj}^{n-1}\cdot \tau_i,\mathbf{e}_{Sj}^n\cdot \tau_i\rangle\nonumber\\
&&\hspace{2mm}+g((\mathbb{K}_j^{-1}-\overline{\mathbb{K}})\mathbf{e}_{Dj}^{n-1},\mathbf{e}_{Dj}^{n})_{D}
+g(k_j^{\mathrm{min}}-\overline{k}^{\mathrm{min}})(\mathrm{div}\hspace{0.5mm}\mathbf{e}_{Dj}^{n-1},\mathrm{div}\hspace{0.5mm}\mathbf{e}_{Dj}^{n})_{D} +\frac{\delta_D-\delta_S}{4}||\mathbf{e}_{Dj}^{n}\cdot\mathbf{n}_{D}||^2_{\Gamma}\Big]\nonumber\\
&\leq& ||\eta_{Sj}^{0}||^2_{\Gamma}+||\eta_{Dj}^{0}||^2_{\Gamma}
+\frac{\delta_S^2-\delta_D^2}{2}\sum_{n=1}^{N}||\mathbf{e}_{Dj}^{n}\cdot\mathbf{n}_{D}||^2_{\Gamma}
+\sum_{i=1}^{d-1}(\delta_S+\delta_D) E_{{\xi}_{i,j}}^{\max} ||\mathbf{e}_{Sj}^{0}\cdot \tau_i||_{\Gamma}^2\nonumber\\
&&+(\delta_S+\delta_D)\Big(\frac{(C_{\mathrm{tr}}C'_{\mathrm{tr}}\xi_j)^2}{C_1\nu}+ gE_{k_j}^{\max}\Big)||\mathbf{e}_{Dj}^{0}||_{\mathrm{div}}^2
-2(\delta_S+\delta_D) \sum_{n=1}^{N} \sum_{i=1}^{d-1} \Big(\overline{\xi}_i-E_{{\xi}_{i,j}}^{\max}\Big) ||\mathbf{e}_{Sj}^{n}\cdot \tau_i||_{\Gamma}^2\nonumber\\
&&-2(\delta_S+\delta_D) \sum_{n=1}^{N} \Big[\Big(\frac{C_1\nu}{2}-\frac{C_2(\delta_D-\delta_S)}{2}\Big)||\mathbf{e}_{Sj}^{n}||_{1}^2 +C_{\alpha} ||\mathbf{e}_{Dj}^{n}||_{\mathrm{div}}^2 \Big],
\end{eqnarray}
where $C_{\alpha}=g\overline{k}^{\min}-gE_{k_j}^{\max}+\frac{\delta_D-\delta_S}{4}-\frac{(C_{\mathrm{tr}}C'_{\mathrm{tr}}\xi_j)^2}{2C_1\nu}$.
Then if enforcing the conditions
\begin{eqnarray}\label{add_condition1}
&&\frac{2(C_{\mathrm{tr}}C'_{\mathrm{tr}}\xi_j)^2}{C_1\nu}-4g(\overline{k}^{\min}-E_{k_j}^{\max})<\delta_D-\delta_S<\frac{C_1\nu}{C_4},\\
\label{conditions=}
&&\overline{\xi}_i \geq E_{{\xi}_{i,j}}^{\max}, \hspace{10mm}
\overline{k}^{\min} > E_{k_j}^{\max},
\end{eqnarray}
we can directly prove that $\mathbf{e}_{Sj}^{n}$, $\mathbf{e}_{Dj}^{n}$ tend to be zeros in $H^1(\Omega_S)^d$ and $H(\mathrm{div};\Omega_D)$ respectively,  as $n\rightarrow \infty$.
For the condition (\ref{add_condition1}), if the hydraulic conductivity $\mathbb{K}_j$ and the viscosity coefficient $\nu$ are not too small, the parameters $\delta_S$ and $\delta_D$ can be found (check Remark 4.2 for details) and more importantly the parameter $\alpha$ in $\xi_{j}=\sum_{i=1}^{d-1}\frac{\alpha}{\sqrt{\tau_i\cdot\mathbb{K}_j\tau_i}}$  is not required to be small enough.

The convergence of series $||\varepsilon_{Dj}^{n}||_D$ will be proved next. Firstly, for given $\varepsilon_{Dj}^{n}\in Q_D$, there exist
$\mathbf{v}^{\varepsilon}_{D}\in \mathbf{X}_D \cap H_0^1(\Omega_D)^d$
and a positive constant $C_{II}>0$ satisfying
\begin{eqnarray*}
	\nabla\cdot \mathbf{v}^{\varepsilon}_{D}=\varepsilon_{Dj}^{n}\ \
	\mathrm{in} \ \Omega_D, \hspace{6mm} \mathbf{v}^{\varepsilon}_{D}=0 \ \ \mathrm{on}\ \partial\Omega_D, \hspace{6mm}
	|| \mathbf{v}^{\varepsilon}_{D}||_{1}\leq
	C_{II}||\varepsilon_{Dj}^{n}||_D,
\end{eqnarray*}
which yields immediately
\begin{eqnarray*}
	b_{D}(\mathbf{v}_{D}^{\varepsilon},\varepsilon_{Dj}^{n})
	=||\varepsilon_{Dj}^{n}||_{D}^2\geq
	1/C_{II}||\mathbf{v}_{D}^{\varepsilon}||_{\mathrm{div}}||\varepsilon_{Dj}^{n}||_{D}.
\end{eqnarray*}
Choosing $\mathbf{v}_{D}^{\varepsilon}\in \mathbf{X}_D$ and $\mathbf{v}_{D}^{\varepsilon}\cdot\mathbf{n}_D=0$ in equation (\ref{err-3}) gives
\begin{eqnarray*}
	\overline{a_{D}}(\mathbf{e}_{Dj}^{n},\mathbf{v}_{D}^{\varepsilon})-b_{D}(\mathbf{v}%
	_{D}^{\varepsilon},\varepsilon_{Dj}^{n})-[a_{Dj}(\mathbf{e}_{Dj}^{n-1},\mathbf{v}_{D}^{\varepsilon})-\overline{a_D}(\mathbf{e}_{Dj}^{n-1},\mathbf{v}_{D}^{\varepsilon})]=0,
\end{eqnarray*}
which further leads to
\begin{eqnarray}\label{PDcontrol}
	||\varepsilon_{Dj}^{n}||_{\Omega_D} &\leq& C_{II}\frac{b_{D}(\mathbf{v}_{D}^{\varepsilon},\varepsilon_{Dj}^{n})}{||\mathbf{v}_{D}^{\varepsilon}||_{\mathrm{div}}}
	=C_{II}\frac{\overline{a_{D}}(\mathbf{e}_{Dj}^{n},\mathbf{v}_{D}^{\varepsilon})-[a_{Dj}(\mathbf{e}_{Dj}^{n-1},\mathbf{v}_{D}^{\varepsilon})-\overline{a_D}(\mathbf{e}_{Dj}^{n-1},\mathbf{v}_{D}^{\varepsilon})]}
	{||\mathbf{v}_{D}^{\varepsilon}||_{\mathrm{div}}}\nonumber\\
	&\leq& C_3g\overline{k}^{\max} ||\mathbf{e}_{Dj}^{n}||_{\mathrm{div}}+C_3 gE_{k_j}^{\max}||\mathbf{e}_{Dj}^{n-1}||_{\mathrm{div}},
\end{eqnarray}
where constant $C_3>0$ only depends on  $\Omega_D$. This follows that  $\varepsilon_{Dj}^{n}$ tends to be zero in $L_0^2(\Omega_D)$.

For the convergence of the pressure $||\varepsilon_{Sj}^{n}||_S$, a similar way is utilized to obtain
\begin{eqnarray}\label{contralp}
	||\varepsilon_{Sj}^{n}||_S \leq C_{I}\frac{b_{S}(\mathbf{v}_{S}^{\varepsilon},\varepsilon_{Sj}^{n})}{||\mathbf{v}_{S}^{\varepsilon}||_{1}}
	=C_{I}\frac{a_{S}(\mathbf{e}_{Sj}^{n},\mathbf{v}_{S}^{\varepsilon})}
	{||\mathbf{v}_{S}^{\varepsilon}||_{1}}
	\leq C_4\nu||\mathbf{e}_{Sj}^{n}||_{1},
\end{eqnarray}
where  $C_4$ is a positive constant.
This implies that the convergence of $||\varepsilon_{Sj}^{n}||_S$ also holds in $L_0^2(\Omega_S)$.

Then the convergence of series $||\eta_{Sj,\tau}^n||_{H^{-1/2}(\Gamma)}$ is proved by using the trace inequality in (\ref{erref2})
\begin{eqnarray}\label{etaStau-control}
	||\eta_{Sj,\tau}^{n}||_{H^{-1/2}(\Gamma)}\leq \sum_{i=1}^{d-1}\xi_{i,j}
	||\mathbf{e}_{Dj}^{n-1} \cdot \tau_{i}||_{H^{-1/2}(\Gamma)}\leq  C_{\mathrm{tr}}\xi_j||\mathbf{e}_{Dj}^{n-1}||_{\mathrm{div}}.
\end{eqnarray}
which means that  $\eta_{Sj,\tau}^{n}$ converges to zero in $H^{-1/2}(\Gamma)$.
Moreover, in order to demonstrate the convergence of $\eta_{Dj}^n$ and $\eta_{Sj}^n$ in $H^{-1/2}(\Gamma)$, we can basically use (\ref{PDcontrol}) and follow the analysis of \cite{Sun2021}.

For the case $\delta_S=\delta_D$, if the condition (\ref{add_condition1})  holds, the parameter $\alpha$ in $\xi_{j}=\sum_{i=1}^{d-1}\frac{\alpha}{\sqrt{\tau_i\cdot\mathbb{K}\tau_i}}$ must be sufficiently small, as the results of existing studies \cite{Cao10, Cao11}.  The convergence analysis for all other quantities is similar to the discussion in the case $\delta_S<\delta_D$ above, and is omitted here.
%
\end{proof}

\begin{rem}\label{4.1}
	For the restriction (\ref{add_condition1}) on parameters $\delta_S$ and $\delta_D$, we need further analysis on the influences of some physical parameters $\alpha$, $\mathbb{K}_j$ and $\nu$.
	If the experimentally determined parameter $\alpha$ is small enough, as assumed in other literature, we can directly have $-4g(\overline{k}^{\min}-E_{k_j}^{\max})<0$. Then we can derive a very similar conclusion with the case of $\delta_S=\delta_D$.
	As we have discussed, real applications require that the positive parameter $\alpha$ should not be small enough. 
    Since $\xi_{j}=\sum_{i=1}^{d-1}\frac{\alpha}{\sqrt{\tau_i\cdot\mathbb{K}_j\tau_i}}$, we have $
	\sum_{i=1}^{d-1}\frac{2(C_{\mathrm{tr}}C'_{\mathrm{tr}}\alpha)^2}{C_1\nu\tau_i\cdot\mathbb{K}_j\tau_i}-4g(\overline{k}^{\min}-E_{k_j}^{\max})<\frac{C_1\nu}{C_4}$.
	So we only require that
	\begin{eqnarray}\label{limtK}
		(C_{\mathrm{tr}}C'_{\mathrm{tr}}\alpha)^2{k}_j^{\max}
		<2gC_1\nu(\overline{k}^{\min}-E_{k_j}^{\max})+\frac{(C_1\nu)^2}{2C_4}.
	\end{eqnarray}
	As shown above, we can suppose that the hydraulic conductivity tensor $\mathbb{K}_j$ and the viscosity coefficient $\nu$ are not sufficiently small, maybe $O(1)$ or larger to overcome the assumption in \cite{Cao10, Cao11}. 
\end{rem}

In the following, we will address the most important contribution that our ensemble DDM algorithm has a mesh-independent convergence rate with some suitable choice of parameters. We first introduce one necessary and important lemma as follows.

\begin{lemma}\label{abc}
	Suppose $a_1,a_2,b_1,b_2,c_1,c_2$ are positive constants, with $a_2 < a_1,\ b_2 < b_1, \ c_2 < c_1$, and $A^n,\ B^n, \ C^n, \ n=1,2,\cdots$ are three different iterative sequences. If $a_1 A^n + b_1 B^n + c_1 C^n \leq a_2 A^{n-1} + b_2 B^{n-1} + c_2 c^{n-1}$, the following estimate holds:
	\begin{eqnarray*}
		&&a_1 A^n + b_1 B^n + c_1 C^n \leq
		\max\Big\{\frac{a_2}{a_1},\frac{b_2}{b_1},\frac{c_2}{c_1} \Big\}^{n-1}\Big(a_2 A^{0} + b_2 B^{0} + c_2 C^{0}\Big).
	\end{eqnarray*}
\end{lemma}
\begin{proof}
Without loss of generality, assuming $\frac{a_2}{a_1}$ is the largest one in $\Big\{\frac{a_2}{a_1},\frac{b_2}{b_1},\frac{c_2}{c_1} \Big\}$. Then we have
\begin{eqnarray*}
 	&&a_1 A^n + b_1 B^n + c_1 C^n \leq\frac{a_2}{a_1} \Big( a_1 A^{n-1} + b_1 B^{n-1} + c_1 c^{n-1} \Big) 
 	+ \Big(b_2-\frac{a_2}{a_1}b_1 \Big) B^{n-1}
 	+ \Big(c_2-\frac{a_2}{a_1}c_1 \Big) C^{n-1}.
\end{eqnarray*}
Since $\frac{a_2}{a_1}$ is the largest one, we know $\frac{a_2}{a_1}\geq \frac{b_2}{b_1}$ and $\frac{a_2}{a_1}\geq \frac{c_2}{c_1}$, which can inform that $b_2-\frac{a_2}{a_1}b_1\leq 0$ and  $c_2-\frac{a_2}{a_1}c_1\leq 0$. Then by recurrence we can derive
\begin{eqnarray*}
	a_1 A^n + b_1 B^n + c_1 C^n \leq  \frac{a_2}{a_1} \Big( a_1 A^{n-1} + b_1 B^{n-1} + c_1 C^{n-1} \Big) 
	&\leq& \frac{a_2}{a_1} \Big( a_2 A^{n-2} + b_2 B^{n-2} + c_2 C^{n-2} \Big)\\ 
	\leq\Big( \frac{a_2}{a_1}\Big)^2 \Big( a_1 A^{n-2} + b_1 B^{n-2} + c_1 C^{n-2} \Big) &\leq& \cdots \leq \Big(\frac{a_2}{a_1} \Big)^{n-1}\Big(a_2 A^{0} + b_2 B^{0} + c_2 C^{0}\Big),
\end{eqnarray*}
which completes the lemma.
\end{proof}

Now, we present the geometric convergence of ensemble DDM for the case $\delta_S<\delta_D$ as follows.
\begin{theorem}\label{ItrDDMf<p}
	For the case of $\delta_S<\delta_D$, if selecting the parameters $\delta_S$ and $\delta_D$ as
	\begin{eqnarray}\label{controlTP}
	&&\frac{2(C_{\mathrm{tr}}C'_{\mathrm{tr}}\xi_j)^2}{C_1\nu}-g(\overline{k}^{\min}-E_{k_j}^{\max})<\delta_D-\delta_S<\frac{C_1\nu}{C_4}, \hspace{8mm} 
		\rho(\delta_S,\delta_D)<1.
	\end{eqnarray}
where $\rho(\delta_S,\delta_D):=\frac{(1+C_3)^2g(\delta_D-\delta_S)[2(\overline{k}^{\max})^2+(E_{k_{j}}^{\max})^2]+2\delta_D^2(\overline{k}^{\min}-E_{k_{j}}^{\max})}{4(\delta_D^2-\delta_S^2)(\overline{k}^{\min}-E_{k_{j}}^{\max})}$, and assuming that 
\begin{eqnarray}\label{perturbation}
	&&E_{{\xi}_{i,j}}^{\mathrm{max}} < \overline{\xi}_j, \hspace{10mm}
	E_{k_j}^{\mathrm{max}} < \overline{k}^{\mathrm{min}},
\end{eqnarray}
which means that the random hydraulic conductivity tensor satisfies a small disturbance. For a given positive constant $C^*$,
the Ensemble DDM algorithm has the following convergence
	\begin{eqnarray}\label{result}
		&& ||\mathbf{e}_{Sj}^N||_1^2+||\varepsilon_{Sj}^{N}||^2
        +||\mathbf{e}_{Dj}^{N-1}||_{\mathrm{div}}^2 + ||\varepsilon_{Dj}^{N}||^2+ ||\eta^{N}_{Sj}||_{\Gamma}^2
          \nonumber\\
       &&\hspace{12mm}
	+ ||\eta^{N}_{Sj,\tau}||_{H^{-\frac{1}{2}}(\Gamma)}^2
       +||\eta^{N}_{Dj}||_{\Gamma}^2
        +\sqrt{\rho(\delta_S,\delta_D)}||\eta^{N-1}_{Dj}||_{\Gamma}^2 \nonumber \\
		&&\leq C^*	\max\Big\{ \sqrt{\rho(\delta_S,\delta_D)},\  \sum_{i=1}^{d-1}\frac{E_{{\xi}_{i,j}}^{\mathrm{max}}}{2\overline{\xi}_i-E_{{\xi}_{i,j}}^{\mathrm{max}}},\ \frac{2g\overline{k}^{\min}}{3g\overline{k}^{\min}-gE_{k_j}^{\max}+\delta_D-\delta_S-\frac{2(C_{\mathrm{tr}}C'_{\mathrm{tr}}\xi_j)^2}{C_1\nu}} \Big\}^{N-2}\nonumber\\
		&&\hspace{2mm} \Big[ \sqrt{\rho(\delta_S,\delta_D)}\Big(||\eta_{Dj}^{1}||^2_{\Gamma}+\sqrt{\rho(\delta_S,\delta_D)}||\eta_{Dj}^{0}||^2_{\Gamma}\Big)\nonumber\\ 
		&&\hspace{10mm}+(\delta_S+\delta_D) \sum_{i=1}^{d-1}E_{{\xi}_{i,j}}^{\max}
		||\mathbf{e}_{Sj}^{1}\cdot \tau_{i}||_{\Gamma}^2+(\delta_S+\delta_D)g\overline{k}^{\min}||\mathbf{e}_{Dj}^{0}||_{\mathrm{div}}^2 \Big].
	\end{eqnarray}
\end{theorem}

\begin{proof}
Firstly, the test function $\mathbf{v}_D$ with the properties
\begin{eqnarray}\label{etacontrolv}
\mathbf{v}_D\cdot\mathbf{n}_D|_{\Gamma}=\eta_{Dj}^{n-2}, \hspace{10mm} ||\mathbf{v}_D||_{\mathrm{div}}\leq||\eta_{Dj}^{n-2}||_{\Gamma}.
\end{eqnarray}
can be constructed, see for instance (3.34) therein of \cite{Sun2021}.

Substituting such $\mathbf{v}_D$ into (\ref{err-3}), and using the Cauchy-Schwarz inequality, Young's inequality and (\ref{PDcontrol}), we can get
\begin{eqnarray}\label{contral_eta}
||\eta_{Dj}^{n-2}||_{\Gamma}^{2}&=&
-\overline{a_{D}}(\mathbf{e}_{Dj}^{n-1},\mathbf{v}_{D})+b_D({\varepsilon}_{Dj}^{n-1},\mathbf{v}_D)
-\delta_D\langle\mathbf{e}_{Dj}^{n-1}\cdot\mathbf{n}_D,\eta_{Dj}^{n-2}\rangle\nonumber\\
&&+g((\mathbb{K}_j^{-1}-\overline{\mathbb{K}})\mathbf{e}_{Dj}^{n-2},\mathbf{v}_{D})_{D}
+g(k_j^{\mathrm{min}}-\overline{k}^{\mathrm{min}})(\mathrm{div}\hspace{0.5mm}\mathbf{e}_{Dj}^{n-2},\mathrm{div}\hspace{0.5mm}\mathbf{v}_{D})_{D}\nonumber\\
&\leq&g\overline{k}^{\max}||\mathbf{e}_{Dj}^{n-1}||_{\mathrm{div}}||\mathbf{v}_D||_{\mathrm{div}}+g||\varepsilon_{Dj}^{n-1}||_D||\nabla\cdot\mathbf{v}_D||_D
\nonumber\\
&&
+\delta_D||\mathbf{e}_{Dj}^{n-1}\cdot\mathbf{n}_D||_{\Gamma}||\eta_{Dj}^{n-2}||_{\Gamma}+gE_{k_j}^{\max}||\mathbf{e}_{Dj}^{n-2}||_{\mathrm{div}}||\mathbf{v}_D||_{\mathrm{div}}
\nonumber\\
&\leq&
(1+C_3)g\overline{k}^{\max}||\mathbf{e}_{Dj}^{n-1}||_{\mathrm{div}}||\eta_{Dj}^{n-2}||_{\Gamma}+(1+C_3)gE_{k_j}^{\max}||\mathbf{e}_{Dj}^{n-2}||_{\mathrm{div}}||\eta_{Dj}^{n-2}||_{\Gamma}\nonumber\\
&&+\delta_D||\mathbf{e}_{Dj}^{n-1}\cdot\mathbf{n}_D||_{\Gamma}||\eta_{Dj}^{n-2}||_{\Gamma}\nonumber\\
&\leq& \frac{1}{2}(\delta_S+\delta_D)g(\overline{k}^{\min}-E_{k_j}^{\max})||\mathbf{e}_{Dj}^{n-1}||_{\mathrm{div}}^2+\frac{\delta_D^2-\delta_S^2}{2}||\mathbf{e}_{Dj}^{n-1}\cdot\mathbf{n}_{D}||^2_{\Gamma}\nonumber\\
&&+(\delta_S+\delta_D)g(\overline{k}^{\min}-E_{k_j}^{\max})||\mathbf{e}_{Dj}^{n-2}||_{\mathrm{div}}^2+\rho(\delta_S,\delta_D)||\eta_{Dj}^{n-2}||_{\Gamma}^2,
\end{eqnarray}
where 
\begin{eqnarray}\label{rhoSD}
\rho(\delta_S,\delta_D):=\frac{(1+C_3)^2g(\delta_D-\delta_S)[2(\overline{k}^{\max})^2+(E_{k_j}^{\max})^2]+2\delta_D^2(\overline{k}^{\min}-E_{k_j}^{\max})}{4(\delta_D^2-\delta_S^2)(\overline{k}^{\min}-E_{k_j}^{\max})}.
\end{eqnarray}
In the equation (\ref{errorD}), replacing the term $||\eta_{Sj}^{n-1}||_{\Gamma}^2$ by the equation (\ref{errorS}), applying the estimates (\ref{bilinear-ieq})-(\ref{term3}) and (\ref{contral_eta}), using the Poincar$\acute{\mathrm{e}}$, Korn's and Young's inequalities, we can conclude
\begin{eqnarray*}
||\eta_{Dj}^{n}||^2_{\Gamma}
&=& ||\eta_{Dj}^{n-2}||^2_{\Gamma}+(\delta_D^2-\delta_S^2)||\mathbf{e}_{Sj}^{n}\cdot\mathbf{n}_{S}||^2_{\Gamma}+(\delta_S^2-\delta_D^2)||\mathbf{e}_{Dj}^{n-1}\cdot\mathbf{n}_{D}||^2_{\Gamma}\\
&&-2(\delta_S+\delta_D)\Big[ a_S(\mathbf{e}_{Sj}^n,\mathbf{e}_{Sj}^n)+ \overline{a_D}(\mathbf{e}_{Dj}^{n-1},\mathbf{e}_{Dj}^{n-1})+\sum_{i=1}^{d-1}\overline{\xi}_i
||\mathbf{e}_{Sj}^n\cdot \tau_{i}||_{\Gamma}^2\\
&&\hspace{20mm} -\sum_{i=1}^{d-1}{\xi}_{i,j}
\langle\mathbf{e}_{Dj}^{n-1}\cdot \tau_{i},\mathbf{e}_{Sj}^n\cdot \tau_{i}\rangle+\sum_{i=1}^{d-1}({\xi}_{i,j}-\overline{\xi}_i)
\langle\mathbf{e}_{Sj}^{n-1}\cdot \tau_{i},\mathbf{e}_{Sj}^n\cdot \tau_{i}\rangle \\
&&\hspace{20mm}+g((\mathbb{K}_j^{-1}-\overline{\mathbb{K}})\mathbf{e}_{Dj}^{n-2},\mathbf{e}_{Dj}^{n-1})_D
+g(k_j^{\mathrm{min}}-\overline{k}^{\mathrm{min}})(\mathrm{div}\hspace{0.5mm}\mathbf{e}_{Dj}^{n-2},\mathrm{div}\hspace{0.5mm}\mathbf{e}_{Dj}^{n-1})_{D}\Big]\\
&\leq&||\eta_{Dj}^{n-2}||^2_{\Gamma}+\frac{\delta_S^2-\delta_D^2}{2}||\mathbf{e}_{Dj}^{n-1}\cdot\mathbf{n}_{D}||^2_{\Gamma}-2(\delta_S+\delta_D)\\
&&\Big[ \Big(\frac{C_1\nu}{2}-\frac{C_4(\delta_D-\delta_S)}{2}\Big)||\mathbf{e}_{Sj}^{n}||_{1}^2+\sum_{i=1}^{d-1}\frac{2\overline{\xi}_i-E_{{\xi}_{i,j}}^{\max}}{2}
||\mathbf{e}_{Sj}^n\cdot \tau_{i}||_{\Gamma}^2 -\sum_{i=1}^{d-1}\frac{E_{{\xi}_{i,j}}^{\max}}{2}
||\mathbf{e}_{Sj}^{n-1}\cdot \tau_{i}||_{\Gamma}^2 \\ 
&&\hspace{2mm}+\Big(g\overline{k}^{\min}-\frac{gE_{k_j}^{\max}}{2}+\frac{\delta_D-\delta_S}{4}-\frac{(C_{\mathrm{tr}}C'_{\mathrm{tr}}\xi_j)^2}{2C_1\nu} \Big) ||\mathbf{e}_{Dj}^{n-1}||_{\mathrm{div}}^2-\frac{gE_{k_j}^{\max}}{2}||\mathbf{e}_{Dj}^{n-2}||_{\mathrm{div}}^2\Big] \\
&\leq&\rho(\delta_S,\delta_D)||\eta_{Dj}^{n-2}||^2_{\Gamma}-2(\delta_S+\delta_D)\Big[ \Big(\frac{C_1\nu}{2}-\frac{C_4(\delta_D-\delta_S)}{2}\Big)||\mathbf{e}_{Sj}^{n}||_{1}^2\\
&&\hspace{6mm}+\sum_{i=1}^{d-1}\frac{2\overline{\xi}_i-E_{{\xi}_{i,j}}^{\max}}{2}
||\mathbf{e}_{Sj}^n\cdot \tau_{i}||_{\Gamma}^2-\sum_{i=1}^{d-1}\frac{E_{{\xi}_{i,j}}^{\max}}{2}
||\mathbf{e}_{Sj}^{n-1}\cdot \tau_{i}||_{\Gamma}^2 \\ 
&&\hspace{6mm}+\Big(\frac{3g\overline{k}^{\min}}{4}-\frac{gE_{k_j}^{\max}}{4}+\frac{\delta_D-\delta_S}{4}-\frac{(C_{\mathrm{tr}}C'_{\mathrm{tr}}\xi_j)^2}{2C_1\nu} \Big) ||\mathbf{e}_{Dj}^{n-1}||_{\mathrm{div}}^2 
-\frac{g\overline{k}^{\min}}{2}||\mathbf{e}_{Dj}^{n-2}||_{\mathrm{div}}^2\Big].
\end{eqnarray*}
Then if improving the additional condition (\ref{add_condition1}) as:
\begin{eqnarray}\label{add_condition2}
\frac{2(C_{\mathrm{tr}}C'_{\mathrm{tr}}\xi_j)^2}{C_1\nu}-g(\overline{k}^{\min}-E_{k_j}^{\max})<\delta_D-\delta_S<\frac{C_1\nu}{C_4},
\end{eqnarray}
so that $\frac{3g\overline{k}^{\min}}{4}-\frac{gE_{k_j}^{\max}}{4}+\frac{\delta_D-\delta_S}{4}-\frac{(C_{\mathrm{tr}}C'_{\mathrm{tr}}\xi_j)^2}{2C_1\nu}>\frac{g\overline{k}^{\min}}{2}$.
To this end, assuming that 
\begin{eqnarray*}
\rho(\delta_S,\delta_D)<1,\hspace{6mm} \overline{k}^{\min}>E_{k_j}^{\max},\hspace{6mm}
\overline{\xi}_i>E_{{\xi}_{i,j}}^{\max},
\end{eqnarray*}
we can derive
\begin{eqnarray*}
	&&\hspace{-6mm}2(\delta_S+\delta_D)\Big[\sum_{i=1}^{d-1}\frac{2\overline{\xi}_i-E_{{\xi}_{i,j}}^{\max}}{2}
	||\mathbf{e}_{Sj}^n\cdot \tau_{i}||_{\Gamma}^2
	+\Big(\frac{3g\overline{k}^{\min}-gE_{k_j}^{\max}+\delta_D-\delta_S}{4}-\frac{(C_{\mathrm{tr}}C'_{\mathrm{tr}}\xi_j)^2}{2C_1\nu} \Big) ||\mathbf{e}_{Dj}^{n-1}||_{\mathrm{div}}^2\Big]\\
	&&\hspace{-4mm}+||\eta_{Dj}^{n}||^2_{\Gamma}\leq\rho(\delta_S,\delta_D)||\eta_{Dj}^{n-2}||^2_{\Gamma}+2(\delta_S+\delta_D) \sum_{i=1}^{d-1}\frac{E_{{\xi}_{i,j}}^{\max}}{2}
	||\mathbf{e}_{Sj}^{n-1}\cdot \tau_{i}||_{\Gamma}^2 +2(\delta_S+\delta_D)\frac{g\overline{k}^{\min}}{2}||\mathbf{e}_{Dj}^{n-2}||_{\mathrm{div}}^2,
\end{eqnarray*}
which can be technically transformed into 
\begin{eqnarray*}
	&&\Big(||\eta_{Dj}^{n}||^2_{\Gamma}+\sqrt{\rho(\delta_S,\delta_D)}||\eta_{Dj}^{n-1}||^2_{\Gamma}\Big)+2(\delta_S+\delta_D)\sum_{i=1}^{d-1}\frac{2\overline{\xi}_i-E_{{\xi}_{i,j}}^{\max}}{2}
	||\mathbf{e}_{Sj}^n\cdot \tau_{i}||_{\Gamma}^2\\
	&&\hspace{4.5mm}+2(\delta_S+\delta_D)\Big(\frac{3g\overline{k}^{\min}-gE_{k_j}^{\max}+\delta_D-\delta_S}{4}-\frac{(C_{\mathrm{tr}}C'_{\mathrm{tr}}\xi_j)^2}{2C_1\nu} \Big) ||\mathbf{e}_{Dj}^{n-1}||_{\mathrm{div}}^2\\
	&&\leq\sqrt{\rho(\delta_S,\delta_D)}\Big(||\eta_{Dj}^{n-1}||^2_{\Gamma}+\sqrt{\rho(\delta_S,\delta_D)}||\eta_{Dj}^{n-2}||^2_{\Gamma}\Big)\\
	&&\hspace{4.5mm}+2(\delta_S+\delta_D)\Big[ \sum_{i=1}^{d-1}\frac{E_{{\xi}_{i,j}}^{\max}}{2}
	||\mathbf{e}_{Sj}^{n-1}\cdot \tau_{i}||_{\Gamma}^2+\frac{g\overline{k}^{\min}}{2}||\mathbf{e}_{Dj}^{n-2}||_{\mathrm{div}}^2\Big].
\end{eqnarray*}
Finally, by utilizing Lemma \ref{abc}, we arrive at
\begin{eqnarray*}
&&\Big(||\eta_{Dj}^{N}||^2_{\Gamma}+\sqrt{\rho(\delta_S,\delta_D)}||\eta_{Dj}^{N-1}||^2_{\Gamma}\Big)+2(\delta_S+\delta_D)\sum_{i=1}^{d-1}\frac{2\overline{\xi}_i-E_{{\xi}_{i,j}}^{\max}}{2}
||\mathbf{e}_{Sj}^N\cdot \tau_{i}||_{\Gamma}^2\\
&&\hspace{14.5mm}+2(\delta_S+\delta_D)\Big(\frac{3g\overline{k}^{\min}-gE_{k_j}^{\max}+\delta_D-\delta_S}{4}-\frac{(C_{\mathrm{tr}}C'_{\mathrm{tr}}\xi_j)^2}{2C_1\nu} \Big) ||\mathbf{e}_{Dj}^{N-1}||_{\mathrm{div}}^2\\
&&\leq	\max\Big\{ \sqrt{\rho(\delta_S,\delta_D)},\  \sum_{i=1}^{d-1}\frac{E_{{\xi}_{i,j}}^{\max}}{2\overline{\xi}_i-E_{{\xi}_{i,j}}^{\max}},\ \frac{2g\overline{k}^{\min}}{3g\overline{k}^{\min}-gE_{k_j}^{\max}+\delta_D-\delta_S-\frac{2(C_{\mathrm{tr}}C'_{\mathrm{tr}}\xi_j)^2}{C_1\nu}} \Big\}^{N-2}\\
&& \Big\{ \sqrt{\rho(\delta_S,\delta_D)}\Big(||\eta_{Dj}^{1}||^2_{\Gamma}+\sqrt{\rho(\delta_S,\delta_D)}||\eta_{Dj}^{0}||^2_{\Gamma}\Big)+(\delta_S+\delta_D)\Big[ \sum_{i=1}^{d-1}E_{{\xi}_{i,j}}^{\max}
||\mathbf{e}_{Sj}^{1}\cdot \tau_{i}||_{\Gamma}^2+g\overline{k}^{\min}||\mathbf{e}_{Dj}^{0}||_{\mathrm{div}}^2 \Big]\Big\},
\end{eqnarray*}
which shows that $ ||\eta_{Dj}^{n}||_{\Gamma}^{2}, \ ||\mathbf{e}_{Sj}^{n}\cdot \tau_{i}||_{\Gamma}^2$, and $ ||\mathbf{e}_{Dj}^{n}||^2_{\mathrm{div}}$ are geometrically convergent.

By (\ref{errorS}), we can further get
\begin{eqnarray*}
	||\eta_{Sj}^{n}||^2_{\Gamma}&\leq&
	||\eta_{Dj}^{n-1}||^2_{\Gamma}
	-(\delta_S+\delta_D)(2g\overline{k}^{\min}-gE_{k_j}^{\max})||\mathbf{e}_{Dj}^n||_{\mathrm{div}}^2 +(\delta_S+\delta_D)gE_{k_j}^{\max}||\mathbf{e}_{Dj}^{n-1}||_{\mathrm{div}}^2\\
	&\leq&||\eta_{Dj}^{n-1}||^2_{\Gamma}
	+(\delta_S+\delta_D)gE_{k_j}^{\max}||\mathbf{e}_{Dj}^{n-1}||_{\mathrm{div}}^2.
\end{eqnarray*}
Hence $||\eta_{Sj}^{n}||^2_{\Gamma}$ also has geometric convergence.

Moreover, the geometric convergence of
$ ||\eta_{Sj}^{n}||^2_{\Gamma}$, and $ ||\mathbf{e}_{Sj}^{n}\cdot \tau_{i}||_{\Gamma}^2$ together with the error equation (\ref{aS+bS}) implies the geometric convergence of $||\mathbf{e}_{Sj}^{n}||_1^2$.
Combining it with the results of (\ref{contralp}), (\ref{PDcontrol}), (\ref{etaStau-control}), we can summarize all the geometric convergence results in (\ref{result}).
\end{proof}

\begin{rem}
	\begin{figure}[htbp] 
		\centering
		\includegraphics[width=53mm,height=40mm]{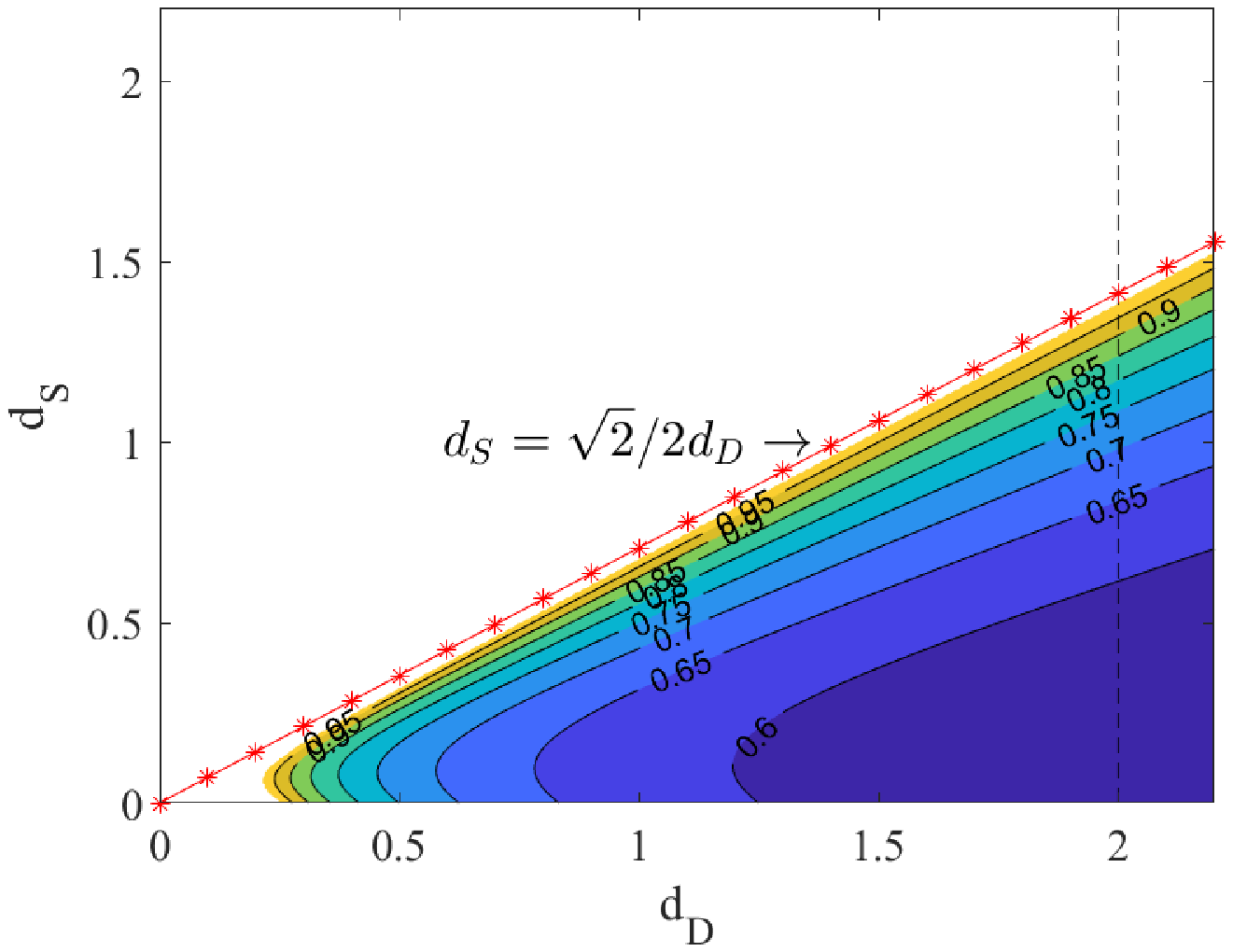} \hspace{3mm}
		\includegraphics[width=53mm,height=40mm]{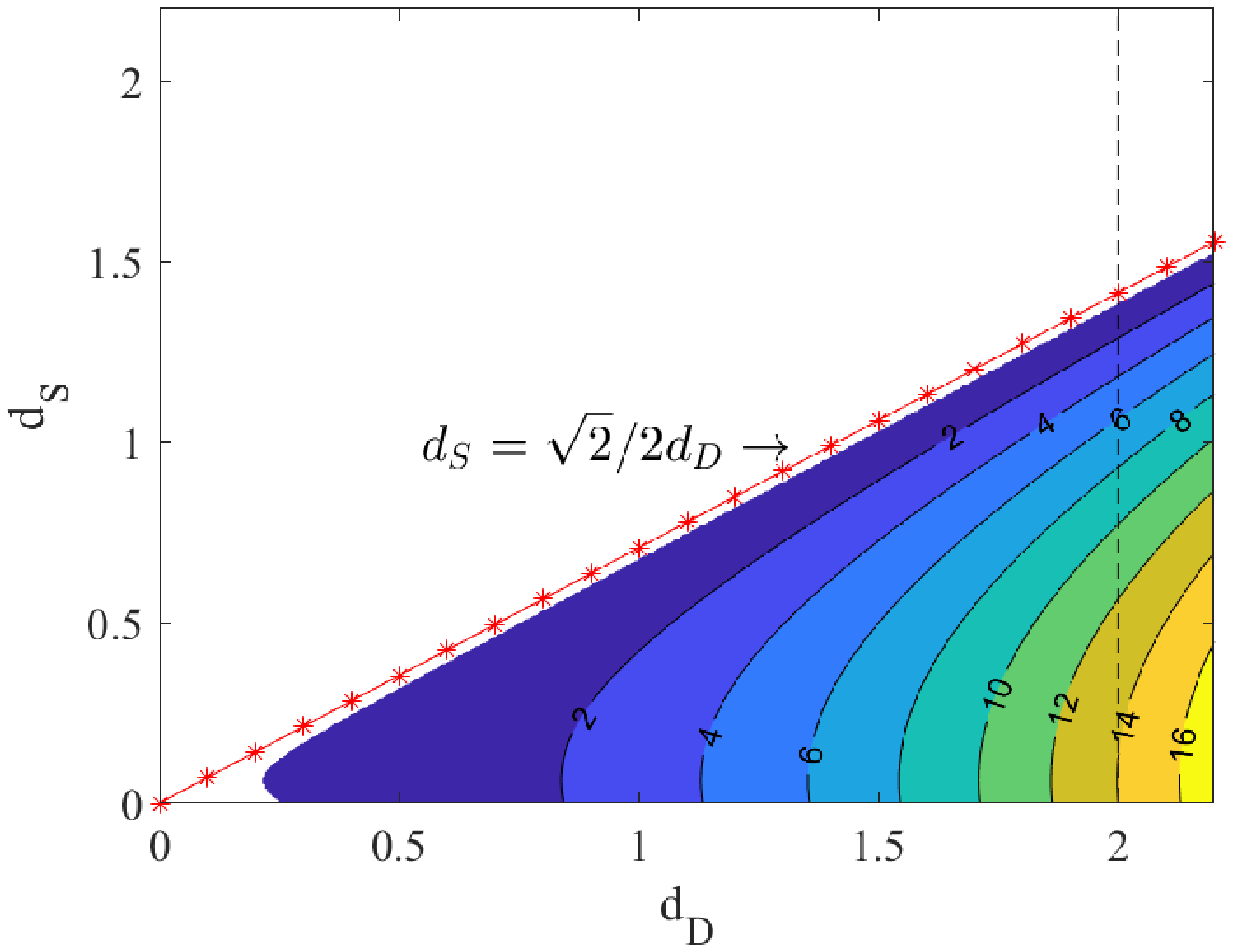}
		\caption{The contour lines of the functions $\hat{\rho}(d_D,d_S)<1$ (left)
			and $\bar{\rho}(d_D,d_S)>0$ (right).}\label{contour}
	\end{figure}
	To show the possible choice of $\rho(\delta_S,\delta_D)<1$, we can utilize a similar argument as Remark 3.4 in \cite{Sun2021} and select	$\delta_D=\frac{2(1+C_3)^2g[2(\overline{k}^{\max})^2+(E_{k_j}^{\max})^2]}{(\overline{k}^{\min}-E_{k_j}^{\max})}d_D, \ \delta_S=\frac{2(1+C_3)^2g[2(\overline{k}^{\max})^2+(E_{k_j}^{\max})^2]}{(\overline{k}^{\min}-E_{k_j}^{\max})}d_S$,
    and also define $\hat{\rho}(d_S,d_D):= \frac{(d_D-d_S)+4d_D^2}{8(d_D^2-d_S^2)}=\rho(\delta_S,\delta_D)$. It is much easier to check that $\rho(\delta_S,\delta_D)=\hat{\rho}(d_S,d_D)<1$ is equivalent to $ \bar{\rho}(d_S,d_D):=4d_D^2-d_D - (8d_S^2-d_S)>0$.
    To this end, for showing the suitable selection of $d_S$ and $d_D$, we display the contours of the functions $\hat{\rho}$ and $\bar{\rho}$ in Fig. \ref{contour} to find the evolution of $\hat{\rho}(d_S,d_D)<0$ and $\bar{\rho}(d_S,d_D)>0$ respectively. In addition, compared with Fig. 3.1 in \cite{Sun2021}, we can see that the present limitations on the choice of the parameters for the ensemble DDM are stricter.
\end{rem}

\section{Optimized Robin Parameters}
For Robin-type transmission conditions, it is of great interest to find optimized Robin parameters to accelerate the convergence of our proposed iterative algorithms. The optimized Schwarz methods have been widely studied  for different PDEs and recently for the Stokes-Darcy model in \cite{Discacciati, Discacciati18, Gander, Liu22}. In this section, we will utilize the Fourier transform in the direction tangential to the interface (corresponding to the $y$-variable in the case studied here), which mainly follows the analysis technique of \cite{Discacciati}. Hence, the same problem domain and hypotheses as used in \cite{Discacciati} are considered, namely, the fluid region is the half-plane $\Omega_{S}=\{(x,y)\in \mathbb{R}^2: x<0 \}$, and the porous media is another half-plane $\Omega_{S}=\{(x,y)\in \mathbb{R}^2: x>0 \}$, separated by the interface $\Gamma=\{(x,y)\in \mathbb{R}^2: x=0 \}$. For the sake of the clearance of the analysis, we similarly assume that $g=1, z=0,$ and $\mathbb{K}_{j}(x)=\hat{k}_j\mathbb{I}$, where $\hat{k}_j$ is a positive constant. Then, we can define $\overline{\hat{k}}=\frac{1}{J}\sum_{j=1}^{J} \hat{k}_j^{-1}$. Under such circumstances, given $\mathbf{e}_{Sj}^0, \varepsilon_{Sj}^0, \mathbf{e}_{Dj}^0 $ and $\varepsilon_{Dj}^0 $, the Stokes-Darcy error equations by the idea of the ensemble domain decomposition can be expressed as:
\begin{eqnarray}
	&&-\nabla \cdot \mathbb{T}(\mathbf{e}_{Sj}^n,\varepsilon_{Sj}^n) = 0 \hspace{67mm} \mathrm{in}~(-\infty,0) \times \mathbb{R},
	\label{OpStokes1} \\
	&&\hspace{14.4mm}\nabla \cdot \mathbf{e}_{Sj}^n = 0\ \ \hspace{61.3mm} \ \ \ \mathrm{in}~(-\infty,0) \times \mathbb{R},\label{OpStokes2}\\
	&&-\mathbf{n}_{S}\cdot (\mathbb{T}(\mathbf{e}_{Sj}^n,\varepsilon_{Sj}^n)\cdot \mathbf{n}_{S})-\delta_S \mathbf{e}_{Sj}^n\cdot \mathbf{n}_S
	=  \varepsilon_{Dj}^{n-1} + \delta_S \mathbf{e}_{Dj}^{n-1} \cdot \mathbf{n}_D \hspace{8.3mm} \mathrm{on} ~ \{0\}\times \mathbb{R},\label{OpInterface1}\\
	&&-\tau_{i}\cdot (\mathbb{T}(\mathbf{e}_{Sj}^n,\varepsilon_{Sj}^n)\cdot \mathbf{n}_{S})
	=\frac{\alpha}{\sqrt{\tau_{i}\cdot\overline{\mathbb{K}}\tau_{i}}} \tau_i \mathbf{e}_{Sj}^n+\Big(\frac{\alpha}{\sqrt{\tau_{i}\cdot\mathbb{K}_j\tau_{i}}}-\frac{\alpha}{\sqrt{\tau_{i}\cdot\overline{\mathbb{K}}\tau_{i}}}\Big) \tau _{i} \cdot (\mathbf{e}_{Sj}^{n-1}-\mathbf{e}_{Dj}^{n-1})\nonumber\\
	&&\hspace{74.3mm}  1\leq i \leq d-1 \hspace{4mm} \mathrm{on} \{0\}\times \mathbb{R},\label{OpInterface2}\\
	&&\overline{\mathbb{K}} \mathbf{e}_{Dj}^n+(\mathbb{K}_j^{-1}-\overline{\mathbb{K}}) \mathbf{e}_{Dj}^{n-1}= -\nabla \varepsilon_{Dj}^n\ \ \ \ \ \ \hspace{39.3mm} 
	\mathrm{in}~(0,\infty) \times \mathbb{R} ,  \label{OpDarcy1} \\
	&&\hspace{25mm}\nabla \cdot \mathbf{e}_{Dj}^{n} = 0\ \ \ \ \ \ \ \ \ \hspace{45.5mm}
	\mathrm{in}~(0,\infty) \times \mathbb{R}. \label{OpDarcy2}\\
	&&\varepsilon_{Dj}^{n} - \delta_D \mathbf{e}_{Dj}^{n} \cdot \mathbf{n}_D =-\mathbf{n}_{S}\cdot (\mathbb{T}(\mathbf{e}_{Sj}^{n-1},\varepsilon_{Sj}^{n-1})\cdot \mathbf{n}_{S})+\delta_D \mathbf{e}_{Sj}^{n-1}\cdot \mathbf{n}_S
	\hspace{4.0mm} \mathrm{on} ~ \{0\}\times \mathbb{R},\label{OpInterface3}
\end{eqnarray}
Then, for $w(x,y)\in L^2(\mathbb{R}^2)$, we can define the Fourier transform as: { $\hat{\mathcal{F}}: w(x,y) \mapsto \hat{w}(x,m)=\int_{\mathbb{R}} e^{-imy}w(x,y)  \mathrm{d} y $}, where $m$ is the frequency variable. We can explicitly characterize the convergence factor of the proposed ensemble domain decomposition algorithm in the following theorem.

\begin{thm}
	Given $( \mathbf{e}_{Sj}^0, \ \varepsilon_{Sj}^0, \ \mathbf{e}_{Dj}^0, \ \varepsilon_{Dj}^0 )$, the convergence factor of the ensemble domain decomposition algorithm does not depend on the iteration, which is precisely given by
	\begin{eqnarray}\label{Covfactor}
		\rho(\delta_S, \delta_D, m) =\left| \frac{2\nu |m|- \delta_D}{2\nu |m| + \delta_S} \right|.
	\end{eqnarray}
\end{thm}

\begin{proof}
	Following the analysis of Proposition 3.1 in \cite{Discacciati}, we can utilize the same arguments to obtain the results for the Stokes error equations (\ref{OpStokes1})-(\ref{OpInterface2}) and the Darcy error equations (\ref{OpDarcy1})-(\ref{OpInterface3}). For any $m$, the pressure error $\hat{\varepsilon}_{Sj}^{n}$ and the $x$-direction component of the velocity error  $\hat{e}_{Sj1}^n$ of the Stokes problem (\ref{OpStokes1})-(\ref{OpStokes2}) satisfy
	\begin{eqnarray}\label{ResultSpu}
		\hat{\varepsilon}_{Sj}^n(x,m)=P^n(m) e^{|m|x}, \hspace{6mm}
		\hat{e}_{Sj1}^n(x,m)=\Big( A^n(m)+\frac{x}{2\nu}P^n(m) \Big) e^{|m|x},
	\end{eqnarray}
	where the value of $P^n(m)$ and $A^n(m)$ are only determined by the interface condition (\ref{OpInterface1}). The piezometric head error $\hat{\varepsilon}_{Dj}$ of the Darcy equation  (\ref{OpDarcy1})-(\ref{OpDarcy2}) is given by
	\begin{eqnarray}\label{ResultDp}
		\hat{\varepsilon}_{Dj}^n(x,m)=Q^n(m) e^{-|m|x},
	\end{eqnarray}
	where the value of $Q^n(m)$ is determined uniquely by the interface condition (\ref{OpInterface3}). Since the equation (\ref{OpDarcy1}) is linear, we can assume the $x$-direction component  of the velocity error $\hat{e}_{Dj1}^n$ of the Darcy system as
	\begin{eqnarray}\label{ResultDu}
		\hat{e}_{Dj1}^n(x,m)=B^n(m) e^{-|m|x}.
	\end{eqnarray}
	This definition of $\hat{e}_{Dj1}^n$ is the unique solution that satisfies equation (\ref{OpDarcy1}), where the value of $B^n(m)$ can be obtained by the interface condition (\ref{OpInterface3}). Substituting (\ref{ResultDp})-(\ref{ResultDu}) into (\ref{OpDarcy1}), then letting $x=0$
 we can get
	\begin{eqnarray}\label{Iter1}
		\overline{\hat{k}}B^n(m)-(\overline{\hat{k}}-\hat{k}_j^{-1})B^{n-1}(m)=|m|Q^{n}(m).
	\end{eqnarray}
	Inserting (\ref{ResultSpu})-(\ref{ResultDu}) into (\ref{OpInterface1}) and (\ref{OpInterface3}), we have
	\begin{eqnarray}\label{Iter2}
		-(2\nu|m|+\delta_S)A^n(m)&=&Q^{n-1}(m)-\delta_S B^{n-1}(m),\\ 
		\label{Iter3}
		Q^n(m)+\delta_D B^n(m) &=&(\delta_D-2\nu |m|)A^{n-1}(m).
	\end{eqnarray}
	To display our analysis much clearer and instructively, we omit the variable $m$ for $A$ and $B$ with suitable superscripts temporarily. We can combine the three equations (\ref{Iter1})-(\ref{Iter3}) for purpose of removing $Q^{n}$ or $Q^{n-1}$ to obtain\begin{eqnarray}\label{Iter4}
		&&(\overline{\hat{k}}+\delta_D |m|) B^n=(\overline{\hat{k}}-\hat{k}_j^{-1})B^{n-1}+|m|(\delta_D-2\nu |m|)A^{n-1},\\ 
		\label{Iter5}
		&&(\overline{\hat{k}}-\delta_S |m|) B^n+|m|(2\nu|m| +\delta_S)A^{n+1}=(\overline{\hat{k}}-\hat{k}_j^{-1})B^{n-1},\\
		\label{Iter6}
		&&(\delta_S+\delta_D)B^{n}-(2\nu |m|+\delta_S) A^{n+1}= (\delta_D-2\nu |m|) A^{n-1}.
	\end{eqnarray}
	Replacing the term $A^{n-1}$ in equation (\ref{Iter4}) by equation (\ref{Iter5}), we can derive
	\begin{eqnarray*}
		(\overline{\hat{k}}+\delta_D |m|) B^n-(\overline{\hat{k}}-\hat{k}_j^{-1})B^{n-1}
		=\frac{2\nu |m|-\delta_D}{2\nu |m|+\delta_S} \Big[ (\overline{\hat{k}}-\delta_S |m|) B^{n-2}-(\overline{\hat{k}}-\hat{k}_j^{-1})B^{n-3} \Big],
	\end{eqnarray*}
	namely,
	\begin{eqnarray*}
		B^n-\frac{\overline{\hat{k}}-\hat{k}_j^{-1}}{\overline{\hat{k}}+\delta_D |m|}B^{n-1}
		=\frac{2\nu |m|-\delta_D}{2\nu |m|+\delta_S} \Big[B^{n-2}-\frac{\overline{\hat{k}}-\hat{k}_j^{-1}}{\overline{\hat{k}}+\delta_D |m|}B^{n-3} \Big]
		- \frac{2\nu |m|-\delta_D}{2\nu |m|+\delta_S} \frac{(\delta_S+\delta_D) |m|}{\overline{\hat{k}}+\delta_D |m|} B^{n-2}.
	\end{eqnarray*}
	Then, we can write equation (\ref{Iter6}) as
	\begin{eqnarray*}
		-\frac{2\nu |m|-\delta_D}{2\nu |m|+\delta_S} \frac{(\delta_S+\delta_D) |m|}{\overline{\hat{k}}+\delta_D |m|} B^{n-2}
		=-\frac{(2\nu |m|-\delta_D)|m|}{\overline{\hat{k}}+\delta_D |m|} A^{n-1}
		+\frac{2\nu |m|-\delta_D}{2\nu |m|+\delta_S} \frac{(2\nu |m|-\delta_D)|m|}{\overline{\hat{k}}+\delta_D |m|} A^{n-2}.
	\end{eqnarray*}
	Consequently, we can deduce
	\begin{eqnarray*}
		&&B^n-\frac{\overline{\hat{k}}-\hat{k}_j^{-1}}{\overline{\hat{k}}+\delta_D |m|}B^{n-1}+\frac{(2\nu |m|-\delta_D)|m|}{\overline{\hat{k}}+\delta_D |m|} A^{n-1}\\
&&\hspace{25mm}=\frac{2\nu |m|-\delta_D}{2\nu |m|+\delta_S} \Big[B^{n-2}-\frac{\overline{\hat{k}}-\hat{k}_j^{-1}}{\overline{\hat{k}}+\delta_D |m|}B^{n-3} + \frac{(2\nu |m|-\delta_D)|m|}{\overline{\hat{k}}+\delta_D |m|} A^{n-2}  \Big].
	\end{eqnarray*}
	Finally, we arrive at 
	\begin{eqnarray*}
	&&	\left|B^{2n}-\frac{\overline{\hat{k}}-\hat{k}_j^{-1}}{\overline{\hat{k}}+\delta_D |m|}B^{2n-1}+\frac{(2\nu |m|-\delta_D)|m|}{\overline{\hat{k}}+\delta_D |m|} A^{2n-1}\right|
\\
&&\hspace{25mm}=\rho^n(\delta_S,\delta_D,m) \left|B^{1}-\frac{\overline{\hat{k}}-\hat{k}_j^{-1}}{\overline{\hat{k}}+\delta_D |m|}B^{0} + \frac{(2\nu |m|-\delta_D)|m|}{\overline{\hat{k}}+\delta_D |m|} A^{0} \right|,
	\end{eqnarray*}
	by which the convergence factor $\rho(\delta_S,\delta_D,m)$ as given by (\ref{Covfactor}) is demonstrated.
\end{proof}

In order to ensure the convergence of the ensemble algorithm for all relevant frequencies, the optimized Robin parameters $\delta_S, \delta_D>0$ need to be characterized. Inspired by \cite{Discacciati}, the relevant frequencies are assumed in a range of $0<m_{\min} \leq m \leq m_{\max}$, where $m_{\min}=\frac{\pi}{L} $ ($L$ is the length of the interface $\Gamma$) and $m_{\max}=\frac{\pi}{h}$ ($h$ is the mesh size). The Robin parameters should satisfy that $\delta_S <  \delta_D$ and $\rho(\delta_S,\delta_D,m)<1$ for all $m \in [m_{\min},m_{\max}]$. It is clear that the optimal Robin parameters can be easily devised from (\ref{Covfactor}): 
{$\delta^{\mathrm{exact}}_D=2\nu m $ and any $\delta^{\mathrm{exact}}_S < \delta^{\mathrm{exact}}_D$, unfortunately, they are not viable.}

We utilize the well-known min-max technique to find the optimal Robin parameters $\delta_S, \delta_D$ by minimizing the convergence rate over all the relevant frequencies, which is equivalent to solving the min-max problem
\begin{eqnarray}\label{min-max}
	\mathop{\min}_{0<\delta_S <  \delta_D} \mathop{\max}_{m \in [m_{\min},m_{\max}] } \rho (\delta_S,\delta_D,m).
\end{eqnarray}
We can deduce the following theorem that provides the solution to the optimization procedure.

\begin{thm} \label{Thm52}
	The solution of the min-max problem (\ref{min-max}) is: for any $\delta^*_S>0 $,
	\begin{eqnarray}\label{delta_D*}
		\delta^*_D=\frac{4\nu^2 m_{\min} m_{\max}+\nu (m_{\min}+m_{\max})\delta^*_S}{\nu(m_{\min}+m_{\max})+\delta^*_S}.
	\end{eqnarray}
	Moreover, $0<\delta^*_S <  \delta^*_D$ and $\rho(\delta^*_S,\delta^*_D,m)<1$ for all $m \in [m_{\min},m_{\max}]$.
\end{thm}

\begin{proof}
	Following the Lemma 3.2 in \cite{Discacciati} and adapting a similar argument, we can easily derive a similar conclusion: for any given $m \in(0,+\infty)$, $\nabla \rho (\delta^{\mathrm{exact}}_S,\delta^{\mathrm{exact}}_D,m)=0$ and the point $(\delta^{\mathrm{exact}}_S,\delta^{\mathrm{exact}}_D)$ is an {absolute minimum}. So, wherever the maximum with respect to $m$ is, the minimum value with respect to $(\delta_S,\delta_D)$  needs to consider the limitation of $0<\delta_S < \delta_D$. Next, we further consider the square of the convergence factor $\rho(\delta_S,\delta_D,m)$,
 which leads to
	\begin{eqnarray}\label{Covfactor2}
		\hat{\rho}(\delta_D, m) = \Big(\frac{2\nu m -\delta_D}{2\nu m+\delta^*_S}\Big)^2, \hspace{10mm} \mathrm{for} ~ \mathrm{any} ~ \mathrm{given} ~ \delta^*_S>0. 
	\end{eqnarray}
	The minimum value of $\hat{\rho}(\delta_D, m)$ {can be taken to be zero at $m=\frac{\delta_D}{2\nu}$}, meanwhile, since
 the function  in  (\ref{Covfactor2}) is continuous, its maximum must be achieved at either endpoint of the interval $[m_{\min},m_{\max}]$. By simple argument, we can show that $\hat{\rho}(0, m_{\min})<\hat{\rho}(0, m_{\max})$ and $\lim\limits_{\delta_D \rightarrow \infty} \frac{\hat{\rho}(\delta_D, m_{\min})}{\hat{\rho}(\delta_D, m_{\max})}>1$ (other properties are similar to \cite{Discacciati}), so we can derive that
	\begin{eqnarray*}
		\mathop{\max}_{m\in[m_{\min},m_{\max}]} \hat{\rho}(\delta_D, m) = \Big\{ \hat{\rho}(\delta_D, m_{\min}),\ \hat{\rho}(\delta_D, m_{\max}) \Big\}=
		\left\{
		\begin{aligned}
		 	\hat{\rho}(\delta_D, m_{\min}) && \mathrm{for} ~ \delta_D > \delta^*_D,\\
		 	\hat{\rho}(\delta_D, m_{\max}) && \mathrm{for} ~ \delta_D \leq \delta^*_D,
		\end{aligned}
		\right.
	\end{eqnarray*}
	where $\delta^*_D>0$ is the value at which the convergence rate exhibits a balance between the minimum and maximum frequencies, i.e., $\hat{\rho}(\delta^*_D, m_{\min})=\hat{\rho}(\delta^*_D, m_{\max})$. Then, the positive value $\delta^*_D$ given in (\ref{delta_D*}) can be easily derived by simple calculation. We can also easily verify that $\rho(\delta^*_D,m)<1$ for all $m \in [m_{\min},m_{\max}]$  by showing that both $\rho(\delta^*_D,m_{\min})<1$ and  $\rho(\delta^*_D,m_{\max})<1$ hold.
\end{proof}

\begin{rem}\label{OpSD1}
    \begin{figure}[htbp]
		\centering
		\includegraphics[width=55mm,height=45mm]{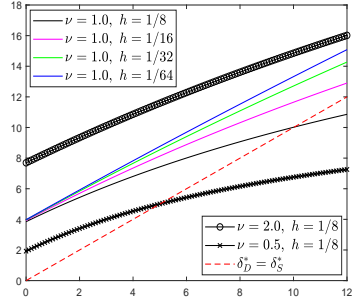}
		\caption{The function $\delta^*_D=\frac{4\nu^2 m_{\min} m_{\max}+\nu (m_{\min}+m_{\max})\delta^*_S}{\nu(m_{\min}+m_{\max})+\delta^*_S}$ with different $\nu$ and $h$.}\label{OpSD}
    \end{figure}
	In Theorem 5.2, the optimal Robin parameters $\delta^*_S$ and $\delta^*_D$ need to satisfy (\ref{delta_D*}) and $0<\delta^*_S <  \delta^*_D$, so their selections should be more cautious. Intuitively, we display some possible choices of $\delta^*_S, \delta^*_D$ with different $\nu$ and $h$ in Fig. 5.1. If the pair of $(\delta^*_S, \delta^*_D)$ locates above the red dotted line, it indicates that such pair is acceptable.
\end{rem}

\section{Finite Element Approximations}
In this section, we will further study the finite element discretization of the Ensemble DDM algorithm. Consider a regular, quasi-uniform triangulation ($d=2$) or tetrahedron ($d=3$) $\mathcal{T}_h$ with mesh scale $h$ for the global domain $\Omega$. For the subdomains $\Omega_S$ and $\Omega_D$, two triangulations ($d=2$) or  tetrahedrons ($d=3$) $\mathcal{T}_{S,h}$, $\mathcal{T}_{D,h}$ are assumed to be compatible at the interface $\Gamma$, on which the triangulation of the meshes on $\Gamma$ is also quasi-uniform. The conforming Stokes velocity, pressure, and Darcy velocity, hydraulic head finite element spaces can be defined to satisfy $\mathbf{X}_{S,h} \subset \mathbf{X}_S, \ Q_{S,h} \subset Q_S, \ \mathbf{X}_{D,h} \subset \mathbf{X}_D, \ Q_{D,h} \subset Q_D$.
Here, the pair of spaces $(\mathbf{X}_{S,h},\ Q_{S,h})$ is  assumed to satisfy the discrete LBB or inf-sup condition, meanwhile the Darcy finite element spaces $(\mathbf{X}_{D,h},\ Q_{D,h})$ are also supposed to permit the standard inf-sup condition. MINI elements $P1b-P1$ for Stokes and $BDM1-P0$ elements for mixed Darcy are one class of suitable choices that will be used in the numerical test.  

The finite element approximation of the decoupled Stokes-Darcy model and the finite element (FE) Ensemble DDM algorithm can be obtained from (\ref{StokesR1})-(\ref{DarcyR2}) and Ensemble DDM algorithm (\ref{decoupled-1})-(\ref{decoupled-comp3}) by adding $h$ in the subscript of all the functions and spaces, and is omitted here.

Next, we will study the convergence of this FE Ensemble DDM algorithm. The error functions and error equations are defined similarly as (\ref{err-1})-(\ref{erref2}), by adding $h$ in the subscripts. Moreover, we have two important equations similarly as (\ref{errorD})-(\ref{errorS}):
\begin{eqnarray}\label{feerrorD}
	||\eta_{Dj,h}^{n}||^2_{\Gamma}
	&=&||\eta_{Sj,h}^{n-1}||^2_{\Gamma}
	+(\delta_D^2-\delta_S^2)||\mathbf{e}_{Sj,h}^{n}\cdot\mathbf{n}_{S}||^2_{\Gamma}
	-2(\delta_S+\delta_D)a_S(\mathbf{e}_{Sj,h}^n,\mathbf{e}_{Sj,h}^n)
	\nonumber\\
	&&-2(\delta_S+\delta_D)\sum_{i=1}^{d-1}\overline{\xi}_i
	||\mathbf{e}_{Sj,h}^n\cdot \tau_{i}||_{\Gamma}^2
	+2(\delta_S+\delta_D)\sum_{i=1}^{d-1}{\xi}_{i,j}
	\langle\mathbf{e}_{Dj,h}^{n-1}\cdot \tau_{i},\mathbf{e}_{Sj,h}^n\cdot \tau_{i}\rangle\nonumber\\
	&&-2(\delta_S+\delta_D)\sum_{i=1}^{d-1}({\xi}_{i,j}-\overline{\xi}_i)
	\langle\mathbf{e}_{Sj,h}^{n-1}\cdot \tau_{i},\mathbf{e}_{Sj,h}^n\cdot \tau_{i}\rangle,\\
	\label{feerrorS}
	||\eta_{Sj,h}^{n}||^2_{\Gamma}&=& ||\eta_{Dj,h}^{n-1}||^2_{\Gamma}
	+(\delta_S^2-\delta_D^2)||\mathbf{e}_{Dj,h}^{n}\cdot\mathbf{n}_{D}||^2_{\Gamma}
	-2(\delta_S+\delta_D)\overline{a_D}(\mathbf{e}_{Dj,h}^n,\mathbf{e}_{Dj,h}^n)\nonumber\\
	&&-2(\delta_S+\delta_D)\Big[a_{Dj}(\mathbf{e}_{Dj,h}^{n-1},\mathbf{e}_{Dj,h}^n)-\overline{a_D}(\mathbf{e}_{Dj,h}^{n-1},\mathbf{e}_{Dj,h}^n)\Big].
\end{eqnarray}

In (\ref{feerrorS}), the term $||\mathbf{e}_{Dj,h}^{n}\cdot\mathbf{n}_{D}||^2_{\Gamma}$ is positive in the case of $\delta_S>\delta_D$, which is one of the main difficult issue encountered in the convergence analysis of the corresponding continuous problem. To resolve such difficulty, we make use of the following trace-inverse inequality \cite{Sun2021}:
\begin{eqnarray}\label{inverse}
	||\mathbf{e}_{Dj,h}^{n}\cdot\mathbf{n}_{D}||^2_{\Gamma}\leq C_5 h^{-1}||\mathbf{e}_{Dj,h}^{n}||_D^2  \leq C_5 h^{-1}||\mathbf{e}_{Dj,h}^{n}||_{\mathrm{div}}^2,
\end{eqnarray}
where the positive constant $C_5$ only depends on the domain $\Omega_D$.

Now we can present the convergence results of the FE Ensemble DDM in the case $\delta_S>\delta_D$.
\newline

\begin{theorem}\label{ItrDDMf>p}
	For the case of $\delta_S>\delta_D$, assume that the hydraulic conductivity tensor $\mathbb{K}_j$ is small enough and the parameters $\delta_S$, $\delta_D$ are chosen to satisfy
	\begin{eqnarray}\label{fecontrolTP}
		&&\delta_S-\delta_D<\frac{hgC_1\nu(\overline{k}^{\min}-E_{k_j}^{\max})-2h(C_{\mathrm{tr}}C'_{\mathrm{tr}}\xi_j)^2}{4C_1 C_5\nu}, \hspace{8mm} 
		\rho'(\delta_S,\delta_D)<1.
	\end{eqnarray}
	where $\rho'(\delta_S,\delta_D):=\frac{(1+C_3)^2g(\delta_S-\delta_D)[2(\overline{k}^{\max})^2+(E_{k_{j}}^{\max})^2]+\delta_D^2(\overline{k}^{\min}-E_{k_{j}}^{\max})}{4(\delta_S^2-\delta_D^2)(\overline{k}^{\min}-E_{k_{j}}^{\max})}$, with the conditions
	\begin{eqnarray}\label{perturbation2}
		&&E_{{\xi}_{i,j}}^{\mathrm{max}} < \overline{\xi}_j, \hspace{10mm}
		E_{k_j}^{\mathrm{max}} < \overline{k}^{\mathrm{min}}.
	\end{eqnarray}
	Then the following convergence estimates for the FE Ensemble DDM algorithm are guaranteed
	\begin{eqnarray}\label{result2>}
		&&||\mathbf{e}_{Sj,h}^N||_1+||\varepsilon_{Sj,h}^{N}||^2+ ||\mathbf{e}_{Dj,h}^{N-1}||_{\mathrm{div}}+ ||\varepsilon_{Dj,h}^{N}||^2\nonumber \\
  &&\hspace{8mm}+||\eta^{N}_{Sj,h}||_{\Gamma}^2
		+ ||\eta^{N}_{Sj,\tau,h}||_{H^{-\frac{1}{2}}(\Gamma)}^2
		+||\eta^{N}_{Dj,h}||_{\Gamma}^2
  +\rho'(\delta_S,\delta_D)||\eta^{N-1}_{Dj,h}||_{\Gamma}^2\nonumber \\
		&&\leq \bar{C}^*	\max\Big\{ \sqrt{\rho'(\delta_S,\delta_D)},\  \sum_{i=1}^{d-1}\frac{E_{{\xi}_{i,j}}^{\max}}{2\overline{\xi}_i-E_{{\xi}_{i,j}}^{\max}},\ 
		CR\Big\}^{N-2} \Big[ \sqrt{\rho'(\delta_S,\delta_D)}\Big(||\eta_{Dj,h}^{1}||^2_{\Gamma}\\ 
&&\hspace{2.2mm} +\sqrt{\rho'(\delta_S,\delta_D)}  ||\eta_{Dj,h}^{0}||^2_{\Gamma}\Big)+(\delta_S+\delta_D) \sum_{i=1}^{d-1} E_{{\xi}_{i,j}}^{\max}
		||\mathbf{e}_{Sj,h}^{1}\cdot \tau_{i}||_{\Gamma}^2+(\delta_S+\delta_D) g\overline{k}^{\min} ||\mathbf{e}_{Dj,h}^{0}||_{\mathrm{div}}^2 \Big],\nonumber
	\end{eqnarray}
	for a given positive constant $\bar{C}^*$, here the convergence rate $CR$ is defined by
	\begin{eqnarray}\label{CR}
	CR=\frac{2hg\overline{k}^{\min}}{h\Big(3g\overline{k}^{\min}-gE_{k_j}^{\max}-\frac{2(C_{\mathrm{tr}}C'_{\mathrm{tr}}\xi_j)^2}{C_1\nu}\Big)-2C_5(\delta_S-\delta_D)}.
	\end{eqnarray}
\end{theorem}

\begin{proof}
	 Similarly with the analysis in (\ref{contral_eta}), substituting the test function $\mathbf{v}_D$ defined as (\ref{etacontrolv}) into (\ref{err-3}), and applying the Cauchy-Schwarz inequality, Young's inequality and (\ref{PDcontrol}) to derive
	\begin{eqnarray}\label{fecontral_eta}
		||\eta_{Dj,h}^{n-2}||_{\Gamma}^{2}
		&\leq& \frac{1}{2}(\delta_S+\delta_D)g(\overline{k}^{\min}-E_{k_j}^{\max})||\mathbf{e}_{Dj,h}^{n-1}||_{\mathrm{div}}^2+(\delta_S^2-\delta_D^2)||\mathbf{e}_{Dj,h}^{n-1}\cdot\mathbf{n}_{D}||^2_{\Gamma}\nonumber\\
		&&+(\delta_S+\delta_D)g(\overline{k}^{\min}-E_{k_j}^{\max})||\mathbf{e}_{Dj,h}^{n-2}||_{\mathrm{div}}^2+\rho'(\delta_S,\delta_D)||\eta_{Dj,h}^{n-2}||_{\Gamma}^2,
	\end{eqnarray}
	where 
	\begin{eqnarray}\label{rhoSD'}
		\rho'(\delta_S,\delta_D):=\frac{(1+C_3)^2g(\delta_S-\delta_D)[2(\overline{k}^{\max})^2+(E_{k_j}^{\max})^2]+\delta_D^2(\overline{k}^{\min}-E_{k_j}^{\max})}{4(\delta_S^2-\delta_D^2)(\overline{k}^{\min}-E_{k_j}^{\max})}.
	\end{eqnarray}
	By replacing the term $||\eta_{Sj,h}^{n-1}||_{\Gamma}^2$ by (\ref{feerrorS}) in the equation (\ref{feerrorD}), and then using (\ref{bilinear-ieq})-(\ref{term3}) and (\ref{contral_eta}), also by the Poincar$\acute{\mathrm{e}}$, Korn's, Young's and the trace-inverse (\ref{inverse}) inequalities, we arrive at
	\begin{eqnarray*}
		&&\hspace{-10mm}||\eta_{Dj,h}^{n}||^2_{\Gamma}
		\leq\rho'(\delta_S,\delta_D)||\eta_{Dj,h}^{n-2}||^2_{\Gamma} \\
		&&\hspace{0mm}-2(\delta_S+\delta_D)\Big[ \frac{C_1\nu}{2}||\mathbf{e}_{Sj,h}^{n}||_{1}^2+\sum_{i=1}^{d-1}\frac{2\overline{\xi}_i-E_{{\xi}_{i,j}}^{\max}}{2}
		||\mathbf{e}_{Sj,h}^n\cdot \tau_{i}||_{\Gamma}^2-\sum_{i=1}^{d-1}\frac{E_{{\xi}_{i,j}}^{\max}}{2}
		||\mathbf{e}_{Sj,h}^{n-1}\cdot \tau_{i}||_{\Gamma}^2 \\ 
		&&\hspace{0mm}+\Big(\frac{3g\overline{k}^{\min}}{4}-\frac{gE_{k_j}^{\max}}{4}-C_5h^{-1}(\delta_S-\delta_D)-\frac{(C_{\mathrm{tr}}C'_{\mathrm{tr}}\xi_j)^2}{2C_1\nu} \Big) ||\mathbf{e}_{Dj,h}^{n-1}||_{\mathrm{div}}^2 -\frac{g\overline{k}^{\min}}{2}||\mathbf{e}_{Dj,h}^{n-2}||_{\mathrm{div}}^2\Big].
	\end{eqnarray*}
	Under the assumption
	\begin{eqnarray}\label{feadd_condition}
		\delta_S-\delta_D<\frac{hgC_1\nu(\overline{k}^{\min}-E_{k_j}^{\max})-2h(C_{\mathrm{tr}}C'_{\mathrm{tr}}\xi_j)^2}{4C_1 C_5\nu},
	\end{eqnarray}
	we can immediately have
$\frac{3g\overline{k}^{\min}}{4}-\frac{gE_{k_j}^{\max}}{4}-C_5h^{-1}(\delta_S-\delta_D)-\frac{(C_{\mathrm{tr}}C'_{\mathrm{tr}}\xi_j)^2}{2C_1\nu}>\frac{g\overline{k}^{\min}}{2}$.
If
	\begin{eqnarray*}
		\rho'(\delta_S,\delta_D)<1,\hspace{6mm} \overline{k}^{\min}>E_{k_j}^{\max},\hspace{6mm}
		\overline{\xi}_i>E_{{\xi}_{i,j}}^{\max},
	\end{eqnarray*}
are assumed, then by applying a similar argument with $\delta_S<\delta_D$ above, we can conclude that
	\begin{eqnarray*}
		&&\hspace{-8mm}\Big(||\eta_{Dj,h}^{N}||^2_{\Gamma}+\sqrt{\rho'(\delta_S,\delta_D)}||\eta_{Dj,h}^{N-1}||^2_{\Gamma}\Big)+2(\delta_S+\delta_D)\sum_{i=1}^{d-1}\frac{2\overline{\xi}_i-E_{{\xi}_{i,j}}^{\max}}{2}
		||\mathbf{e}_{Sj,h}^N\cdot \tau_{i}||_{\Gamma}^2\\
		&&\hspace{6mm}+2(\delta_S+\delta_D)\Big(\frac{3g\overline{k}^{\min}-gE_{k_j}^{\max}}{4}-C_5h^{-1}(\delta_S-\delta_D)-\frac{(C_{\mathrm{tr}}C'_{\mathrm{tr}}\xi_j)^2}{2C_1\nu} \Big) ||\mathbf{e}_{Dj,h}^{N-1}||_{\mathrm{div}}^2\\
		&&\hspace{-6mm}\leq	\max\Big\{ \sqrt{\rho'(\delta_S,\delta_D)},\  \sum_{i=1}^{d-1}\frac{E_{{\xi}_{i,j}}^{\max}}{2\overline{\xi}_i-E_{{\xi}_{i,j}}^{\max}},\ \frac{2hg\overline{k}^{\min}}{h\Big(3g\overline{k}^{\min}-gE_{k_j}^{\max}-\frac{2(C_{\mathrm{tr}}C'_{\mathrm{tr}}\xi_j)^2}{C_1\nu}\Big)-2C_5(\delta_S-\delta_D)} \Big\}^{N-2}\\
		&&\hspace{6mm}\times \Big\{ \sqrt{\rho'(\delta_S,\delta_D)}\Big(||\eta_{Dj,h}^{1}||^2_{\Gamma}+\sqrt{\rho'(\delta_S,\delta_D)}||\eta_{Dj,h}^{0}||^2_{\Gamma}\Big)\\
 &&\hspace{15mm} +(\delta_S+\delta_D)\Big[ \sum_{i=1}^{d-1}E_{{\xi}_{i,j}}^{\max}
		||\mathbf{e}_{Sj,h}^{1}\cdot \tau_{i}||_{\Gamma}^2+g\overline{k}^{\min}||\mathbf{e}_{Dj,h}^{0}||_{\mathrm{div}}^2 \Big]\Big\},
	\end{eqnarray*}
	which clearly yields the geometric convergence of $ ||\eta_{Dj,h}^{n}||_{\Gamma}^{2}, \ ||\mathbf{e}_{Sj,h}^{n}\cdot \tau_{i}||_{\Gamma}^2$, and $ ||\mathbf{e}_{Dj,h}^{n}||^2_{\mathrm{div}}$ when the hydraulic conductivity tensor $\mathbb{K}_j$ is small enough.
	
	Finally by (\ref{feerrorS}), the geometric convergence of $||\eta_{Sj,h}^{n}||^2_{\Gamma}$ can be further obtained
	\begin{eqnarray*}
		||\eta_{Sj,h}^{n}||^2_{\Gamma}&\leq&
		||\eta_{Dj,h}^{n-1}||^2_{\Gamma}
		-(\delta_S+\delta_D)(2g\overline{k}^{\min}-gE_{k_j}^{\max})||\mathbf{e}_{Dj,h}^n||_{\mathrm{div}}^2 +(\delta_S+\delta_D)gE_{k_j}^{\max}||\mathbf{e}_{Dj,h}^{n-1}||_{\mathrm{div}}^2\\
		&\leq&||\eta_{Dj,h}^{n-1}||^2_{\Gamma}
		+(\delta_S+\delta_D)gE_{k_j}^{\max}||\mathbf{e}_{Dj,h}^{n-1}||_{\mathrm{div}}^2.
	\end{eqnarray*}	
 The geometric convergence of $||\mathbf{e}_{Sj,h}^{n}||_1^2$ follows directly by (\ref{aS+bS}). The proof is complete.
\end{proof}

\begin{rem}\label{smallK11}
	The assumption (\ref{feadd_condition}) is clearly tenable when the hydraulic conductivity tensor $\mathbb{K}_j$ and the $\alpha$ in  $\xi_{j}=\sum_{i=1}^{d-1}\frac{\alpha}{\sqrt{\tau_i\cdot\mathbb{K}_j\tau_i}}$  are small enough. If $\mathbb{K}_j$ is small enough, then $\overline{k}^{\min}$ is larger than $h^{-1}$, hence can control $h$, otherwise, with the refinement of the triangulation, the occurrence of smaller $h$ will lead to $\delta_S=\delta_D$. Moreover, if we want to obtain geometric convergence, the convergence rate $CR$ as defined by (\ref{CR}) should be independent of the mesh size $h$, therefore we need to enforce the hydraulic conductivity tensor $\mathbb{K}_j$ to be small enough to control $h$. In particular, in order to acquire the condition (\ref{fecontrolTP}), one reasonable  suggestion is selecting Robin parameters $\delta_S$ and $\delta_D$ to have the same order of magnitude with $\overline{k}^{\min}$.
\end{rem}


\begin{rem}
	 A small perturbation constraint of the random hydraulic conductivity $\mathbb{K}(\mathbf{x},\omega)$ is required for the convergence of the Ensemble DDM, which has been presented in both (\ref{perturbation}) and (\ref{perturbation2}).  We also note that this assumption commonly occurs in practical applications. If the random $\mathbb{K}(\mathbf{x},\omega)$ has a larger disturbance, we can also modify our algorithm for numerical simulation. Actually, we can separate it into different levels according to the order of magnitude, and then solve the corresponding parts. Meanwhile, if the $\mathbb{K}(\mathbf{x},\omega)$  possesses a normal distribution, we can use the `3-$\sigma$' principle to choose the disturbance range before applying our algorithm, where the $\sigma$ is the standard deviation of the normal distribution.
\end{rem}

\section{Numerical Experiments}
In this section, we present two numerical experiments to illustrate the approximate accuracy and efficiency of the proposed ensemble DDM for the fully-mixed random Stokes-Darcy fluid flow model. In the first numerical experiment, we test a smooth problem to check the convergence of our ensemble DDM. We will show the performance of the combinations of the ensemble DDM and Monte Carlo method in the second example, where a random Stokes-Darcy model problem with a random hydraulic conductivity tensor is used. In both tests, the finite element spaces are constructed by well-known MINI (\emph{P1b-P1}) elements for the Stokes problem and \emph{Brezzi-Douglas-Marini} (\emph{BDM1-P0}) elements for the mixed Darcy problem. 

All the numerical tests are implemented by the open software FreeFEM++ \cite{F18}. The stopping criteria for the iterative process of ensemble DDM are usually selected as a fixed tolerance of $10^{-6}$  between two successive solution components of the Stokes and Darcy velocities in the sense of $L^2$-norm, $i.e.$,
	$\Bigl(||\mathbf{u}^{n+1}_{Sj,h}-\mathbf{u}_{Sj,h}^n||_S^2+||\mathbf{u}^{n+1}_{Dj,h}-\mathbf{u}_{Dj,h}^n||_D^2 \Bigr)^{1/2}\leq 10^{-6}$.

\subsection{Smooth Problem with Convergence Test}
This testing example with an exact solution is adapted from \cite{JiangN19} to verify the convergence and check the feasibility of Ensemble DDM. The free fluid flow region $\Omega_{S}= [0,\pi]\times[0,1]$ and the porous medium region $\Omega_{D}=[0,\pi]\times[-1,0]$ are considered, including the interface  $\Gamma= \{ 0 \leq x \leq \pi, y=0\}$. For the computational convenience, we assume $z=0$, and  other physical parameters $ \nu, g$ and $\alpha$ to be $ 1.0 $. The random hydraulic conductivity tensor $\mathbb{K}$ will be assumed as
\begin{eqnarray*}
	\mathbb{K}=\mathbb{K}_{j}=\left[\begin{array}{cc}
		k_{11}^{j} & 0 \\
		0 & k_{22}^{j}
	\end{array}\right], \quad j=1, \ldots, J,
\end{eqnarray*}
where $\mathbb{K}_{j}$ is $j$th sample of $\mathbb{K}$.
The exact solution is selected as:
\begin{eqnarray*}
	\begin{aligned}
	\mathbf{u}_D&= -\mathbb{K}_{j} \nabla \phi_D,\quad	\phi_{D} =\left(e^{y}-e^{-y}\right) \sin (x),\\
		\mathbf{u}_{S} &=\left[\frac{k_{11}^{j}}{\pi} \sin (2 \pi y) \cos (x),\Big(-2 k_{22}^{j}+\frac{k_{22}^{j}}{\pi^{2}} \sin ^{2}(\pi y)\Big) \sin (x)\right]^{T},\quad  	p_{S} =0.
	\end{aligned}
\end{eqnarray*}
Coincidentally, the above exact solution from \cite{JiangN19} satisfies both the BJS and BJ interface conditions.

\begin{figure}[htbp]
	\centering
	\includegraphics[width=48mm,height=40mm]{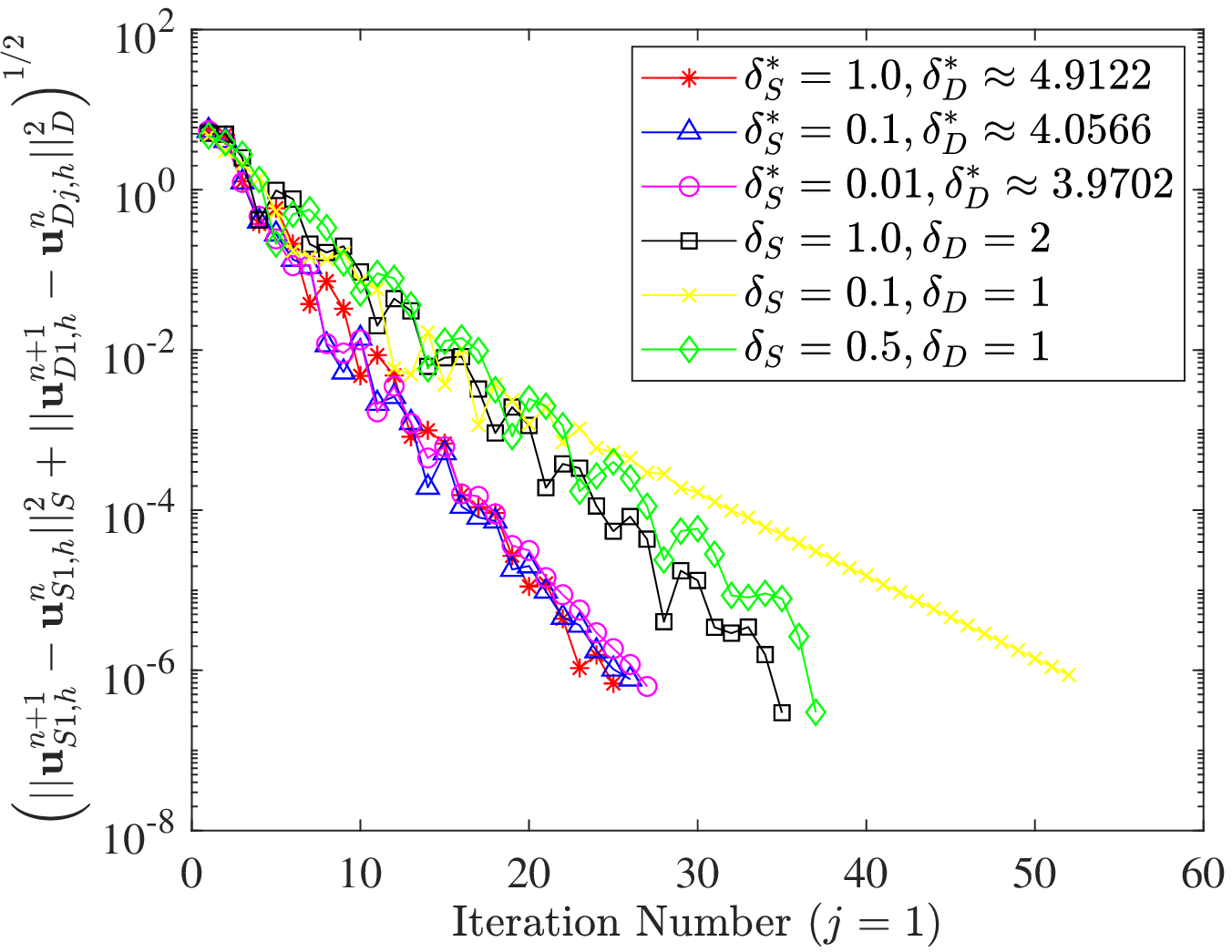} 
	\includegraphics[width=48mm,height=40mm]{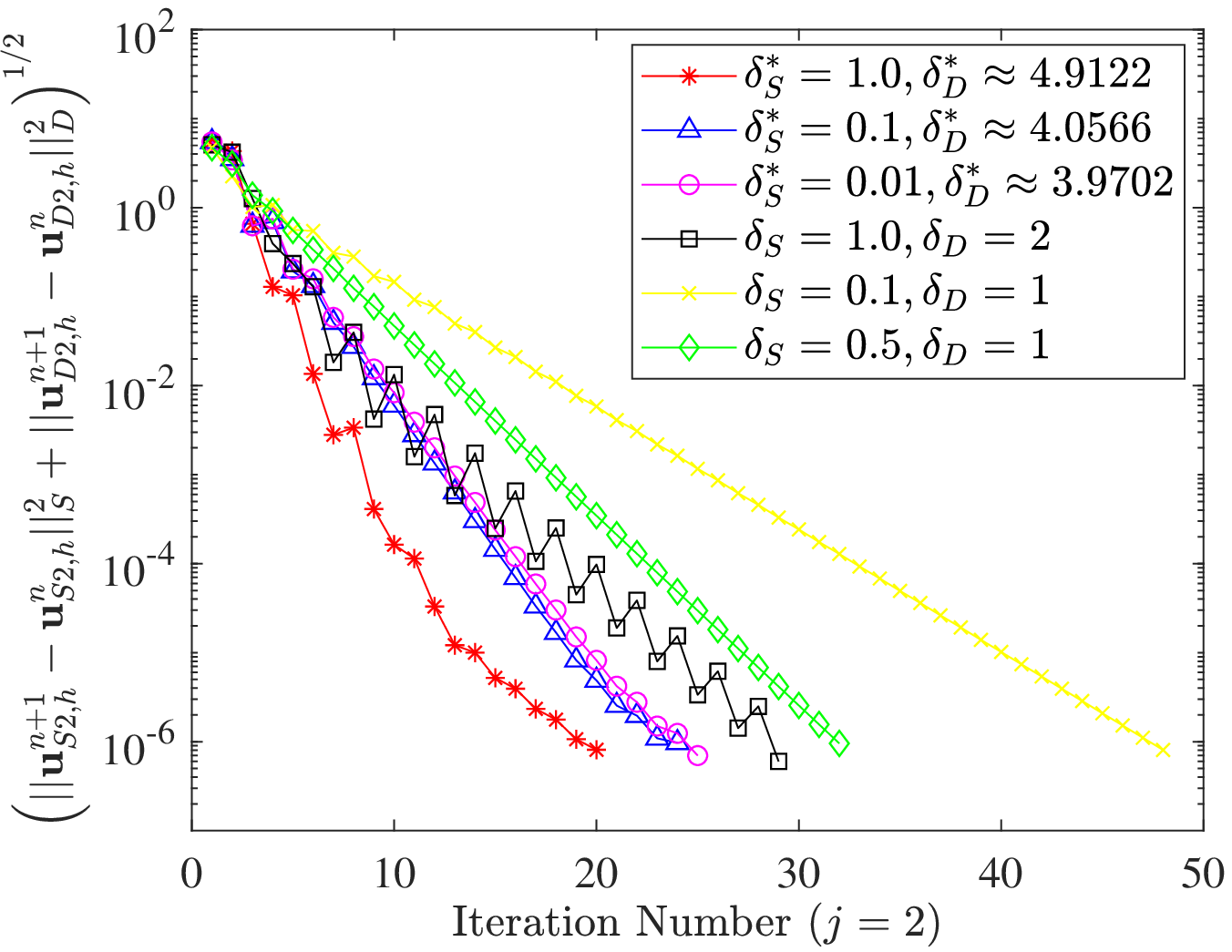} 
	\includegraphics[width=48mm,height=40mm]{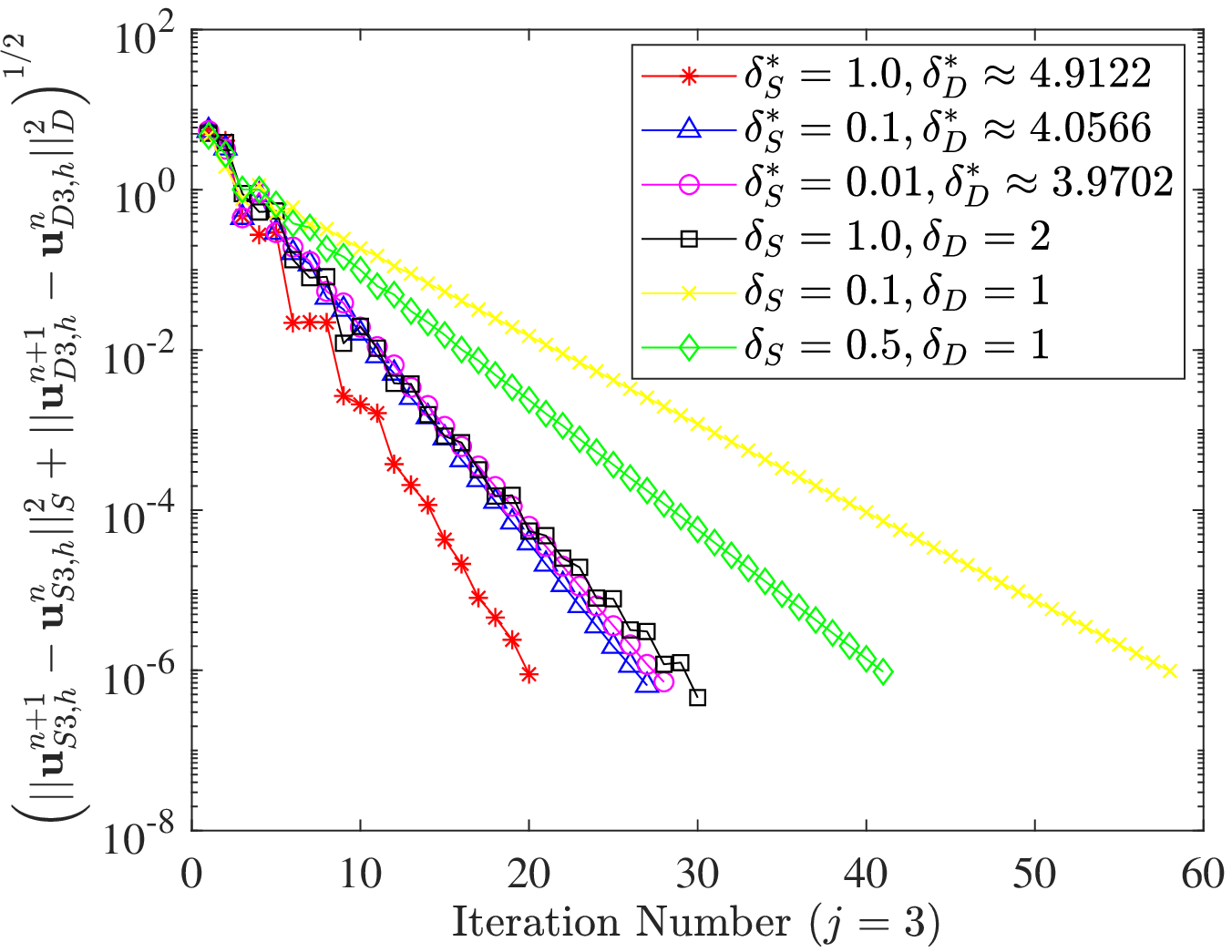}
	\caption{The iterates of the Ensemble DDM ($J=3$) with different Robin parameters $\delta_S$ and $\delta_D$ with $h=\frac{1}{32}$.}
	\label{OpSDerr}
\end{figure}

To verify the theoretical analysis regarding the optimal Robin parameters, we test the proposed Ensemble DDM with different Robin parameters $\delta_S, \delta_D$ while $h=\frac{1}{32}$, and 
compare the corresponding numbers of the iterations. We display the evolution of the $L^2$-error between two successive solution components of the Stokes and Darcy velocities with increasing iteration steps in Fig. \ref{OpSDerr}. In order to demonstrate that actually $\delta_S$  can be chosen arbitrarily, as explained in Remark \ref{OpSD1} and also shown in Fig. \ref{OpSD}, we select $\delta^*_S=1.0, 0.1, 0.01$, and then compute corresponding $\delta^*_D \approx 4.9122, 4.0566,3.9702$ to carry out numerical experiments. Meanwhile we use three groups of the Robin parameters for comparisons, including $(\delta_S,\delta_D)=(1.0,2.0), (0.1,1.0)$ and $(0.5,1.0)$. Among the compared Robin parameters, the pair  selected by Theorem 5.2 will converge fast and have fewer iterations, which positively supports our theoretical analysis.

For better comparison with \cite{JiangN19}, we carry out the same group of simulations with $J=3$, by the selections of  hydraulic conductivity as $k_{11}^{1}=k_{22}^{1}=2.21,~ k_{11}^{2}=k_{22}^{2}=4.11,~ k_{11}^{3}=k_{22}^{3}=6.21$. Based on the discussion in Remark 4.2, the hydraulic conductivity and the viscosity coefficient are almost $O(1)$, we can only check the geometric convergence of Ensemble DDM when the Robin parameters satisfy $\delta_S<\delta_D$. We present numerical errors in Table \ref{OrderEDDME1}, including the approximate accuracy of the velocity components in $L^2$-norm, $H^1$-norm, and $H_{\mathrm{div}}$-norm, and that of the pressure components in $L^2$-norm in both free fluid flow and porous medium regions. From this table, We can observe that our algorithm has a $h$-independent convergence rate in the case of $\delta_S<\delta_D$, which supports Theorem \ref{ItrDDMf<p}. Numerical results in Table \ref{OrderEDDME1} verify the optimal convergence orders for both velocity and pressure, and also demonstrate very similar approximate accuracy as \cite{JiangN19}.
\begin{table}[!h]
	\caption{Convergence performance for Ensemble DDM (J=3) with different $\mathbb{K}$ while $\delta_S=1.0$ and $\delta_D=\frac{5\pi+h}{\pi+2h}$.}
	\label{OrderEDDME1}\tabcolsep 0pt \vspace*{-10pt}
	\par
	\begin{center}
		\def\temptablewidth{1.0\textwidth}
		{\rule{\temptablewidth}{1pt}}
		\begin{tabular*}{\temptablewidth}{@{\extracolsep{\fill}}cc|cccccccc}
			\hline
			$k_{11}^j$&$k_{22}^j$&$h $& iteration & $\frac{||\mathbf u_S-\mathbf u_{S,h}||_{S}}{||\mathbf u_S||_{S}}$ & $\frac{||\mathbf u_S-\mathbf u_{S,h}||_{1}}{||\mathbf u_S||_{1}}$ & $\frac{||p_S-p_{S,h} ||_{S}}{||p_S||_S}$  &  $\frac{||\phi_D-\phi_{D,h}||_{D}}{||\phi_D||_D}$   & $\frac{||\mathbf u_D-\mathbf u_{D,h}||_{D}}{||\mathbf u_D||_D} $ & $\frac{||\mathbf u_D-\mathbf u_{D,h}||_{\mathrm{div}}}{||\mathbf u_D||_{\mathrm{div}}} $\\
			\hline
			2.21 & 2.21 &  $\frac{1}{16}$ & 25   & 0.0019640   &      0.0690629   &     0.0427929   &     0.0358742    &     0.0021811    &    0.0006461
			\\ 
			& & $\frac{1}{32}$ &25      &    0.0004933     &    0.0345968   &     0.0136990   &     0.0179326    &     0.0005484   &     0.0001625
			\\
			& & $\frac{1}{64}$ &25     &      0.0001234    &     0.0173043    &    0.0045482     &   0.0089656  &       0.0001376  &      0.0000408
			\\
			\hline
			4.11 & 4.11  &  $\frac{1}{16}$  &22     &   0.0019656     &    0.0690630   &     0.0422920   &     0.0358749     &    0.0021903   &     0.0006489
			\\ 
			& & $\frac{1}{32}$ &20     &     0.0004936    &     0.0345969   &     0.0135924    &    0.0179327   &      0.0005509     &   0.0001632 
			\\
			& & $\frac{1}{64}$ & 17    &    0.0001234    &     0.0173043   &     0.0045271    &    0.0089656   &      0.0001382   &     0.0000409
			\\
			\hline
			6.21 & 6.21  & $\frac{1}{16}$  &22     &     0.0019662  &       0.0690630    &    0.0420951    &    0.0358751      &   0.0021941    &    0.0006499
			\\ 
			& & $\frac{1}{32}$ &20   &      0.0004938     &    0.0345969     &   0.0135503      &  0.0179327   &      0.0005516     &   0.0001634
			\\
			& & $\frac{1}{64}$ &20    &      0.0001235  &       0.0173043   &     0.0045186   &     0.0089657   &      0.0001382    &    0.0000409
			\\
			\hline
		\end{tabular*}%
	\end{center}
\end{table}

\begin{figure}[htbp]
	\centering
	\includegraphics[width=50mm,height=42mm]{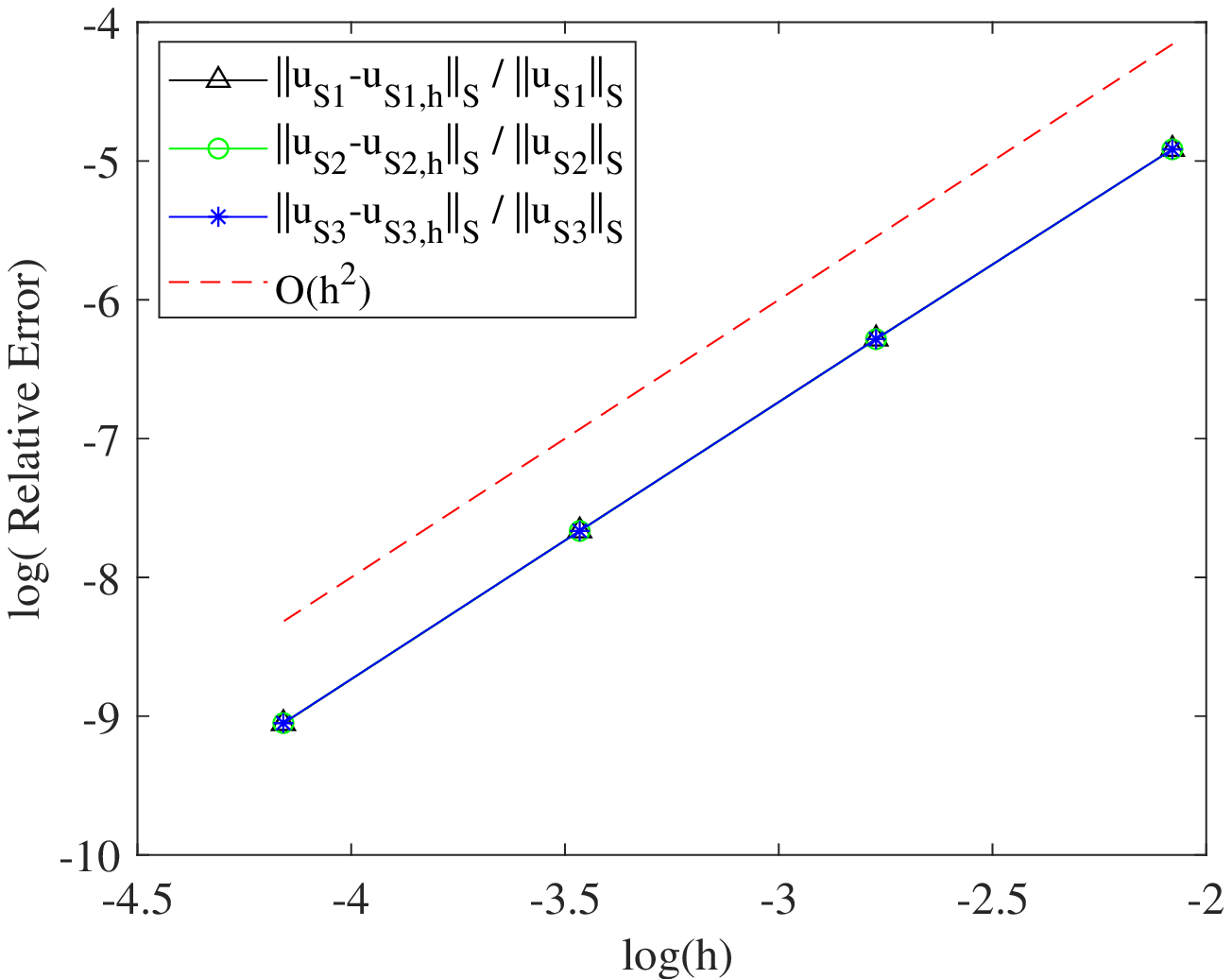} \hspace{3mm}
	\includegraphics[width=50mm,height=42mm]{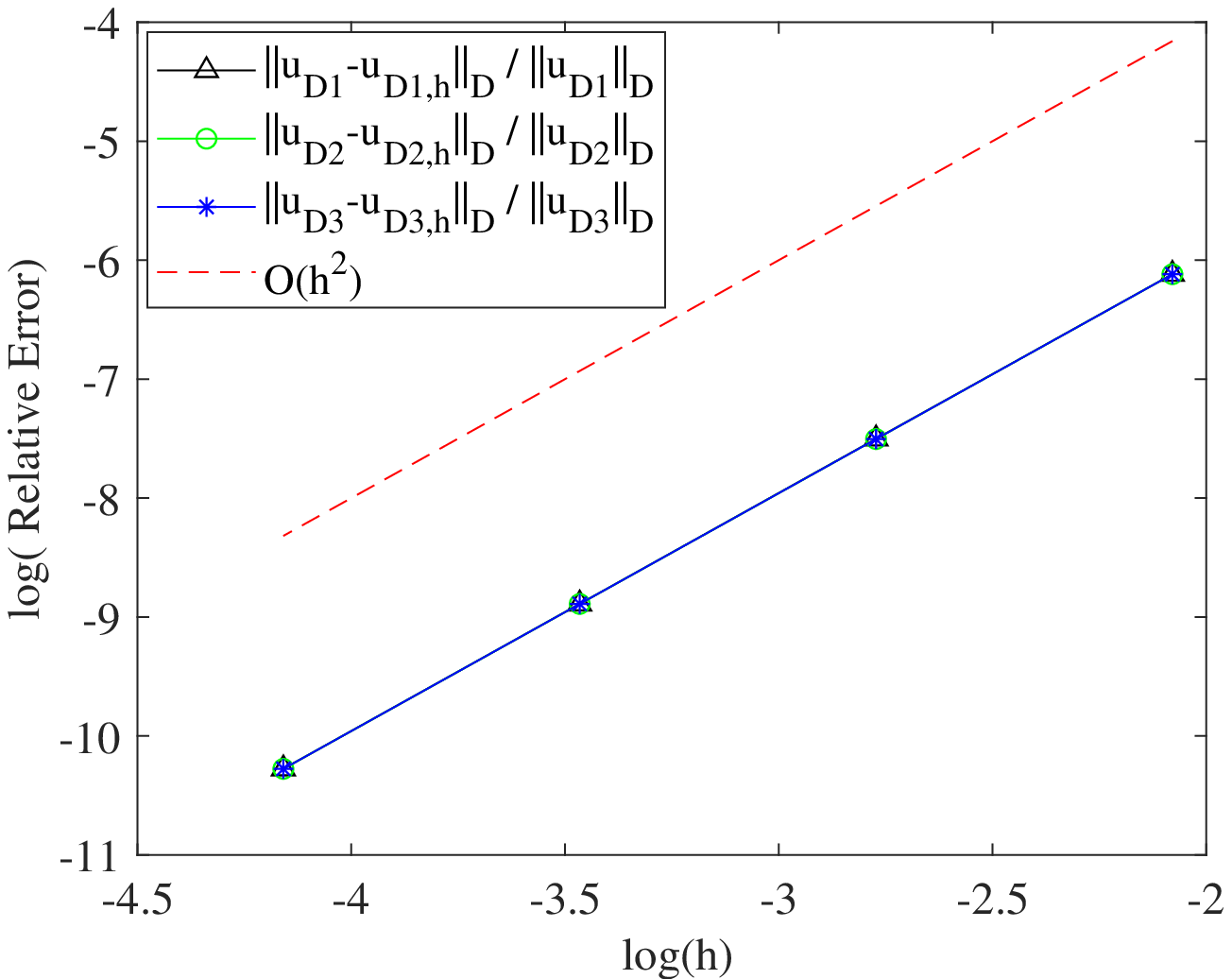}
	\caption{The $L^2$-errors of the Stokes (left) and Darcy (right) velocities for each sample with  $\delta_S=100, \ \delta_D=50$ while the hydraulic conductivity $\mathbb{K}_j$ are $O(10^{-4})$.}
	\label{L2smallK}
\end{figure}

Furthermore, it would be more interesting to verify the convergence orders of the  Ensemble DDM when the hydraulic conductivity $\mathbb{K}$ is significantly small. Noting that the exact solution above admits the variation of the hydraulic conductivity $\mathbb{K}$. We choose $k_{11}^{1}=k_{22}^{1}=1.e-4,~ k_{11}^{2}=k_{22}^{2}=2.e-4,~ k_{11}^{3}=k_{22}^{3}=3.e-4$ to further test our Ensemble DDM under the case of the realistic physical parameters $\mathbb{K}$. According to the theoretical analysis in Section 6, we choose the Robin parameters $\delta_S=100, \ \delta_D=50$ and modify the stopping criteria to a tolerance of $10^{-9}$, due to a smaller order of magnitude of the solutions caused by small permeability coefficients. Under the computational circumstances above, for simplicity, we execute our algorithm and display the convergence orders of the velocities in $L^2$-norms for both regions in Fig. \ref{L2smallK}.
We can observe that the $L^2$-errors of the Stokes and Darcy velocities for each sample have achieved the optimal orders, which indicates the effectiveness of our algorithm.

\subsection{``Shallow Water" System with Random Hydraulic Conductivity}
\begin{figure}[htbp]
	\centering
	\includegraphics[width=33mm,height=38.5mm]{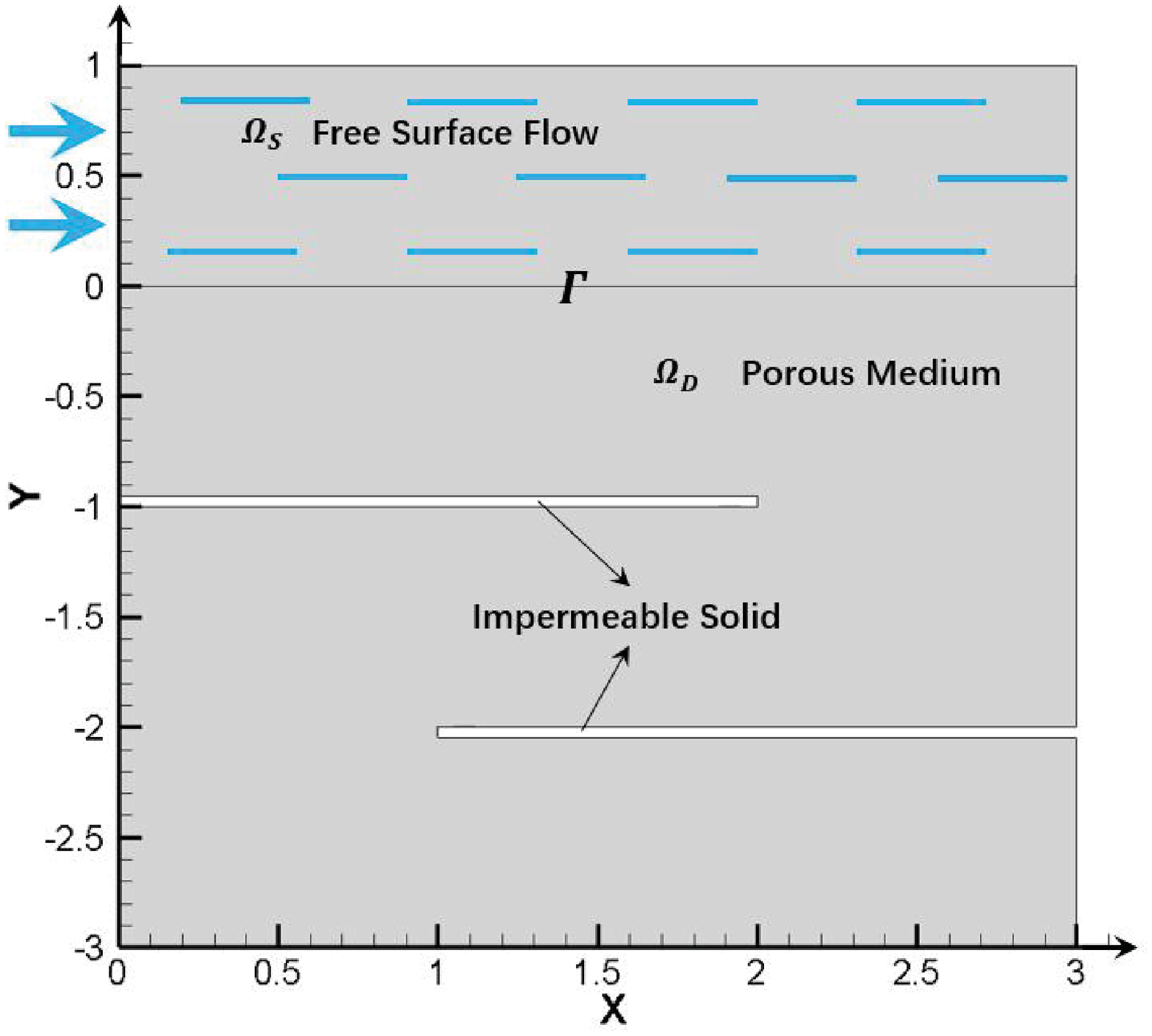} \hspace{2mm}
	\includegraphics[width=45mm,height=40mm]{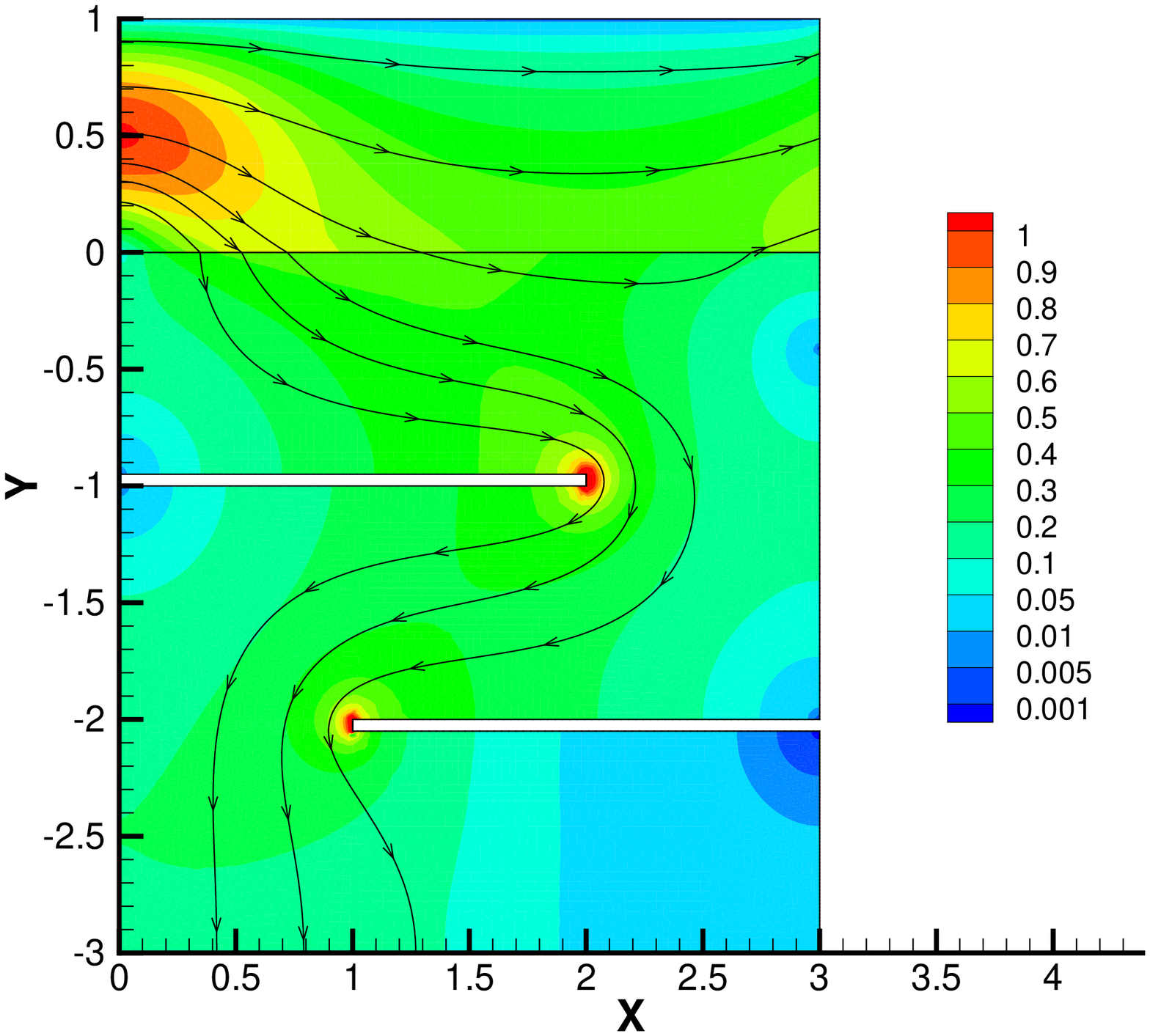}
	\includegraphics[width=45mm,height=40mm]{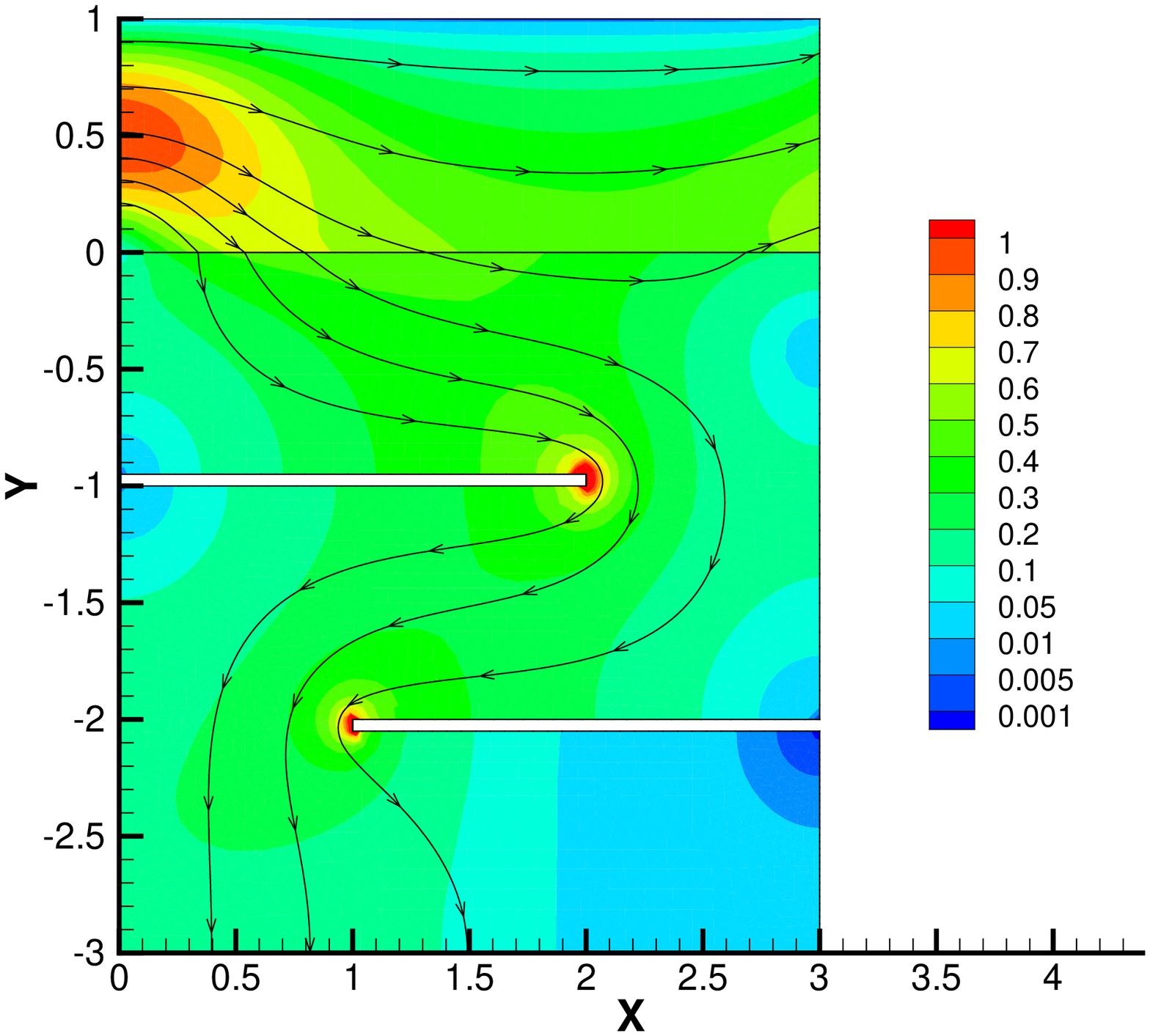}
	\caption{The cross-section of the coupling 3D Shallow Water system $\Omega_S$ with the porous medium $\Omega_D$ (right) and the velocity streamlines of
	the  velocity expectations of Ensemble DDM (left) and traditional DDM (right).}
	\label{3D}
\end{figure}

 The second test is to show the efficiency of the proposed ensemble DDM under much practical permeability $\mathbb{K}$. We simulate a more complicated couple fluid flow model, 
as illustrated in Fig. \ref{3D} (right). We assume a slightly different situation for the computational domain, namely the water channel of length $3.0$ ($X$ axis) and height $1.0$ ($Y$ axis) covered by the porous media with a depth of $3.0$ ($Y$ axis).  The two impermeable solids can be expressed explicitly as $[0,2.0]\times [-1.0,-0.95] \cup [1.0,3.0]\times [-2.05,-2.0]$ in the conceptual domain. We impose  $\mathbf{u}_S=[4y(1-y),0]$  as the inflow surface velocity and non-reflective boundary condition on the outlet surface. The rest of the boundaries except the interface of the Stokes subdomain are treated with no-slip boundary conditions.  Impermeable boundary condition $\mathbf{u}_D \cdot \mathbf{n}_D=0$ is considered on the vertical ($Y$ direction) face of the porous medium. Homogeneous Dirichlet boundary condition  $\phi_D=0$ is applied on the bottom surface of the porous medium. Let $\nu=1.0$, $z=0$, $g=1.0$, $\mathbf{f}_S=\mathbf{0}$ and $f_D=0$.  Actually, this model is the cross-section of the coupling 3D  Shallow Water system with the porous medium in \cite{DMQ02, Sun2021}.

We construct the random $\mathbb{K}$ that varies in the vertical direction as follows
\begin{eqnarray*}
	\begin{aligned}
		&\mathbb{K}(\mathbf{x}, \omega)=\left[\begin{array}{cc}
			k_{11}(\mathbf{x}, \omega) & 0 \\
			0 & k_{22}(\mathbf{x}, \omega)
		\end{array}\right],\quad \text { with }\,\, k_{11}(\mathbf{x}, \omega)=k_{22}(\mathbf{x}, \omega)=k(\mathbf{x}, \omega) \quad \text { and } \\
		&k(\mathbf{x}, \omega)=a_{0}+\sigma \sqrt{\lambda_{0}} Y_{0}(\omega)+\sum_{i=1}^{n_{f}} \sigma \sqrt{\lambda_{i}}\left[Y_{i}(\omega) \cos (i \pi y)+Y_{n_{f}+i}(\omega) \sin (i \pi y)\right],
	\end{aligned}
\end{eqnarray*}
where $\mathbf{x}=(x, y)^{T}, \lambda_{0}=\frac{\sqrt{\pi L_{c}}}{2}, \lambda_{i}=\sqrt{\pi} L_{c} e^{-\frac{\left(i \pi L_{c}\right)^{2}}{4}}$ for $i=1, \ldots, n_{f}$ and $Y_{0}, \ldots, Y_{2 n_{f}}$ are uncorrelated random variables with unit variance and zero mean. In the numerical test, we choose the desired physical correlation length $L_{c}=0.25$ for the random field and $a_{0}=1, \sigma=0.15, n_{f}=3$. The random variables $Y_{0}, \ldots, Y_{2 n_{f}}$ are assumed to be independent and uniformly distributed in the interval $[-\sqrt{3}, \sqrt{3}] $. With above assumptions, we note that  $\mathbb{K}$ is symmetric positive definite, since the random functions $k_{11}(\mathbf{x}, \omega)$ and $k_{22}(\mathbf{x}, \omega)$ are both positive. Moreover, we can follow Theorem \ref{Thm52} to choose $\delta_S=1$ and $\delta_D=\frac{4\nu^2 m_{\min} m_{\max}+\nu (m_{\min}+m_{\max})\delta_S}{\nu(m_{\min}+m_{\max})+\delta_S}$ ($m_{\min}=\frac{\pi}{3}, m_{\max}=\frac{\pi}{h}$) for simulating this random model faster. 

\begin{figure}[htbp]
	\centering
	\includegraphics[width=50mm,height=42mm]{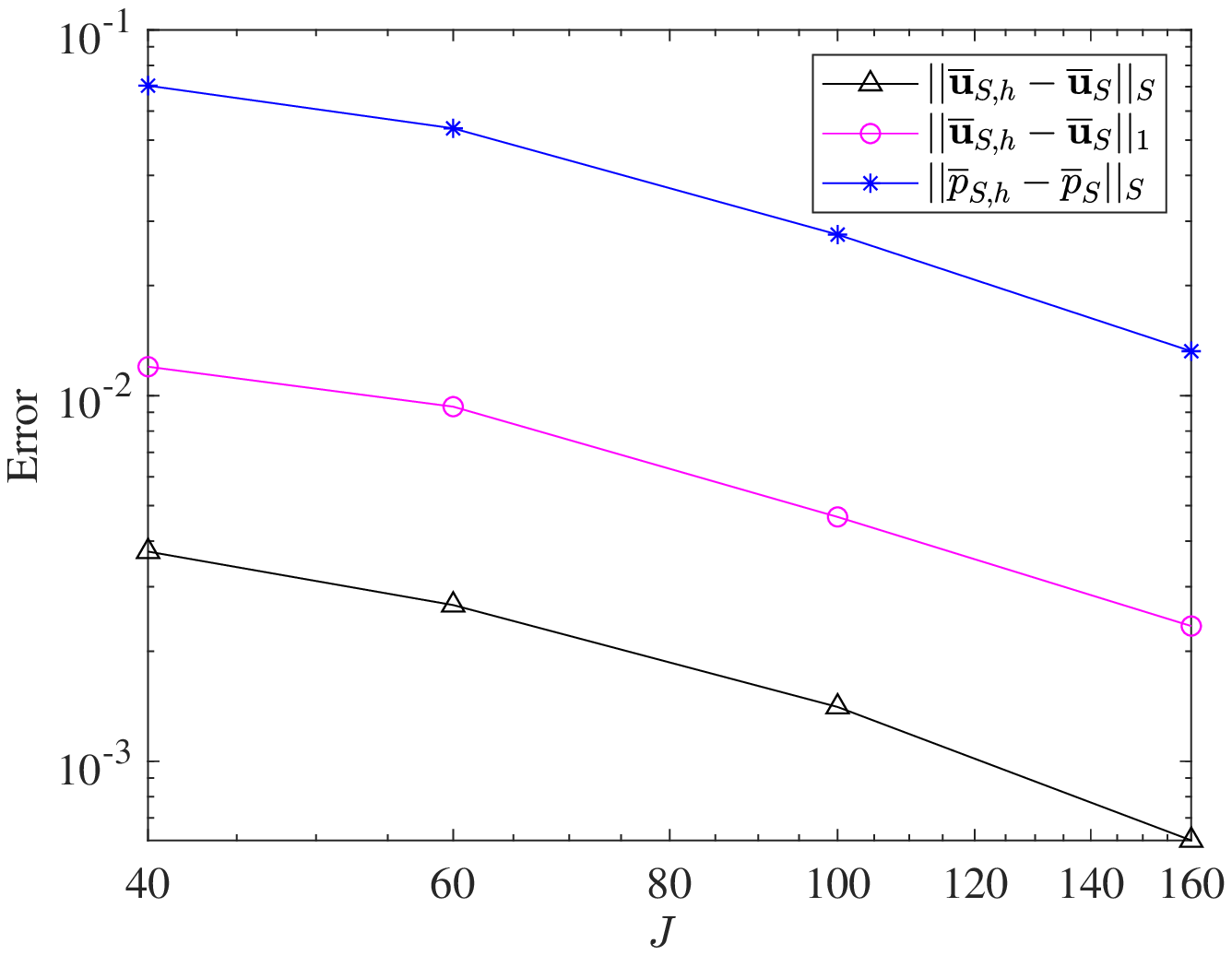} \hspace{3mm}
	\includegraphics[width=50mm,height=42mm]{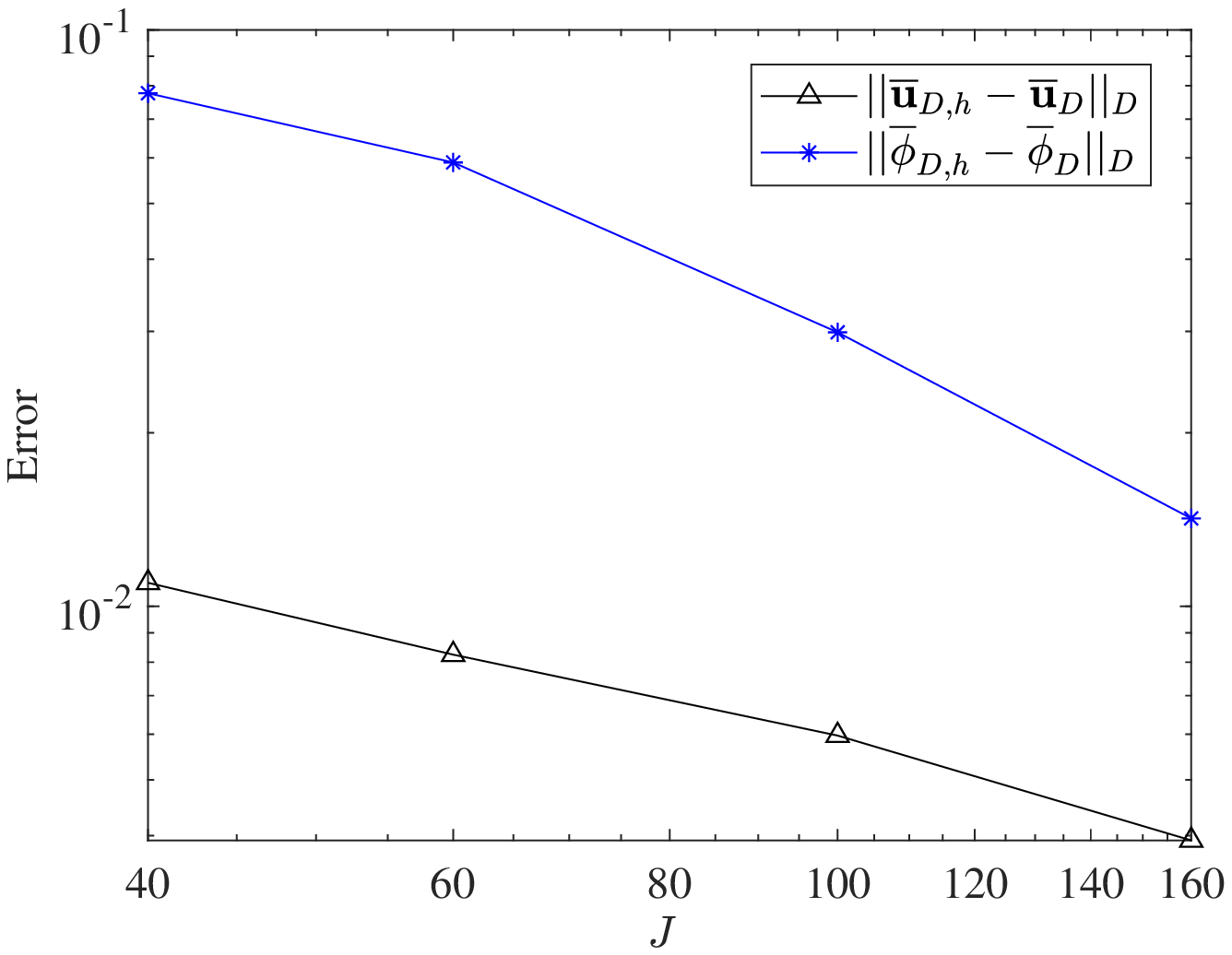}
	\caption{The numerical errors with $J=40,60,100,160$ realizations while $h=\frac{1}{32}$.}
	\label{errorMC}
\end{figure}
We first check the convergence rate of the Monte Carlo method with respect to $J$, the number of samples. Due to the unknown exact solution of this random model, we take the expectation of numerical solutions of $J_0=500$ realizations as our exact solution, for instance, we refer such expectation of  the Stokes velocity as $\overline{\mathbf{u}}_S$, and evaluate the approximation errors based on it. Here, we fix the mesh size to $h=\frac{1}{32}$. Then, the numerical errors with $J=40,60,100,160$ realizations are displayed in Fig. \ref{errorMC}. From this figure, we can see that along with the increasing $J$, the numerical approximation will be more accurate, which coincides with the classical Monte Carlo method.

To get useful statistical information from the solutions of random PDEs, we need to choose a large number of realizations. So, we compare the computational efficiency of our proposed algorithm and the traditional DDM under the selected $J=1, 10, 20, 40, 80, 160$ realizations. The elapsed CPU time of both methods is shown in Table \ref{CPUtime} for comparison, from which we can clearly see that the Ensemble DDM is meaningfully faster than the traditional DDM except for the case $J=1$.
\begin{table}[!h]
	\caption{The comparison of the elapsed CPU time while the mesh size $h=\frac{1}{32}$.}
	\label{CPUtime}\tabcolsep 0pt \vspace*{-10pt}
	\par
	\begin{center}
		\renewcommand\arraystretch{1.2}
		\def\temptablewidth{1.0\textwidth}
		{\rule{\temptablewidth}{1pt}}
		\begin{tabular*}{\temptablewidth}{@{\extracolsep{\fill}}ccccccc}
			\hline
			$J$ & 1  &  10  &   20 &  40  &   80   &160\\
			\hline
			Traditional DDM & 13.3 &  141.3  & 285.1  &  561.5 &  1120.5 & 2261.8
			\\ 
			Ensemble DDM  & 18.7 & 72.4  &  142.2   &  279.7  &  559.4 & 1140.0
			\\ 
			\hline
		\end{tabular*}%
	\end{center}
\end{table}

With the choices of $J=80$ realizations and the mesh size $h=\frac{1}{32}$, we present the velocity streamlines in Fig. \ref{3D} of the random Stokes-Darcy model. From Fig. \ref{3D}, the velocity streamlines are regular across the interface, and the magnitudes of the velocity are scattered among the random Stokes-Darcy domain reasonably for both methods. 
More importantly, both methods capture the same behaviors while the Ensemble DDM saves $50.0\%$ of the computation time.


It is of practical interest to examine the application of
 the Ensemble DDM to the realistic physical parameters. We multiply the previously defined random permeability coefficients $k_{11}(\mathbf{x}, \omega)$ and $k_{22}(\mathbf{x}, \omega)$ by $10^{-6}$ respectively, i.e $\bar{k}_{11}(\mathbf{x}, \omega)=10^{-6}k_{11}(\mathbf{x}, \omega)$ and $\bar{k}_{22}(\mathbf{x}, \omega)=10^{-6}k_{22}(\mathbf{x}, \omega)$. According to the analysis results in Section 6, especially Remark \ref{smallK11}, the $\alpha$ in this experiment also needs to be small enough up to the scale of $10^{-6}$. Moreover, the Robin parameters are selected as $\delta_S=10^6$ and $\delta_D=\frac{1}{5}\delta_S$, which have nearly same orders of $\frac{1}{\bar{k}_{11}(\mathbf{x}, \omega)}, \frac{1}{\bar{k}_{22}(\mathbf{x}, \omega)}$. The stopping criteria for the iterative process are also changed to a tolerance of $10^{-12}$ due to a smaller order of magnitude of the solutions caused by small permeability coefficients.  In Table \ref{Itstep}, we select the iteration steps at several samples, $j=1,10,20,40,60$, which supports our theoretical analysis.
\begin{table}[!h]
	\caption{The iteration steps of $j=1,10,20,40,60$ samples while $J=80$.}
	\label{Itstep}\tabcolsep 0pt \vspace*{-10pt}
	\par
	\begin{center}
		\renewcommand\arraystretch{1.2}
		\def\temptablewidth{0.8\textwidth}
		{\rule{\temptablewidth}{1.0pt}}
		\begin{tabular*}{\temptablewidth}{@{\extracolsep{\fill}}c||cccccc}
			\hline
			mesh size $h$ \hspace{2mm} & $j=1$ & $j=10$ & $j=20$ & $j=40$ & $j=60$   \\
			\hline
			$\frac{1}{8}$  & 30 & 21 &  22  &  23 &  20 
			\\ 
			$\frac{1}{16}$  & 27 & 19 & 20  &  21  &  19  
			\\ 
			$\frac{1}{32}$  &27 & 19  &  20   &   21   & 19  
			\\
			\hline
		\end{tabular*}%
	\end{center}
\end{table}

\begin{figure}[htbp]
	\centering
	\includegraphics[width=50mm,height=42mm]{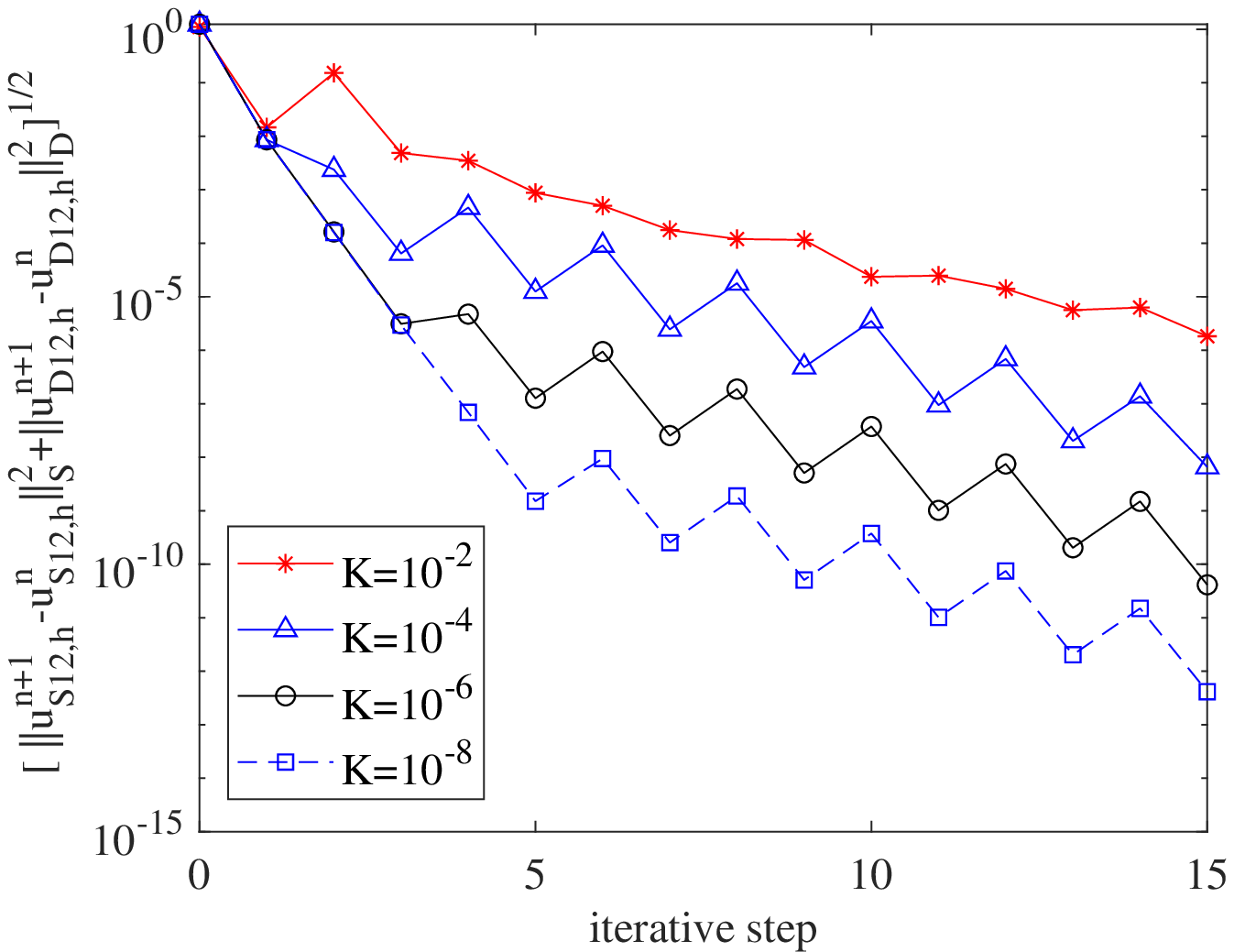} \hspace{3mm}
	\includegraphics[width=50mm,height=42mm]{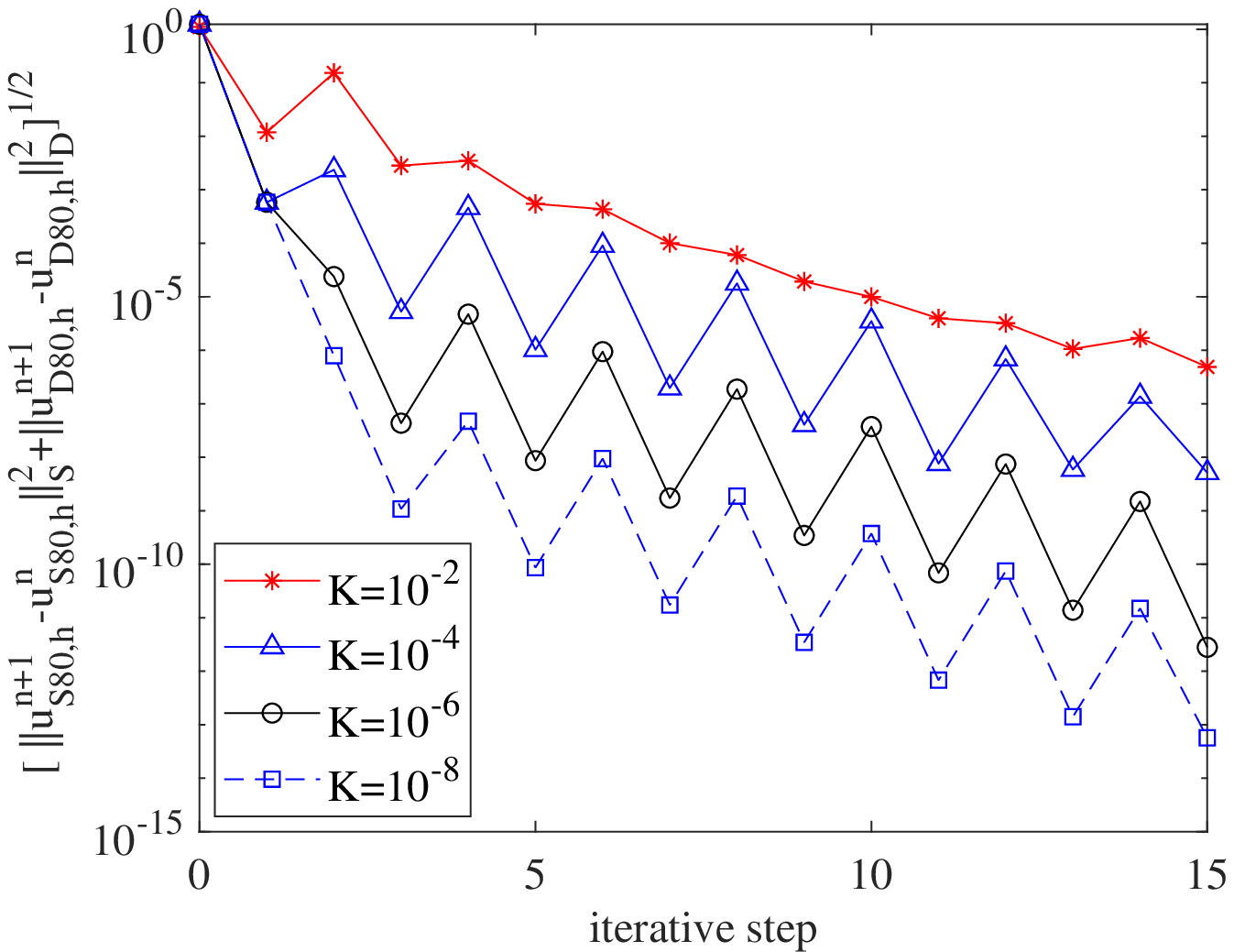}
	\caption{The $L^2$ absolute errors for the numerical velocities of both domains at $j=12$ sample (left) and $j=80$ sample  (right) while $J=80$ and $h=\frac{1}{32}$.}
	\label{DKfig}
\end{figure}

Finally, to verify that the Ensemble DDM for different sizes of permeability tensors, we test the convergence of the proposed algorithm for several cases of $\hat{k}_{ii}(\mathbf{x}, \omega)=10^{-j}k_{ii}(\mathbf{x}, \omega)$, with $i=1,2$ and $j=2,4,6,8$, in Fig. \ref{DKfig}. Here we choose the fixed stopping criteria for the iterative process and compute the $L^2$-norm of two successive iterative solutions of the Stokes and Darcy velocities. Then Fig. \ref{DKfig} clearly shows that our approach is convergent for realistic applications while $\delta_S>\delta_D$.

\section{Conclusions}
In this paper, an efficient Ensemble DDM algorithm is proposed to solve the fully-mixed random Stokes-Darcy model with BJ interface conditions. We utilize the Monte Carlo method for the coupled model with random inputs to derive some deterministic Stokes-Darcy models.  With a small disturbance of the physical parameter, the mesh-independent convergence rates are derived rigorously by choosing suitable Robin parameters. Optimized Robin parameters are derived by the optimized Schwarz method to accelerate the convergence. Moreover, we obtain the almost optimal geometric convergence for small hydraulic conductivity in practice.  This kind of idea can also be extended to the Navier-Stokes/Darcy problem, and fluid-fluid problem in parallel.

\section{Acknowledgments}
The authors would like to thank the anonymous referees for their very valuable comments and suggestions.


\begin{thebibliography}{999}


\bibitem{DMQ02} { Discacciati M,  Miglio E., Quarteroni A:} {Mathematical and numerical models for coupling surface and groundwater flows}, Appl. Numer. Math., 43 (2002) 57-74.

\bibitem{LMLayton2003} { Layton WJ, Schieweck F, Yotov I:}  {Coupling fluid flow with porous media flow}. SIAM J. Numer. Anal., 40, 2195-2218 (2003).

\bibitem{Badea} {Badea L, Discacciati M, Quarteroni A:}  {Numerical analysis of the Navier-Stokes/Darcy coupling}. Numer.
Math., 115, 195-227  (2010).


\bibitem{Cao10} {Cao Y, Gunzburger M, Hu X, Hua F, Wang X, Zhao W:}  {Finite element approximations for Stokes-Darcy flow with Beavers-Joseph interface conditions}. SIAM. J. Numer. Anal., 47, 4239-4256  (2010).













\bibitem{Zhang20} {Zhang Y, Shan L, Hou Y:} {Well-posedness and finite element approximation for the convection model in superposed fliud and porous layers}. SIAM J. Numer. Anal., 58(1), 541-564 (2020).





\bibitem{HouJY16} { Hou J, Qiu M, He X, Guo C, Wei M:}  {A dual-Porosity-Stokes model and finite element method for coupling dual-Porosity flow and free flow}. SIAM J. Sci. Comput., 38(5),  B710-B739 (2016).



















\bibitem{Discacciati07} {Discacciati M, Quarteroni A, Valli A:}  {Robin-Robin domain decomposition methods for the Stokes-Darcy coupling}. SIAM J. Numer. Anal., 45, 1246-1268 (2007).

\bibitem{Chen} {Chen W, Gunzburger M, Hua F, Wang X:}  {A parallel Robin-Robin domain decomposition method for the Stokes-Darcy system}. SIAM. J. Numer. Anal., 49, 1064-1084  (2011).




\bibitem{Cao11} {Cao Y, Gunzburger M,  He XM, Wang X:} {Robin-Robin domain decomposition methods for the steady-state Stokes-Darcy system with Beaver-Joseph interface condition}. Numer. Math., 117, 601-629  (2011).

\bibitem{He15} {He XM, Li J, Lin YP, Ming J:}  {A domain decomposition method for the steady-state Navier-Stokes-Darcy model with the Beavers-Joseph interface condition}. SIAM J. Sci. Comput., 37, S264-S290 (2015).


\bibitem{Vassilev1} {Vassilev D, Wang C, Yotov I:}  {Domain decomposition for coupled Stokes and Darcy flows}. Comput. Methods Appl. Mech. Engrg., 268, 264-283  (2014).



\bibitem{Sun2021} {Sun YZ, Sun WW, Zheng HB}. {Domain decomposition method for the fully-mixed Stokes-Darcy coupled problem}. Comput. Methods Appl. Mech. Engrg., 374, 113578 (2021).



















\bibitem{Discacciati} {Discacciati M, Gerardo-Giorda L:} {Optimized Schwarz methods for the Stokes-Darcy coupling}. IMA J. Numer. Anal., 38, 1959-1983 (2018).

\bibitem{Discacciati18} {Discacciati M, Gerardo-Giorda L:}  {Is minimizing the convergence rate a good choice for efficient Optimized Schwarz preconditioning in heterogeneous coupling the Stokes-Darcy Case}. Domain Decomposition Methods in Science and Engineering XXIV. Lecture Notes in Computational Science and Engineering (LNCSE 125), 233-41  (2018).

\bibitem{Gander} {Gander MJ, Vanzan T:} {Multilevel Optimized Schwarz Methods}. SIAM J. Sci. Comput., 42(5), A3180-A3209 (2020).

\bibitem{Liu22} {Liu YZ, Boubendir Y, He XM, He YN:} {New Optimized Robin-Robin Domain Decomposition Methods using Krylov Solvers for the Stokes-Darcy System}. SIAM J. Sci. Comput., 44(4), B1068-B1095 (2022).

\bibitem{Liu21} {Liu YZ, He YN, Li, XJ, He XM:}  {A novel convergence analysis of Robin-Robin domain decomposition method for Stokes-Darcy system with Beavers-Joseph interface condition}. Appl. Math. Letters, 119, 107181 (2021).

\bibitem{WangZ18} {Luo Y, Wang Z:} {An ensemble algorithm for numerical solutions to deterministic and random parabolic PDEs}. SIAM J. Numer. Anal., 56(2), 859-876 (2018).

\bibitem{WangZ19} {Luo Y, Wang Z:} {A multilevel Monte Carlo ensemble scheme for solving random parabolic PDEs}. SIAM J. Sci. Comput., 41, A622-A642 (2019).

\bibitem{Feng} {Feng XB, Luo Y, Vo L, Wang Z:}  {An efficient iterative method for solving parameter-dependent and random diffusion problems}. arXiv preprint arXiv, 2105.11901v1 (2021).

\bibitem{JiangN14} {Jiang N, Layton WJ:}  {An algorithm for fast calculation of flow ensembles}. Int. J. Uncertain. Quanti., 4, 273-301 (2014).

\bibitem{Max17} {Gunzburger M, Jiang N, Schneier M:}  {An ensemble-proper orthogonal decomposition method for the nonstationary Navier-Stokes equations}. SIAM J. Numer. Anal., 55, 286-304 (2017).


\bibitem{Mohebujjaman} {Mohebujjaman M, Rebholz L:} {An efficient algorithm for computation of MHD flow ensembles}. Comput. Methods Appl. Math., 17, 121-137 (2017).

\bibitem{Max19} {Gunzburger M, Jiang N, Wang Z:}  {An efficient algorithm for simulating ensembles of parameterized flow problems}. SIMA J. Numer. Anal., 39, 1180-1205 (2019).

\bibitem{JiangN19} {Jiang N, Qiu CX:}  {An efficient ensemble algorithm for the numerical approximation of stochastic Stokes-Darcy equations}. Comput. Methods Appl. Mech. Engrg., 343, 249-275 (2019).

\bibitem{JiangN21} {Jiang N, Yang HH:}  {SAV decoupled ensemble algorithms for fast computation of Stokes-Darcy flow ensembles}. Comput. Methods Appl. Mech. Engrg., 387, 114150 (2021).





\bibitem{Beavers} {Beavers G, Joseph D:}  {Boundary conditions at a naturally permeable wall}. J. Fluid Mech., 30, 197-207 (1967).


\bibitem{Jones}{Jones, IP:}  {Low Reynolds number flow past a porous spherical
shell}. Proc. Camb. Philol. Soc., 73, 231-238 (1973).

 \bibitem{Lions} {Lions PL:} {On the Schwarz alternating method III: A variant for nonoverlapping subdomains}. In: Chan TF, Glowinski R, Perianx J, Widlund OB, eds. Third International Symposium on Domain Decomposition Methods for Partial Differential Equations. Philadephia: SIAM, 202–223 (1990).

\bibitem{GR86} {Girault V,  Raviart PA:} {Finite element methods for Navier-Stokes equations}. Springer-Verlag, 1986.

\bibitem{trace} {Thom$\acute{e}$e V:}  {Galerkin finite element methods for parabolic problems}. Springer-Verlag, Berlin, second edition (2006).



\bibitem{F18}{Hecht F:} New development in freefem++. J. Numer. Math., 20, 251-265 (2012).


\end{thebibliography}
\end{document}